\documentclass[12pt,a4paper,reqno,dvipsnames]{amsart}

\usepackage{pstricks}
\usepackage{pst-bezier}
\usepackage{graphicx}
\usepackage{color}
\usepackage{multido}
\usepackage{pst-plot}
\usepackage{pst-math}
\usepackage{pstricks-add}

\makeatletter
\def\psthyperbola{\pst@object{psthyperbola}}
\def\psthyperbola@i#1#2{%
  \pst@killglue
  \begingroup
  \use@par
  \psthyperbola@ii{#1}{#2}
  \psthyperbola@ii{#1 neg}{#2}
  \endgroup
}
\def\psthyperbola@ii#1#2{
  \addto@pscode{%
    /a {#1} bind def
    /b #2 def
    /d {1 t dup mul sub} bind def }
  \parametricplot{-.99}{.99}{%
    a t dup mul 1 add d div mul
    b t 2 mul d div mul}
}
\makeatother

\newrgbcolor{slategray4}{0.423529 0.482353 0.545098}
\newrgbcolor{darkseagreen1}{0.756863 1.0 0.756863}
\newrgbcolor{peachpuff}{1.0 0.854902 0.72549}
\newrgbcolor{lemonchiffon}{1.0 0.980392 0.803922}
\newrgbcolor{honeydew}{0.941176 1.0 0.941176}
\newrgbcolor{hellrot}{0.737255 0.560784 0.560784}
\newrgbcolor{sandybrown}{0.956863 0.643137 0.376471}
\newrgbcolor{aquamarine}{0.498039 1 0.831373}
\newrgbcolor{khakifour}{0.545098 0.52549 0.305882}
\newrgbcolor{lightcoral}{0.941176 0.501961 0.501961}
\newrgbcolor{palegreen}{0.560784 0.737255 0.560784}
\newrgbcolor{grey88}{0.878431 0.878431 0.878431}
\newrgbcolor{royalblue}{0.254902 0.411765 0.882353}
\newrgbcolor{lightyellow}{1 1 0.878431}
\newrgbcolor{mfzarthimmelblau}{0.6431372549019608 0.8274509803921568 0.9333333333333333}
\newrgbcolor{mfzartstahlblau}{0.792156862745098 0.8823529411764706 1.0}
\newrgbcolor{mfzarttuerkis}{0.7333333333333333 1.0 1.0}
\newrgbcolor{mfzartazur}{0.8784313725490196 0.9333333333333333 0.9333333333333333}
\newrgbcolor{mfgedeckttuerkis}{0.25098039215686274 0.8784313725490196 0.8156862745098039}
\newrgbcolor{mfzartbraun}{1.0 0.8274509803921568 0.6078431372549019}
\newrgbcolor{mfolivgruen}{0.6352941176470588 0.803921568627451 0.35294117647058826}
\newrgbcolor{mfzartgruen}{0.7568627450980392 1.0 0.7568627450980392}
\newrgbcolor{mfkhaki}{1.0 0.9647058823529412 0.5607843137254902}
\newrgbcolor{mfzartorange}{1.0 0.8941176470588236 0.7686274509803922}
\newrgbcolor{mfzartpfirsich}{0.9333333333333333 0.8352941176470589 0.7176470588235294}
\newrgbcolor{mfhonigtau}{0.9411764705882353 1.0 0.9411764705882353}
\newrgbcolor{mfzartrosa}{1.0 0.8941176470588236 0.8823529411764706}
\newrgbcolor{mfgedeckttomate}{0.803921568627451 0.30980392156862746 0.2235294117647059}
\newrgbcolor{mflavendel}{0.9019607843137255 0.9019607843137255 0.9803921568627451}
\newrgbcolor{mfgedecktlila}{0.803921568627451 0.4117647058823529 0.788235294117647}
\newrgbcolor{mfgedecktzwetschge}{0.803921568627451 0.5882352941176471 0.803921568627451}
\newrgbcolor{mfzartweizen}{0.9333333333333333 0.8470588235294118 0.6823529411764706}
\newrgbcolor{mfzartzitrone}{1.0 0.9803921568627451 0.803921568627451}
\newrgbcolor{mfzartgold}{0.9333333333333333 0.9098039215686274 0.6666666666666666}
\newrgbcolor{mfzartgruengelb}{0.803921568627451 0.803921568627451 0.0}

\newgray{backgroundgray}{0.97}
\newgray{verylightgray}{0.9}
\newgray{ratherlightgray}{0.8}

\def\bit{\begin{itemize}}
\def\eit{\end{itemize}}
\def\beq#1{\begin{equation}\label{#1}}
\def\eeq{\end{equation}}
\def\bpf{\begin{proof}\noindent}
\def\epf{\end{proof}}
\def\ben{\begin{enumerate}}
\def\een{\end{enumerate}}
\def\bmat{\begin{bmatrix}}
\def\emat{\end{bmatrix}}

\def\BF#1{{\bf #1}}
\def\EM#1{{\em #1}}

\def\heute{\number\day.~\ifcase\month\or
 J\"anner\or Februar\or M\"arz\or April\or Mai\or Juni\or
 Juli\or August\or September\or Oktober\or November\or Dezember\fi
 \space\number\year}

\title[Hankel determinants and lattice paths]{Hankel Determinants of convoluted Catalan numbers and nonintersecting lattice paths: A bijective proof of
Cigler's Conjecture}
\author{Markus Fulmek}

\def\of#1{\left(#1\right)} 
\def\pas#1{\left(#1\right)} 
\def\setof#1{\left\{#1\right\}}

\def\abs#1{{\left\lvert{#1}\right\rvert}}

\def\strich{^\prime}

\def\defeq{\stackrel{\text{\tiny def}}{=}}

\def\defiff{{\stackrel{\text{def}}{\iff}}}

\def\fdg{\colon\;} 

\def\N{{\mathbb N}}

\def\R{{\mathbb R}}

\def\Z{{\mathbb Z}}

\def\0{{\mathbf 0}} 
\def\1{{\mathbf 1}} 

\def\Iverson#1{\left[#1\right]}

\def\signum#1{\mathop{sgn}\of{#1}}

\def\inv#1{\mathop{inv}\of{#1}}

\def\initvertex#1{\alpha\of{e}}
\def\xxxinitvertex#1{\alpha\of{e}}
\def\termvertex#1{\omega\of{e}}
\def\xxxtermvertex#1{\omega\of{e}}

\def\floor#1{\left\lfloor #1\right\rfloor}

\def\CT1{CT1}
\def\xxxCT1{CT1}

\let\phi\varphi
\let\cps\centerpspicture

\def\inv{\mathop{inv}}

\def\S{{\mathfrak S}}

\def\nilp{non\-in\-ter\-sec\-ting lat\-ti\-ce path}

\def\secA#1{\section{#1}}
\def\secB#1{\subsection{#1}}
\def\secC#1{\subsubsection{#1}}
\newtheorem{thm}{Theorem}
\newtheorem{lem}{Lemma}

\newtheorem{cor}{Corollary}

\def\cps#1{\begin{center}\input graphics/#1 \end{center}}

\def\mal{\cdot}

\def\LGV{Lind\-str\"om--Ges\-sel--Vien\-not}
\def\rhs{\BF{rhs}}
\def\lhs{\BF{lhs}}

\date{\today}
\begin{document}
\bibliographystyle{plain}
\maketitle

\begin{abstract}
In recent preprints, Cigler considered certain Hankel determinants of
convoluted Ca\-ta\-lan numbers and conjectured identities
for these determinants. In this note, we shall give a bijective proof
of Cigler's Conjecture by interpreting determinants as generating functions of
nonintersecting lattice paths: this proof employs the reflection principle,
the \LGV--method and a certain construction involving reflections and
overlays of \nilp s. Shortly after the bijective proof was presented here,
Cigler provided a shorter proof based on earlier results.
\end{abstract}

\secA{Introduction}
In recent preprints \cite{cigler:2023:HDCCN,cigler:2023:SHDOCN}, Cigler considered
Hankel determinants
\begin{equation}
\label{eq:DKMN}
D_{K,M}\of N\defeq\det\of{C_{K,i+j+M}}_{i,j=0}^{N-1},
\end{equation}
for $M\in\Z$ and $K,N\in\N$, $K\geq 1$,
where the entries are convolution powers of Catalan numbers 
for $p\in\N$
\begin{align*}
C_{K,p}&\defeq \binom{2p+K-1}p-\binom{2p+K-1}{p-1} \\
&= \frac K{p+K}\binom{2p+K-1}{p} =
\frac K{2p+K}\binom{2p+K}p,
\end{align*}
and $C_{K,p} = 0 $ for $p\in\Z$ with $p<0$.

Cigler presented the following conjecture \cite[Conjecture 1, equations (7) and (8)]{cigler:2023:HDCCN} regarding these determinants, which we formulate as a theorem
since we will give a (bijective) proof for it. This is justified even more, since shortly
after this proof was presented, Cigler \cite{cigler:2024:HDOCPOCNR}
provided an elegant (and much shorter) proof, based on earlier results of
Cigler \cite[Proposition 2.5]{cigler:2013:ASCOHD} and Andrews and Wimp \cite{Andrews-Wimp:2002:SQOPARHD}.

\begin{thm}[Cigler's Conjecture]
\label{con:cigler}
Let $m,k \in\N$, with $m,k>0$. Then for $K=2k$ we have the \EM{even identities} 
\begin{align}
D_{2k,1-k-m}\of{0} &= 1 \text{ (true by definition)}, \notag\\
D_{2k,1-k-m}\of{N} &= 0 \text{ for }N=1,2,\dots,m+k-1, \label{eq:cigler-even-simple}
\end{align}
and for all $n\in\N$
\begin{equation}
\label{eq:cigler-even}
D_{2k,1-k-m}\of{n+m+k} = \pas{-1}^{\binom{m+k}2}D_{2k,1-k+m}\of{n}
\end{equation}

Moreover, for $K=2k-1$ we have the \EM{odd identities}
\begin{align}
D_{2k-1,2-k-m}\of{0} &= 1 \text{ (true by definition)}, \notag\\
D_{2k-1,2-k-m}\of{N} &= 0 \text{ for }N=1,2,\dots,m+k-2, \label{eq:cigler-odd-simple}
\end{align}
and for all $n\in\N$
\begin{equation}
\label{eq:cigler-odd}
D_{2k-1,2-k-m}\of{n+m+k-1} = \pas{-1}^{\binom{m+k-1}2}D_{2k-1,1-k+m}\of{n}
\end{equation}
\end{thm}

We shall prove this 
by interpreting
the determinants as \EM{generating functions} of tuples of \EM{lattice paths},
for which we apply the well--known \LGV--method, i.e., we construct
a sign--reversing involution for both sides of the identities:
for the set of ``survivors'' of this involution
on the left--hand side of the identities we shall construct a \EM{second} sign--reversing
involution, and we shall show that all lattice paths which ``survive'' also this
second involution have a very special shape, from which one immediately sees a bijective correspondence with the ``survivors'' of the right--hand side of the identities. But while
this final argument is immediately obvious from a picture, the road to get there is
unfortunately long: we will therefore try to explain it through numerous pictures.

In order to keep things as brief as possible, from now on we shall refer to
\bit
\item the \EM{determinants} to the left and to the right in identities
	\eqref{eq:cigler-even-simple}, \eqref{eq:cigler-even},
	\eqref{eq:cigler-odd-simple} or \eqref{eq:cigler-odd},
\item and/or their parameters $K,M,N$, which depend on $k$, $m$, $n$,
	and on the \EM{parity} $p$ (even or odd) of the identity
	$$
	K=K\of{k,m,n,p}, M=M\of{k,m,n,p}, N=N\of{k,m,n,p}
	$$
\item and/or the lattice path interpretation for these
	determinants and their parameters (yet to be explained)
\eit
simply as 
\bit
\item \EM{right--hand side} (which we abbreviate as \rhs)
\item and \EM{left--hand side} (which we abbreviate as \lhs),
\eit
respectively.
Moreover, we shall refer to relations which are ``geometrically obvious''
in the illustrations we provide (like ``above'', ``between'', or ``to the left'')
in simple terms (instead of stating everything
with coordinates or equations).


In section~\ref{sec:preliminaries}, we briefly recall basic concepts,
and in section~\ref{sec:proof} we present the bijective proof of Cigler's Conjecture.

\secA{Preliminary considerations}
\label{sec:preliminaries}
We recall basic facts and make some elementary observations.

\secB{Binomial coefficients and lattice paths}
If we restrict the binomial coefficient $\binom{2m+k-1}m$ to the case
$m, 2m+k-1\in\N$, then we
may interpret it as the number of \EM{lattice paths} in the integer lattice $\Z^2$
consisting of 
\bit
\item horizontal steps to the right, i.e., from $\pas{x,y}$ to $\pas{x+1,y}$
\item or vertical steps upwards, i.e., from $\pas{x,y}$ to $\pas{x,y+1}$,
\eit
which start at $\pas{a-m,a-m-k+1}\in\Z^2$ and end at $\pas{a,a}\in\Z^2$ for some $a\in\Z$,
see the following picture
illustrating such path for $a=0$,  $m=3$ and $k=2$:
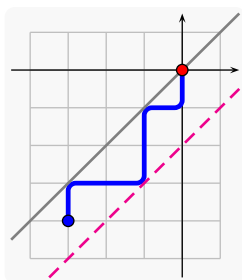
\cps{refprinc1}
If we want to count only paths which do never touch the \EM{forbidden line} $y=x-k$ (this is the
dashed line in the above picture),
then we must subtract
\bit
\item from the number $\binom{2m+k-1}m$ of \EM{all} lattice paths
\item the number of all lattice paths which \EM{do} touch the \EM{forbidden line},
\eit
and by the well--known \EM{reflection principle} \cite{andre:1887} we have
that the latter equals $\binom{2m+k-1}{m-1}$, see the following picture where
the initial segment of the path up to the first point of intersection
(indicated by a small square in the picture)
with the forbidden line is \EM{reflected} on the forbidden line:
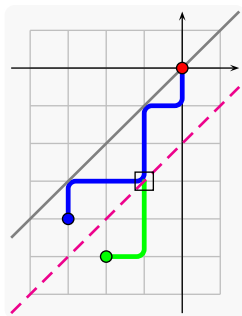
\cps{refprinc2}

\secB{Determinants interpreted as generating functions}
So, we may interpret the entry
$$
C_{K,i+j+M} = \binom{2\pas{i+j+M}+K-1}{i+j+M} - \binom{2\pas{i+j+M}+K-1}{i+j+M-1}
$$
of the determinant \eqref{eq:DKMN} as the \EM{number} of lattice paths
which
\bit
\item start at \EM{initial} point $A_i = \pas{-i-M,-i-M-K+1}$ 
\item and end at \EM{terminal} point $B_j = \pas{j,j}$ 
\item and do not touch the \EM{forbidden line} $y = x-K$,
\eit
as long as $2\pas{i+j+M}+K-1 \geq 0$ (otherwise, the number of such lattice
paths clearly is zero, while the binomial coefficients might be nonzero\footnote{For instance,
if $M=-1$, $K=1$ and $i=j=0$, then $\binom{2\pas{i+j+M}+K-1}{i+j} = \binom{-2}{0} = 1$.}).

By assigning constant weight $1$ to each single lattice path, we may interpret the
\EM{number} $C_{K,i+j+M}$ as the \EM{generating function} of the set of all such lattice paths.

Denoting the permutation group on $\setof{0,1,\dots,N-1}$ by $\S_N$,
the determinant $D_{K,M}\of N$ 
\begin{equation}
\det\of{C_{K,i+j+M}}_{i,j=0}^{N-1}
=
\sum_{\pi\in\S_N}\signum\pi\mal\prod_{i=0}^{N-1}C_{K,i+\pi\of i+M}\label{eq:det-expansion}
\end{equation}
may be viewed as the \EM{generating function} of the set $\mathcal P$ of
all $N$--tuples $t$ of lattice paths
\bit
\item from initial points
$A_i$ 
\item to terminal points $B_{\pi\of i}$,
\eit
for some $\pi\of t = \pi\in\S_N$ which we call the \EM{permutation of $t$},
where the
\EM{weight} 
of such $N$--tuple $t$ equals $\signum{\pi\of t}$:
\begin{equation}
\det\of{C_{K,i+j+M}}_{i,j=0}^{N-1}
=
\sum_{t\in\mathcal P}\signum{\pi\of t}.\label{eq:det-expansion2}
\end{equation}

The following pictures illustrate this idea by two
triples of lattice paths for $K=M=N=3$ and $\pi=\pas{1,0,2}$ (in one--line--notation):
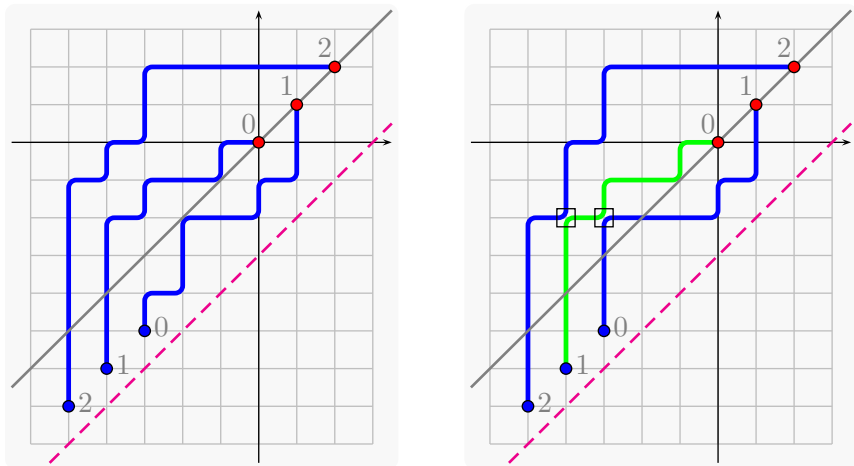
\begin{center}
\psset{unit=0.5cm}
\begin{pspicture}(-6.7,-8.7)(3.7,3.7)
\pspolygon[linecolor=white,fillstyle=solid,fillcolor=backgroundgray,linearc=0.3](-6.7,-8.7)(3.7,-8.7)(3.7,3.7)(-6.7,3.7)

\psset{linewidth=0.5pt,linecolor=gray,linestyle=solid,fillstyle=none}
\psset{linewidth=0.5pt,linecolor=gray,linestyle=solid,fillstyle=none}
\psline[linecolor=lightgray](-6,-8)(-6,3)
\psline[linecolor=lightgray](-5,-8)(-5,3)
\psline[linecolor=lightgray](-4,-8)(-4,3)
\psline[linecolor=lightgray](-3,-8)(-3,3)
\psline[linecolor=lightgray](-2,-8)(-2,3)
\psline[linecolor=lightgray](-1,-8)(-1,3)
\psline[linecolor=lightgray](0,-8)(0,3)
\psline[linecolor=lightgray](1,-8)(1,3)
\psline[linecolor=lightgray](2,-8)(2,3)
\psline[linecolor=lightgray](3,-8)(3,3)
\psline[linecolor=lightgray](-6,-8)(3,-8)
\psline[linecolor=lightgray](-6,-7)(3,-7)
\psline[linecolor=lightgray](-6,-6)(3,-6)
\psline[linecolor=lightgray](-6,-5)(3,-5)
\psline[linecolor=lightgray](-6,-4)(3,-4)
\psline[linecolor=lightgray](-6,-3)(3,-3)
\psline[linecolor=lightgray](-6,-2)(3,-2)
\psline[linecolor=lightgray](-6,-1)(3,-1)
\psline[linecolor=lightgray](-6,0)(3,0)
\psline[linecolor=lightgray](-6,1)(3,1)
\psline[linecolor=lightgray](-6,2)(3,2)
\psline[linecolor=lightgray](-6,3)(3,3)
\psline[linecolor=black,arrows=->](-6.5,0)(3.5,0)
\psline[linecolor=black,arrows=->](0,-8.5)(0,3.5)
\psset{linewidth=1pt,linecolor=black,linestyle=solid,fillstyle=none}
\psline[linewidth=0.125,linearc=0.2,arrows=c-c,linecolor=blue](-5,-7)(-5,-6)(-5,-5)(-5,-4)(-5,-3)(-5,-2)(-5,-1)(-4,-1)(-4,0)(-3,0)(-3,1)(-3,2)(-2,2)(-1,2)(0,2)(1,2)(2,2)
\psline[linewidth=0.125,linearc=0.2,arrows=c-c,linecolor=blue](-4,-6)(-4,-5)(-4,-4)(-4,-3)(-4,-2)(-3,-2)(-3,-1)(-2,-1)(-1,-1)(-1,0)(0,0)
\psline[linewidth=0.125,linearc=0.2,arrows=c-c,linecolor=blue](-3,-5)(-3,-4)(-2,-4)(-2,-3)(-2,-2)(-1,-2)(0,-2)(0,-1)(1,-1)(1,0)(1,1)
\psline[linecolor=gray](-6.5,-6.5)(3.5,3.5)
\psline[linecolor=magenta,linestyle=dashed](-5.5,-8.5)(3.5,0.5)
\pscircle[linecolor=black,linewidth=0.5pt,fillstyle=solid,fillcolor=blue](-3,-5){0.15}
\rput(-2.7,-4.9){ {\gray\small $0$}}
\pscircle[linecolor=black,linewidth=0.5pt,fillstyle=solid,fillcolor=red](0,0){0.15}
\rput(-0.4,0.5){ {\gray\small $0$}}
\pscircle[linecolor=black,linewidth=0.5pt,fillstyle=solid,fillcolor=blue](-4,-6){0.15}
\rput(-3.7,-5.9){ {\gray\small $1$}}
\pscircle[linecolor=black,linewidth=0.5pt,fillstyle=solid,fillcolor=red](1,1){0.15}
\rput(0.6,1.5){ {\gray\small $1$}}
\pscircle[linecolor=black,linewidth=0.5pt,fillstyle=solid,fillcolor=blue](-5,-7){0.15}
\rput(-4.7,-6.9){ {\gray\small $2$}}
\pscircle[linecolor=black,linewidth=0.5pt,fillstyle=solid,fillcolor=red](2,2){0.15}
\rput(1.6,2.5){ {\gray\small $2$}}
\end{pspicture}\hfil
\psset{unit=0.5cm}
\begin{pspicture}(-6.7,-8.7)(3.7,3.7)
\pspolygon[linecolor=white,fillstyle=solid,fillcolor=backgroundgray,linearc=0.3](-6.7,-8.7)(3.7,-8.7)(3.7,3.7)(-6.7,3.7)

\psset{linewidth=0.5pt,linecolor=gray,linestyle=solid,fillstyle=none}
\psset{linewidth=0.5pt,linecolor=gray,linestyle=solid,fillstyle=none}
\psline[linecolor=lightgray](-6,-8)(-6,3)
\psline[linecolor=lightgray](-5,-8)(-5,3)
\psline[linecolor=lightgray](-4,-8)(-4,3)
\psline[linecolor=lightgray](-3,-8)(-3,3)
\psline[linecolor=lightgray](-2,-8)(-2,3)
\psline[linecolor=lightgray](-1,-8)(-1,3)
\psline[linecolor=lightgray](0,-8)(0,3)
\psline[linecolor=lightgray](1,-8)(1,3)
\psline[linecolor=lightgray](2,-8)(2,3)
\psline[linecolor=lightgray](3,-8)(3,3)
\psline[linecolor=lightgray](-6,-8)(3,-8)
\psline[linecolor=lightgray](-6,-7)(3,-7)
\psline[linecolor=lightgray](-6,-6)(3,-6)
\psline[linecolor=lightgray](-6,-5)(3,-5)
\psline[linecolor=lightgray](-6,-4)(3,-4)
\psline[linecolor=lightgray](-6,-3)(3,-3)
\psline[linecolor=lightgray](-6,-2)(3,-2)
\psline[linecolor=lightgray](-6,-1)(3,-1)
\psline[linecolor=lightgray](-6,0)(3,0)
\psline[linecolor=lightgray](-6,1)(3,1)
\psline[linecolor=lightgray](-6,2)(3,2)
\psline[linecolor=lightgray](-6,3)(3,3)
\psline[linecolor=black,arrows=->](-6.5,0)(3.5,0)
\psline[linecolor=black,arrows=->](0,-8.5)(0,3.5)
\psset{linewidth=1pt,linecolor=black,linestyle=solid,fillstyle=none}
\psline[linewidth=0.125,linearc=0.2,arrows=c-c,linecolor=blue](-5,-7)(-5,-6)(-5,-5)(-5,-4)(-5,-3)(-5,-2)(-4,-2)(-4,-1)(-4,0)(-3,0)(-3,1)(-3,2)(-2,2)(-1,2)(0,2)(1,2)(2,2)
\psline[linewidth=0.125,linearc=0.2,arrows=c-c,linecolor=green](-4,-6)(-4,-5)(-4,-4)(-4,-3)(-4,-2)(-3,-2)(-3,-1)(-2,-1)(-1,-1)(-1,0)(0,0)
\psline[linewidth=0.125,linearc=0.2,arrows=c-c,linecolor=blue](-3,-5)(-3,-4)(-3,-3)(-3,-2)(-2,-2)(-1,-2)(0,-2)(0,-1)(1,-1)(1,0)(1,1)
\psline[linecolor=gray](-6.5,-6.5)(3.5,3.5)
\psline[linecolor=magenta,linestyle=dashed](-5.5,-8.5)(3.5,0.5)
\pscircle[linecolor=black,linewidth=0.5pt,fillstyle=solid,fillcolor=blue](-3,-5){0.15}
\rput(-2.7,-4.9){ {\gray\small $0$}}
\pscircle[linecolor=black,linewidth=0.5pt,fillstyle=solid,fillcolor=red](0,0){0.15}
\rput(-0.4,0.5){ {\gray\small $0$}}
\pscircle[linecolor=black,linewidth=0.5pt,fillstyle=solid,fillcolor=blue](-4,-6){0.15}
\rput(-3.7,-5.9){ {\gray\small $1$}}
\pscircle[linecolor=black,linewidth=0.5pt,fillstyle=solid,fillcolor=red](1,1){0.15}
\rput(0.6,1.5){ {\gray\small $1$}}
\pscircle[linecolor=black,linewidth=0.5pt,fillstyle=solid,fillcolor=blue](-5,-7){0.15}
\rput(-4.7,-6.9){ {\gray\small $2$}}
\pscircle[linecolor=black,linewidth=0.5pt,fillstyle=solid,fillcolor=red](2,2){0.15}
\rput(1.6,2.5){ {\gray\small $2$}}
\pspolygon[linecolor=black,linewidth=0.5pt](-4.24,-2.24)(-3.76,-2.24)(-3.76,-1.76)(-4.24,-1.76)
\pspolygon[linecolor=black,linewidth=0.5pt](-3.24,-2.24)(-2.76,-2.24)(-2.76,-1.76)(-3.24,-1.76)
\end{pspicture}\end{center}
Note that in the right picture (where different colours for the
paths are used just to make their course more visible), there are two points of
intersections of lattice paths, indicated in the picture by small squares: 
if some $N$--tuple of lattice paths \EM{has} such points of intersections, then
it is called \EM{intersecting}, otherwise it is called \EM{nonintersecting} (so the
left picture above shows a \EM{nonintersecting} triple of lattice paths).

In the following considerations, we will also encounter initial and terminal
points $A^\prime_i, B^\prime_j$ in \EM{other} positions, and we will sloppily refer to
any such configuration
of initial and terminal points as \EM{situation}. 
If we \EM{know} in some situation that
initial point $A^\prime_i$ is (or must be) connected with terminal point
$B^\prime_{\pi\of i}$
by a lattice path, we call
$\pi$ the 
permutation \EM{corresponding to this situation}. (So $\pi=\pas{1,0,2}$
is the permutation corresponding to the situation shown in the above
pictures.)

\secB{The Lindstr\"om--Gessel--Viennot method}
The well--known Lindstr\"om--Gessel--Viennot method \cite{Lindstroem:1973:OTVROIM,Gessel-Viennot:1998:DPAPP} constructs
an \EM{involution} $\phi$ on the set $\mathcal P$ of
all $N$--tuples of lattice paths, i.e.,
$$
\phi\fdg{\mathcal P}\to{\mathcal P}\text{ is bijective and } \phi^{-1} = \phi,
$$
which is \EM{sign--reversing}, i.e.,
$$
\phi\of t \neq t \implies \signum {\phi\of t} = - \signum  t.
$$
Clearly, such \EM{sign--reversing involution} $\phi$ results in a \EM{cancellation} of summands in \eqref{eq:det-expansion2},
the only ``survivors'' of which are the fixed points of $\phi$ (i.e., $N$--tuples $t$ of lattice paths with
$\phi\of t = t$).

Consider the \EM{lexicographic order on $\Z^2$}:
$$
\pas{r,s} \succ \pas{u,v} \defiff r > u \text{ or }\pas{r=u\text{ and } s>v}.
$$
The \LGV--involution $\phi$ is defined as follows: if $t$ is a  \EM{nonintersecting}
$N$--tuple of lattice paths, then $\phi\of t = t$. Otherwise
choose the \EM{maximal}
point of intersection $P$ \EM{in lexicographic order}
(in the right picture above, this is the point with coordinates $\pas{-3,-2}$):
%
Assuming that (precisely)
the two paths ending in $B_{j_1}$ and $B_{j_2}$
meet in $P$, we obtain another (again intersecting) $N$--tuple $t^\prime$ of lattice paths with
the \EM{opposite} sign
by \EM{exchanging} the \EM{terminal} path segments from $P$ up to $B_{j_1}$ and $B_{j_2}$,
respectively; and we set $\phi\of t = t^\prime$.

For instance, after applying this involution to the situation of the right picture
above, we obtain the following picture for $\pi\strich = \pas{0,1,2}$ (in one--line--notation,
again):
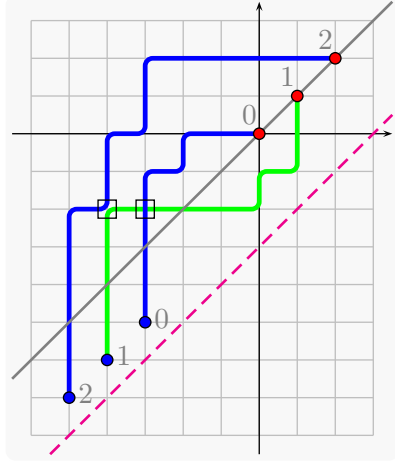
\cps{ilps2}

Clearly, the fixed points of this involution are precisely the \EM{nonintersecting}
$N$--tuples of lattice paths.
For the rest of this paper, we shall only consider such $N$--tuples of \nilp s,
which we call \EM{survivors} (since they ``survive'' the cancellation
implied by the \LGV--involution): note that the generating of
these survivors equals the determinant $D_{K,M}\of N$ in
\eqref{eq:det-expansion2}, where the \EM{sign} $\signum o$ of a
survivor $o$ connecting
initial points $A_i$ to terminal points $B_{\pi\of i}$ is defined
as  the sign $\signum\pi$ of the permutation $\pi$. 

In the typical application
of the Lindstr\"om--Gessel--Viennot method, there is only one permutation $\pi$ (typically:
the identity) for which a tuple of nonintersecting paths is possible: note that this
is \EM{not necessarily} the case in the situation we are considering here, see again the
above pictures.

\secB{Enforced segments for survivors}
Observe that the \nilp s
starting in $A_i$, $i=0,1,\dots,N-1$, 
\EM{must not} start with a horizontal step to the right
(since they must not touch the \EM{forbidden line}): so each lattice path starts with
an \EM{enforced initial segment} of vertical upward steps (which might be empty in
special cases), and since
the lattice paths are \EM{nonintersecting}, the length of the \EM{enforced initial segment}
starting in $A_{i+1}$ exceeds the length of the \EM{enforced initial segment}
starting in $A_{i}$ by \EM{two}, \EM{except} in cases where this enforced initial segment
runs into one of the end points $B_j$, where it is ``stopped'' immediately (so if
some initial point \EM{coincides} with some terminal point, the enforced initial segment
starting there has length $0$). Following these enforced initial segments from their
(original) initial points leads to (new) \EM{enforced initial points}: so instead of considering
the \EM{original} initial points, we may always consider
these \EM{enforced} initial points.

The meaning of this will become clear by looking at the 
following pictures (where $K=4$, $M=-2$ and $N=4$). Here,
the \EM{enforced initial segments} are shown as blue vertical lines in the middle picture,
and the \EM{enforced initial points} are either blue or \EM{bicoloured} (blue \EM{and} red;
if the enforced initial point \EM{coincides} with a terminal point):
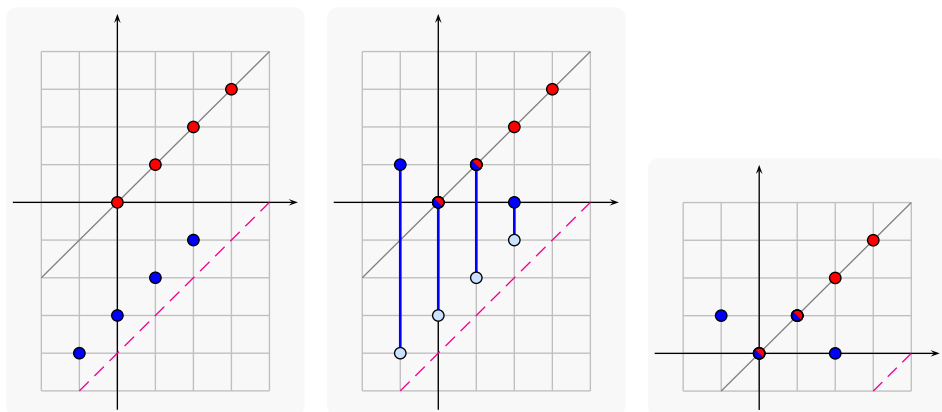
\begin{center}
\psset{unit=0.5cm}
\begin{pspicture}(-2.95,-5.7)(4.95,5.2)
\pspolygon[linecolor=white,fillstyle=solid,fillcolor=backgroundgray,linearc=0.3](-2.95,-5.7)(4.95,-5.7)(4.95,5.2)(-2.95,5.2)

\psset{linewidth=0.5pt,linecolor=gray,linestyle=solid,fillstyle=none}
\psline[linecolor=lightgray](-2,-5)(-2,4)
\psline[linecolor=lightgray](-1,-5)(-1,4)
\psline[linecolor=lightgray](0,-5)(0,4)
\psline[linecolor=lightgray](1,-5)(1,4)
\psline[linecolor=lightgray](2,-5)(2,4)
\psline[linecolor=lightgray](3,-5)(3,4)
\psline[linecolor=lightgray](4,-5)(4,4)
\psline[linecolor=lightgray](-2,-5)(4,-5)
\psline[linecolor=lightgray](-2,-4)(4,-4)
\psline[linecolor=lightgray](-2,-3)(4,-3)
\psline[linecolor=lightgray](-2,-2)(4,-2)
\psline[linecolor=lightgray](-2,-1)(4,-1)
\psline[linecolor=lightgray](-2,0)(4,0)
\psline[linecolor=lightgray](-2,1)(4,1)
\psline[linecolor=lightgray](-2,2)(4,2)
\psline[linecolor=lightgray](-2,3)(4,3)
\psline[linecolor=lightgray](-2,4)(4,4)
\psline[linecolor=black,arrows=->](-2.75,0)(4.75,0)
\psline[linecolor=black,arrows=->](0,-5.5)(0,5)
\psline(-2,-2)(4,4)
\psline[linecolor=magenta,linestyle=dashed](-1,-5)(4,0)
\psset{linewidth=1pt,linecolor=black,linestyle=solid,fillstyle=none}
\pscircle[linecolor=black,linewidth=0.5pt,fillstyle=solid,fillcolor=red](0,0){0.15}
\pscircle[linecolor=black,linewidth=0.5pt,fillstyle=solid,fillcolor=red](1,1){0.15}
\pscircle[linecolor=black,linewidth=0.5pt,fillstyle=solid,fillcolor=red](2,2){0.15}
\pscircle[linecolor=black,linewidth=0.5pt,fillstyle=solid,fillcolor=red](3,3){0.15}
\pscircle[linecolor=black,linewidth=0.5pt,fillstyle=solid,fillcolor=blue](2,-1){0.15}
\pscircle[linecolor=black,linewidth=0.5pt,fillstyle=solid,fillcolor=blue](1,-2){0.15}
\pscircle[linecolor=black,linewidth=0.5pt,fillstyle=solid,fillcolor=blue](0,-3){0.15}
\pscircle[linecolor=black,linewidth=0.5pt,fillstyle=solid,fillcolor=blue](-1,-4){0.15}
\end{pspicture}\hfil
\psset{unit=0.5cm}
\begin{pspicture}(-2.95,-5.7)(4.95,5.2)
\pspolygon[linecolor=white,fillstyle=solid,fillcolor=backgroundgray,linearc=0.3](-2.95,-5.7)(4.95,-5.7)(4.95,5.2)(-2.95,5.2)

\psset{linewidth=0.5pt,linecolor=gray,linestyle=solid,fillstyle=none}
\psline[linecolor=lightgray](-2,-5)(-2,4)
\psline[linecolor=lightgray](-1,-5)(-1,4)
\psline[linecolor=lightgray](0,-5)(0,4)
\psline[linecolor=lightgray](1,-5)(1,4)
\psline[linecolor=lightgray](2,-5)(2,4)
\psline[linecolor=lightgray](3,-5)(3,4)
\psline[linecolor=lightgray](4,-5)(4,4)
\psline[linecolor=lightgray](-2,-5)(4,-5)
\psline[linecolor=lightgray](-2,-4)(4,-4)
\psline[linecolor=lightgray](-2,-3)(4,-3)
\psline[linecolor=lightgray](-2,-2)(4,-2)
\psline[linecolor=lightgray](-2,-1)(4,-1)
\psline[linecolor=lightgray](-2,0)(4,0)
\psline[linecolor=lightgray](-2,1)(4,1)
\psline[linecolor=lightgray](-2,2)(4,2)
\psline[linecolor=lightgray](-2,3)(4,3)
\psline[linecolor=lightgray](-2,4)(4,4)
\psline[linecolor=black,arrows=->](-2.75,0)(4.75,0)
\psline[linecolor=black,arrows=->](0,-5.5)(0,5)
\psline(-2,-2)(4,4)
\psline[linecolor=magenta,linestyle=dashed](-1,-5)(4,0)
\psset{linewidth=1pt,linecolor=black,linestyle=solid,fillstyle=none}
\psline[linecolor=blue](2,-1)(2,0)
\psline[linecolor=blue](1,-2)(1,-1)(1,0)(1,1)
\psline[linecolor=blue](0,-3)(0,-2)(0,-1)(0,0)
\psline[linecolor=blue](-1,-4)(-1,-3)(-1,-2)(-1,-1)(-1,0)(-1,1)
\pscircle[linecolor=black,linewidth=0.5pt,fillstyle=solid,fillcolor=red](0,0){0.15}
\pscircle[linecolor=black,linewidth=0.5pt,fillstyle=solid,fillcolor=red](1,1){0.15}
\pscircle[linecolor=black,linewidth=0.5pt,fillstyle=solid,fillcolor=red](2,2){0.15}
\pscircle[linecolor=black,linewidth=0.5pt,fillstyle=solid,fillcolor=red](3,3){0.15}
\pscircle[linecolor=black,linewidth=0.5pt,fillstyle=solid,fillcolor=blue](2,0){0.15}
\psarc[linecolor=black,linewidth=0.5pt,fillstyle=solid,fillcolor=red](1,1){0.15}{-45.0}{135.0}
\psarc[linecolor=black,linewidth=0.5pt,fillstyle=solid,fillcolor=blue](1,1){0.15}{135.0}{315.0}
\pscircle[linecolor=black,linewidth=0.5pt](1,1){0.15}
\psarc[linecolor=black,linewidth=0.5pt,fillstyle=solid,fillcolor=red](0,0){0.15}{-45.0}{135.0}
\psarc[linecolor=black,linewidth=0.5pt,fillstyle=solid,fillcolor=blue](0,0){0.15}{135.0}{315.0}
\pscircle[linecolor=black,linewidth=0.5pt](0,0){0.15}
\pscircle[linecolor=black,linewidth=0.5pt,fillstyle=solid,fillcolor=blue](-1,1){0.15}
\pscircle[linecolor=black,linewidth=0.5pt,fillstyle=solid,fillcolor=mfzartstahlblau](2,-1){0.15}
\pscircle[linecolor=black,linewidth=0.5pt,fillstyle=solid,fillcolor=mfzartstahlblau](1,-2){0.15}
\pscircle[linecolor=black,linewidth=0.5pt,fillstyle=solid,fillcolor=mfzartstahlblau](0,-3){0.15}
\pscircle[linecolor=black,linewidth=0.5pt,fillstyle=solid,fillcolor=mfzartstahlblau](-1,-4){0.15}
\end{pspicture}\hfil
\psset{unit=0.5cm}
\begin{pspicture}(-2.95,-1.7)(4.95,5.2)
\pspolygon[linecolor=white,fillstyle=solid,fillcolor=backgroundgray,linearc=0.3](-2.95,-1.7)(4.95,-1.7)(4.95,5.2)(-2.95,5.2)

\psset{linewidth=0.5pt,linecolor=gray,linestyle=solid,fillstyle=none}
\psline[linecolor=lightgray](-2,-1)(-2,4)
\psline[linecolor=lightgray](-1,-1)(-1,4)
\psline[linecolor=lightgray](0,-1)(0,4)
\psline[linecolor=lightgray](1,-1)(1,4)
\psline[linecolor=lightgray](2,-1)(2,4)
\psline[linecolor=lightgray](3,-1)(3,4)
\psline[linecolor=lightgray](4,-1)(4,4)
\psline[linecolor=lightgray](-2,-1)(4,-1)
\psline[linecolor=lightgray](-2,0)(4,0)
\psline[linecolor=lightgray](-2,1)(4,1)
\psline[linecolor=lightgray](-2,2)(4,2)
\psline[linecolor=lightgray](-2,3)(4,3)
\psline[linecolor=lightgray](-2,4)(4,4)
\psline[linecolor=black,arrows=->](-2.75,0)(4.75,0)
\psline[linecolor=black,arrows=->](0,-1.5)(0,5)
\psline(-1,-1)(4,4)
\psline[linecolor=magenta,linestyle=dashed](3,-1)(4,0)
\psset{linewidth=1pt,linecolor=black,linestyle=solid,fillstyle=none}
\pscircle[linecolor=black,linewidth=0.5pt,fillstyle=solid,fillcolor=red](0,0){0.15}
\pscircle[linecolor=black,linewidth=0.5pt,fillstyle=solid,fillcolor=red](1,1){0.15}
\pscircle[linecolor=black,linewidth=0.5pt,fillstyle=solid,fillcolor=red](2,2){0.15}
\pscircle[linecolor=black,linewidth=0.5pt,fillstyle=solid,fillcolor=red](3,3){0.15}
\pscircle[linecolor=black,linewidth=0.5pt,fillstyle=solid,fillcolor=blue](2,0){0.15}
\psarc[linecolor=black,linewidth=0.5pt,fillstyle=solid,fillcolor=red](1,1){0.15}{-45.0}{135.0}
\psarc[linecolor=black,linewidth=0.5pt,fillstyle=solid,fillcolor=blue](1,1){0.15}{135.0}{315.0}
\pscircle[linecolor=black,linewidth=0.5pt](1,1){0.15}
\psarc[linecolor=black,linewidth=0.5pt,fillstyle=solid,fillcolor=red](0,0){0.15}{-45.0}{135.0}
\psarc[linecolor=black,linewidth=0.5pt,fillstyle=solid,fillcolor=blue](0,0){0.15}{135.0}{315.0}
\pscircle[linecolor=black,linewidth=0.5pt](0,0){0.15}
\pscircle[linecolor=black,linewidth=0.5pt,fillstyle=solid,fillcolor=blue](-1,1){0.15}
\end{pspicture}\end{center}
Note that in the situation shown in the above pictures, for any nonintersecting $4$--tuple
the lattice path
\bit
\item starting in $A_0$ or $A_3$ \EM{must} end in $B_2$ \EM{or} in $B_3$,
\item starting in $A_1$ \EM{must} end in $B_1$,
\item starting in $A_2$ \EM{must} end in $B_0$,
\eit
and the corresponding \EM{enforced sub--permutation} $$
\pi\fdg 1\mapsto 1,\;2\mapsto 0
$$
is of length $2$ is \EM{descending} and thus 
contributes $\pas{-1}^{\binom{2}2} = -1$ to the sign of the
overall permutation (which depends on whether $A_0$ is connected to $B_2$ or $B_3$;
for both possibilities there are nonintersecting $4$--tuples of lattice paths).

\secB{Positions of enforced initial points}
As a preparation for the proof, we
examine the \EM{enforced initial points} 
of the \nilp s corresponding to the \lhs\ and \rhs\ of
\bit
\item the even identities \eqref{eq:cigler-even-simple} and \eqref{eq:cigler-even},
\item and the odd identities \eqref{eq:cigler-odd-simple} and \eqref{eq:cigler-odd}.
\eit
As examples we consider
\bit
\item the even identity for $k=2$, $m=4$ and $n=2$,
\item and the odd identity for $k=3$, $m=4$ and $n=3$;
\eit
see the following pictures. As can be seen in the pictures,
there are
\bit
\item enforced initial points which are not terminal points (coloured blue),
\item terminal points which are not enforced initial points (coloured red), 
\item and enforced initial points which are \EM{also} terminal points (with
	colours blue \EM{and} red,
	we call such points \EM{two--faced}).
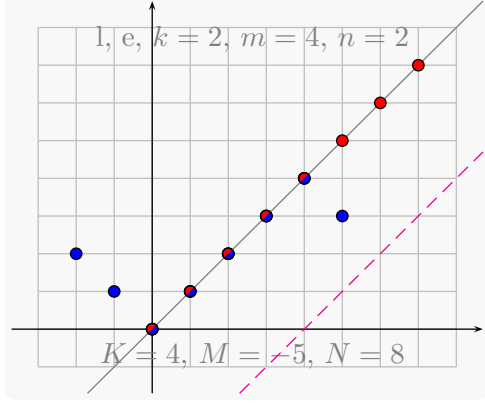
\eit
\begin{center}
\psset{unit=0.5cm}
\begin{pspicture}(-3.9,-1.9)(8.9,8.9)
\pspolygon[linecolor=white,fillstyle=solid,fillcolor=backgroundgray,linearc=0.3](-3.9,-1.9)(8.9,-1.9)(8.9,8.9)(-3.9,8.9)

\psset{linewidth=0.5pt,linecolor=gray,linestyle=solid,fillstyle=none}
\psline[linecolor=lightgray](-3,-1)(-3,8)
\psline[linecolor=lightgray](-2,-1)(-2,8)
\psline[linecolor=lightgray](-1,-1)(-1,8)
\psline[linecolor=lightgray](0,-1)(0,8)
\psline[linecolor=lightgray](1,-1)(1,8)
\psline[linecolor=lightgray](2,-1)(2,8)
\psline[linecolor=lightgray](3,-1)(3,8)
\psline[linecolor=lightgray](4,-1)(4,8)
\psline[linecolor=lightgray](5,-1)(5,8)
\psline[linecolor=lightgray](6,-1)(6,8)
\psline[linecolor=lightgray](7,-1)(7,8)
\psline[linecolor=lightgray](8,-1)(8,8)
\psline[linecolor=lightgray](-3,-1)(8,-1)
\psline[linecolor=lightgray](-3,0)(8,0)
\psline[linecolor=lightgray](-3,1)(8,1)
\psline[linecolor=lightgray](-3,2)(8,2)
\psline[linecolor=lightgray](-3,3)(8,3)
\psline[linecolor=lightgray](-3,4)(8,4)
\psline[linecolor=lightgray](-3,5)(8,5)
\psline[linecolor=lightgray](-3,6)(8,6)
\psline[linecolor=lightgray](-3,7)(8,7)
\psline[linecolor=lightgray](-3,8)(8,8)
\psline[linecolor=black,arrows=->](-3.7,0)(8.7,0)
\psline[linecolor=black,arrows=->](0,-1.7)(0,8.7)
\rput(2.5,7.7){ {\small\gray l, e, $k=2$, $m=4$, $n=2$}}
\rput(2.5,-0.7){ {\small\gray $K=4$, $M=-5$, $N=8$}}
\psset{linewidth=0.5pt,linecolor=gray,linestyle=solid,fillstyle=none}
\psline(-1.7,-1.7)(8.7,8.7)
\psline[linecolor=magenta,linestyle=dashed](2.3,-1.7)(8.7,4.7)
\pscircle[linecolor=black,linewidth=0.5pt,fillstyle=solid,fillcolor=red](0,0){0.15}
\pscircle[linecolor=black,linewidth=0.5pt,fillstyle=solid,fillcolor=red](1,1){0.15}
\pscircle[linecolor=black,linewidth=0.5pt,fillstyle=solid,fillcolor=red](2,2){0.15}
\pscircle[linecolor=black,linewidth=0.5pt,fillstyle=solid,fillcolor=red](3,3){0.15}
\pscircle[linecolor=black,linewidth=0.5pt,fillstyle=solid,fillcolor=red](4,4){0.15}
\pscircle[linecolor=black,linewidth=0.5pt,fillstyle=solid,fillcolor=red](5,5){0.15}
\pscircle[linecolor=black,linewidth=0.5pt,fillstyle=solid,fillcolor=red](6,6){0.15}
\pscircle[linecolor=black,linewidth=0.5pt,fillstyle=solid,fillcolor=red](7,7){0.15}
\psarc[linecolor=black,linewidth=0.5pt,fillstyle=solid,fillcolor=red](4,4){0.15}{45.0}{225.0}
\psarc[linecolor=black,linewidth=0.5pt,fillstyle=solid,fillcolor=blue](4,4){0.15}{225.0}{405.0}
\pscircle[linecolor=black,linewidth=0.5pt](4,4){0.15}
\psarc[linecolor=black,linewidth=0.5pt,fillstyle=solid,fillcolor=red](3,3){0.15}{45.0}{225.0}
\psarc[linecolor=black,linewidth=0.5pt,fillstyle=solid,fillcolor=blue](3,3){0.15}{225.0}{405.0}
\pscircle[linecolor=black,linewidth=0.5pt](3,3){0.15}
\psarc[linecolor=black,linewidth=0.5pt,fillstyle=solid,fillcolor=red](2,2){0.15}{45.0}{225.0}
\psarc[linecolor=black,linewidth=0.5pt,fillstyle=solid,fillcolor=blue](2,2){0.15}{225.0}{405.0}
\pscircle[linecolor=black,linewidth=0.5pt](2,2){0.15}
\psarc[linecolor=black,linewidth=0.5pt,fillstyle=solid,fillcolor=red](1,1){0.15}{45.0}{225.0}
\psarc[linecolor=black,linewidth=0.5pt,fillstyle=solid,fillcolor=blue](1,1){0.15}{225.0}{405.0}
\pscircle[linecolor=black,linewidth=0.5pt](1,1){0.15}
\psarc[linecolor=black,linewidth=0.5pt,fillstyle=solid,fillcolor=red](0,0){0.15}{45.0}{225.0}
\psarc[linecolor=black,linewidth=0.5pt,fillstyle=solid,fillcolor=blue](0,0){0.15}{225.0}{405.0}
\pscircle[linecolor=black,linewidth=0.5pt](0,0){0.15}
\pscircle[linecolor=black,linewidth=0.5pt,fillstyle=solid,fillcolor=blue](5,3){0.15}
\pscircle[linecolor=black,linewidth=0.5pt,fillstyle=solid,fillcolor=blue](-1,1){0.15}
\pscircle[linecolor=black,linewidth=0.5pt,fillstyle=solid,fillcolor=blue](-2,2){0.15}
\end{pspicture}\hfil
\psset{unit=0.5cm}
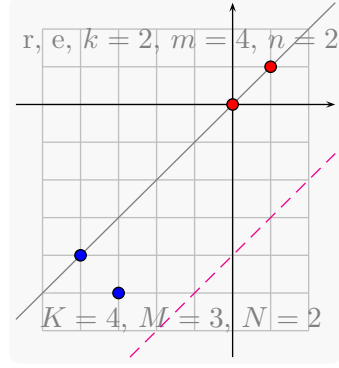
\begin{pspicture}(-5.9,-6.9)(2.9,2.9)
\pspolygon[linecolor=white,fillstyle=solid,fillcolor=backgroundgray,linearc=0.3](-5.9,-6.9)(2.9,-6.9)(2.9,2.9)(-5.9,2.9)

\psset{linewidth=0.5pt,linecolor=gray,linestyle=solid,fillstyle=none}
\psline[linecolor=lightgray](-5,-6)(-5,2)
\psline[linecolor=lightgray](-4,-6)(-4,2)
\psline[linecolor=lightgray](-3,-6)(-3,2)
\psline[linecolor=lightgray](-2,-6)(-2,2)
\psline[linecolor=lightgray](-1,-6)(-1,2)
\psline[linecolor=lightgray](0,-6)(0,2)
\psline[linecolor=lightgray](1,-6)(1,2)
\psline[linecolor=lightgray](2,-6)(2,2)
\psline[linecolor=lightgray](-5,-6)(2,-6)
\psline[linecolor=lightgray](-5,-5)(2,-5)
\psline[linecolor=lightgray](-5,-4)(2,-4)
\psline[linecolor=lightgray](-5,-3)(2,-3)
\psline[linecolor=lightgray](-5,-2)(2,-2)
\psline[linecolor=lightgray](-5,-1)(2,-1)
\psline[linecolor=lightgray](-5,0)(2,0)
\psline[linecolor=lightgray](-5,1)(2,1)
\psline[linecolor=lightgray](-5,2)(2,2)
\psline[linecolor=black,arrows=->](-5.7,0)(2.7,0)
\psline[linecolor=black,arrows=->](0,-6.7)(0,2.7)
\rput(-1.5,1.7){ {\small\gray r, e, $k=2$, $m=4$, $n=2$}}
\rput(-1.5,-5.7){ {\small\gray $K=4$, $M=3$, $N=2$}}
\psset{linewidth=0.5pt,linecolor=gray,linestyle=solid,fillstyle=none}
\psline(-5.7,-5.7)(2.7,2.7)
\psline[linecolor=magenta,linestyle=dashed](-2.7,-6.7)(2.7,-1.3)
\pscircle[linecolor=black,linewidth=0.5pt,fillstyle=solid,fillcolor=red](0,0){0.15}
\pscircle[linecolor=black,linewidth=0.5pt,fillstyle=solid,fillcolor=red](1,1){0.15}
\pscircle[linecolor=black,linewidth=0.5pt,fillstyle=solid,fillcolor=blue](-3,-5){0.15}
\pscircle[linecolor=black,linewidth=0.5pt,fillstyle=solid,fillcolor=blue](-4,-4){0.15}
\end{pspicture}\end{center}
\begin{center}
\psset{unit=0.5cm}
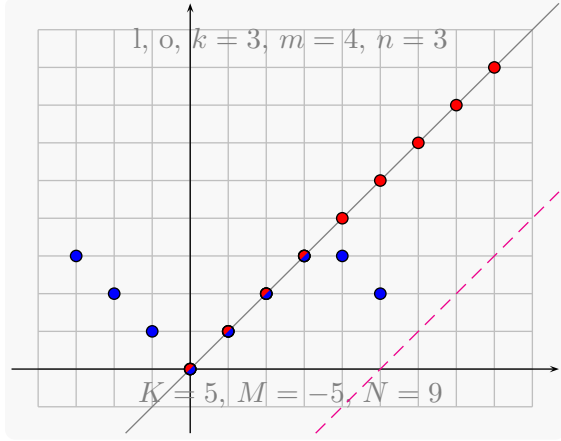
\begin{pspicture}(-4.9,-1.9)(9.9,9.9)
\pspolygon[linecolor=white,fillstyle=solid,fillcolor=backgroundgray,linearc=0.3](-4.9,-1.9)(9.9,-1.9)(9.9,9.9)(-4.9,9.9)

\psset{linewidth=0.5pt,linecolor=gray,linestyle=solid,fillstyle=none}
\psline[linecolor=lightgray](-4,-1)(-4,9)
\psline[linecolor=lightgray](-3,-1)(-3,9)
\psline[linecolor=lightgray](-2,-1)(-2,9)
\psline[linecolor=lightgray](-1,-1)(-1,9)
\psline[linecolor=lightgray](0,-1)(0,9)
\psline[linecolor=lightgray](1,-1)(1,9)
\psline[linecolor=lightgray](2,-1)(2,9)
\psline[linecolor=lightgray](3,-1)(3,9)
\psline[linecolor=lightgray](4,-1)(4,9)
\psline[linecolor=lightgray](5,-1)(5,9)
\psline[linecolor=lightgray](6,-1)(6,9)
\psline[linecolor=lightgray](7,-1)(7,9)
\psline[linecolor=lightgray](8,-1)(8,9)
\psline[linecolor=lightgray](9,-1)(9,9)
\psline[linecolor=lightgray](-4,-1)(9,-1)
\psline[linecolor=lightgray](-4,0)(9,0)
\psline[linecolor=lightgray](-4,1)(9,1)
\psline[linecolor=lightgray](-4,2)(9,2)
\psline[linecolor=lightgray](-4,3)(9,3)
\psline[linecolor=lightgray](-4,4)(9,4)
\psline[linecolor=lightgray](-4,5)(9,5)
\psline[linecolor=lightgray](-4,6)(9,6)
\psline[linecolor=lightgray](-4,7)(9,7)
\psline[linecolor=lightgray](-4,8)(9,8)
\psline[linecolor=lightgray](-4,9)(9,9)
\psline[linecolor=black,arrows=->](-4.7,0)(9.7,0)
\psline[linecolor=black,arrows=->](0,-1.7)(0,9.7)
\rput(2.5,8.7){ {\small\gray l, o, $k=3$, $m=4$, $n=3$}}
\rput(2.5,-0.7){ {\small\gray $K=5$, $M=-5$, $N=9$}}
\psset{linewidth=0.5pt,linecolor=gray,linestyle=solid,fillstyle=none}
\psline(-1.7,-1.7)(9.7,9.7)
\psline[linecolor=magenta,linestyle=dashed](3.3,-1.7)(9.7,4.7)
\pscircle[linecolor=black,linewidth=0.5pt,fillstyle=solid,fillcolor=red](0,0){0.15}
\pscircle[linecolor=black,linewidth=0.5pt,fillstyle=solid,fillcolor=red](1,1){0.15}
\pscircle[linecolor=black,linewidth=0.5pt,fillstyle=solid,fillcolor=red](2,2){0.15}
\pscircle[linecolor=black,linewidth=0.5pt,fillstyle=solid,fillcolor=red](3,3){0.15}
\pscircle[linecolor=black,linewidth=0.5pt,fillstyle=solid,fillcolor=red](4,4){0.15}
\pscircle[linecolor=black,linewidth=0.5pt,fillstyle=solid,fillcolor=red](5,5){0.15}
\pscircle[linecolor=black,linewidth=0.5pt,fillstyle=solid,fillcolor=red](6,6){0.15}
\pscircle[linecolor=black,linewidth=0.5pt,fillstyle=solid,fillcolor=red](7,7){0.15}
\pscircle[linecolor=black,linewidth=0.5pt,fillstyle=solid,fillcolor=red](8,8){0.15}
\psarc[linecolor=black,linewidth=0.5pt,fillstyle=solid,fillcolor=red](3,3){0.15}{45.0}{225.0}
\psarc[linecolor=black,linewidth=0.5pt,fillstyle=solid,fillcolor=blue](3,3){0.15}{225.0}{405.0}
\pscircle[linecolor=black,linewidth=0.5pt](3,3){0.15}
\psarc[linecolor=black,linewidth=0.5pt,fillstyle=solid,fillcolor=red](2,2){0.15}{45.0}{225.0}
\psarc[linecolor=black,linewidth=0.5pt,fillstyle=solid,fillcolor=blue](2,2){0.15}{225.0}{405.0}
\pscircle[linecolor=black,linewidth=0.5pt](2,2){0.15}
\psarc[linecolor=black,linewidth=0.5pt,fillstyle=solid,fillcolor=red](1,1){0.15}{45.0}{225.0}
\psarc[linecolor=black,linewidth=0.5pt,fillstyle=solid,fillcolor=blue](1,1){0.15}{225.0}{405.0}
\pscircle[linecolor=black,linewidth=0.5pt](1,1){0.15}
\psarc[linecolor=black,linewidth=0.5pt,fillstyle=solid,fillcolor=red](0,0){0.15}{45.0}{225.0}
\psarc[linecolor=black,linewidth=0.5pt,fillstyle=solid,fillcolor=blue](0,0){0.15}{225.0}{405.0}
\pscircle[linecolor=black,linewidth=0.5pt](0,0){0.15}
\pscircle[linecolor=black,linewidth=0.5pt,fillstyle=solid,fillcolor=blue](5,2){0.15}
\pscircle[linecolor=black,linewidth=0.5pt,fillstyle=solid,fillcolor=blue](4,3){0.15}
\pscircle[linecolor=black,linewidth=0.5pt,fillstyle=solid,fillcolor=blue](-1,1){0.15}
\pscircle[linecolor=black,linewidth=0.5pt,fillstyle=solid,fillcolor=blue](-2,2){0.15}
\pscircle[linecolor=black,linewidth=0.5pt,fillstyle=solid,fillcolor=blue](-3,3){0.15}
\end{pspicture}\hfil
\psset{unit=0.5cm}
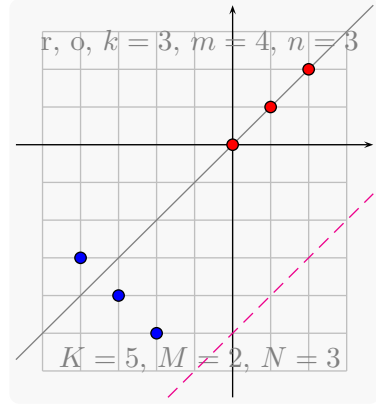
\begin{pspicture}(-5.9,-6.9)(3.9,3.9)
\pspolygon[linecolor=white,fillstyle=solid,fillcolor=backgroundgray,linearc=0.3](-5.9,-6.9)(3.9,-6.9)(3.9,3.9)(-5.9,3.9)

\psset{linewidth=0.5pt,linecolor=gray,linestyle=solid,fillstyle=none}
\psline[linecolor=lightgray](-5,-6)(-5,3)
\psline[linecolor=lightgray](-4,-6)(-4,3)
\psline[linecolor=lightgray](-3,-6)(-3,3)
\psline[linecolor=lightgray](-2,-6)(-2,3)
\psline[linecolor=lightgray](-1,-6)(-1,3)
\psline[linecolor=lightgray](0,-6)(0,3)
\psline[linecolor=lightgray](1,-6)(1,3)
\psline[linecolor=lightgray](2,-6)(2,3)
\psline[linecolor=lightgray](3,-6)(3,3)
\psline[linecolor=lightgray](-5,-6)(3,-6)
\psline[linecolor=lightgray](-5,-5)(3,-5)
\psline[linecolor=lightgray](-5,-4)(3,-4)
\psline[linecolor=lightgray](-5,-3)(3,-3)
\psline[linecolor=lightgray](-5,-2)(3,-2)
\psline[linecolor=lightgray](-5,-1)(3,-1)
\psline[linecolor=lightgray](-5,0)(3,0)
\psline[linecolor=lightgray](-5,1)(3,1)
\psline[linecolor=lightgray](-5,2)(3,2)
\psline[linecolor=lightgray](-5,3)(3,3)
\psline[linecolor=black,arrows=->](-5.7,0)(3.7,0)
\psline[linecolor=black,arrows=->](0,-6.7)(0,3.7)
\rput(-1,2.7){ {\small\gray r, o, $k=3$, $m=4$, $n=3$}}
\rput(-1,-5.7){ {\small\gray $K=5$, $M=2$, $N=3$}}
\psset{linewidth=0.5pt,linecolor=gray,linestyle=solid,fillstyle=none}
\psline(-5.7,-5.7)(3.7,3.7)
\psline[linecolor=magenta,linestyle=dashed](-1.7,-6.7)(3.7,-1.3)
\pscircle[linecolor=black,linewidth=0.5pt,fillstyle=solid,fillcolor=red](0,0){0.15}
\pscircle[linecolor=black,linewidth=0.5pt,fillstyle=solid,fillcolor=red](1,1){0.15}
\pscircle[linecolor=black,linewidth=0.5pt,fillstyle=solid,fillcolor=red](2,2){0.15}
\pscircle[linecolor=black,linewidth=0.5pt,fillstyle=solid,fillcolor=blue](-2,-5){0.15}
\pscircle[linecolor=black,linewidth=0.5pt,fillstyle=solid,fillcolor=blue](-3,-4){0.15}
\pscircle[linecolor=black,linewidth=0.5pt,fillstyle=solid,fillcolor=blue](-4,-3){0.15}
\end{pspicture}\end{center}

In order to show that there are
\bit
\item \EM{always} two--faced points in the \lhs,
\item but \EM{never} two--faced points in the \rhs,
\eit
(as the above pictures suggest),
note that the \EM{first enforced initial points} lie on the line 
$$\pas{y = -x - 2 M - K + 2}$$
with slope $-1$ through the point $\pas{-M,-M-K+2}$.
This line intersects the diagonal $\pas{y = x}$
(which contains the terminal points) in the point 
$$S=\pas{-M+1-\frac K2,-M+1-\frac K2},$$
so terminal points which are \EM{also} initial points only exist if
the following inequality holds:
\begin{equation}
\label{eq:tp=ip}
-M+1-\frac K2\geq 0 \iff \frac K2\leq -M+1.
\end{equation}
If \eqref{eq:tp=ip} holds, then the number of two--faced points
is $\floor{-M+1-\frac K2} + 1$.

\secC{Even identities}
In the even case, the parameters for the \lhs\ are
$$
K=2 k, M=1-k-m,N=n+m+k,
$$
so the point of intersection is $S=\pas{m,m}$, and inequality \eqref{eq:tp=ip} reads
$$
k\leq m+k,
$$
which is true for $m\geq 0$ (and there are $m+1$ two--faced points),
and the parameters for the \rhs\ are
$$
K=2k, M=1-k+m, N=n
$$
so the point of intersection is $S=\pas{-m,-m}$, and \eqref{eq:tp=ip} reads
$$
k\leq k-m,
$$
which is false for $m\geq 1$ (and there are no two--faced points).

\secC{Odd identities}
In the odd case, the parameters for the \lhs\ are
$$
K=2k-1, M=2-k-m,N=n+m+k-1,
$$
so the point of intersection is $S=\pas{m-\frac12,m-\frac12}$, and
\eqref{eq:tp=ip} reads
$$
k-\frac12\leq m+k-1
$$
which is true for $m\geq 1$ (and there are $m$ two--faced points),
and the parameters for the \rhs\ are
$$
K=2k-1, M=1-k+m,N=n, 
$$
so so the point of intersection is $S=\pas{-m+\frac12,-m+\frac12}$, and \eqref{eq:tp=ip} reads
$$
k-\frac12\leq -m+k,
$$
which is false for $m\geq 1$ (and there are no two--faced points).

\secC{Even and odd identities}
\label{sec:even-odd}
It is easy to see that for the \lhs\ \EM{and} the \rhs, the number of initial points \EM{strictly below} the diagonal $\pas{y=x}$ is always
$\min\of{k-1,N}$, and that for the \lhs\ and $n\geq 0$,
the number of initial points 
\bit
\item \EM{strictly above} the diagonal $\pas{y=x}$ is always $n$,
\item \EM{on} the diagonal (these are the \EM{two--faced} points) is
\bit
\item $m+1$ for the even case
\item and $m$ for the odd case,
\eit
which we express uniformly as $\pas{m+\Iverson{\text{even}}}$, using
Iverson's bracket.
\eit
(\EM{Iverson's bracket} $\Iverson{\text{``some assertion''}}$ is defined
as $1$ if ``some assertion'' is true, otherwise as $0$: we shall employ
this useful notation also in the following.)

The following pictures give a schematic illustration of the situation
(here, the two--faced points for the \lhs\ are indicated by a box
coloured blue \EM{and} red,
and all other initial and terminal points are indicated by blue and
red boxes, respectively):

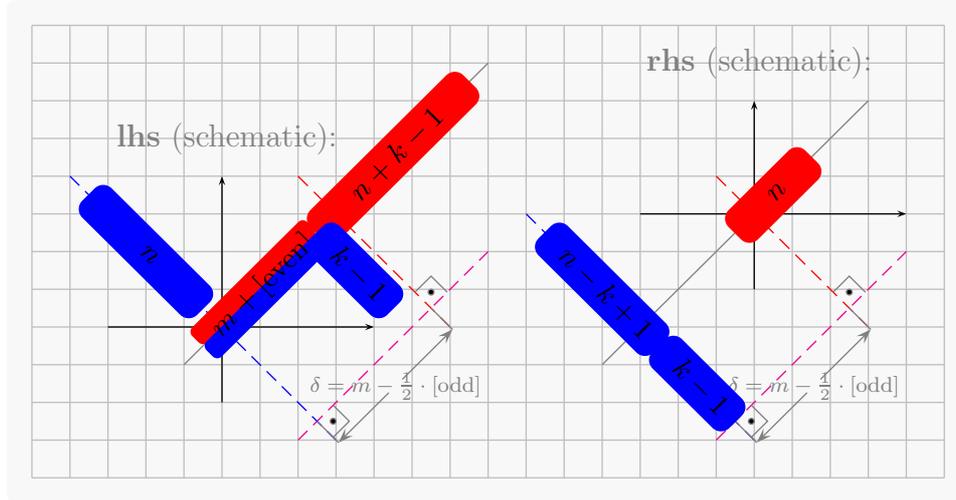
\cps{lhs_rhs}

Note that the distance $\delta$ between the orthogonal projections
\bit
\item of the first initial point (for the
	\lhs: \EM{above} the diagonal $\pas{y=x}$)
\item and of the first
	terminal point (for the \lhs: which is \EM{not two--faced})
\eit
on the forbidden line is the \EM{same} for the \lhs\ and the \rhs\ 
(namely $\delta=m+\frac12\cdot\Iverson{\text{odd}}$, see the above pictures).

\secC{A first consequence}
\label{sec:simple-equations}
These simple observations already prove identities \eqref{eq:cigler-even-simple}
and \eqref{eq:cigler-odd-simple}: if for the \lhs\
\bit
\item $N<m+k$ in the even case
\item or $N<m+k-1$ in the odd case,
\eit
then the terminal point $\pas{0,0}$ cannot be reached by \EM{any} initial point, and therefore
the generating function of the corresponding \nilp s is zero. So it only remains to prove
identities
\eqref{eq:cigler-even} and \eqref{eq:cigler-odd}

\secB{A warmup--exercise: Cigler's Conjecture for $k=1$}
The simple observations we made so far already give a simple bijective proof (essentially by looking at suitable pictures)
for
the special case $k=1$ of Theorem~\ref{con:cigler} (the \EM{odd identities} for $k=1$ are, in fact,
equivalent to Cigler's earlier result \cite[equation (5)]{cigler:2023:HDCCN},
see also \cite[Theorem 1, equation (6)]{cigler:2023:SHDOCN}):
\begin{proof}[Proof of Cigler's Conjecture for special case $k=1$]
We start with the \EM{odd identities} with parameters $K=2k-1 = 1$ and
\bit
\item $M=1-m$ and $N=n+m$ for the \lhs,
\item $M=m$ and $N=n$ for the \rhs.
\eit
Observe that in this situation
the initial points $A_i$ and the terminal points $B_j$ are
\EM{collinear}.

Recall that \eqref{eq:cigler-odd-simple} was already proved in
section~\ref {sec:simple-equations}.

For $N=m$ (this is the special case $n=0$ for the \lhs\
of \eqref{eq:cigler-odd}), we have
$$
B_i= \pas{i,i} = A_{m-1-i} \text{ for } i=0,1,\dots,m-1,
$$
and there is \EM{precisely one} $N$--tuple of \nilp s
(where all paths
have length $0$), and the corresponding permutation
$$
\pi\fdg i \mapsto m-1-i \text{ for } i=0,1,\dots,m-1
$$
has sign $\signum{\pi} = \pas{-1}^{\binom{m}2}$. Since the \rhs\
of \eqref{eq:cigler-odd} equals $1$ by convention for $N=n=0$, equality holds
in this special case.

The following pictures illustrate these simple observations for $m=3$ (i.e.,
$M=2-k-m =-2$) and $N=n=1,2,3$
(again, two--faced points 
are coloured blue \EM{and} red):
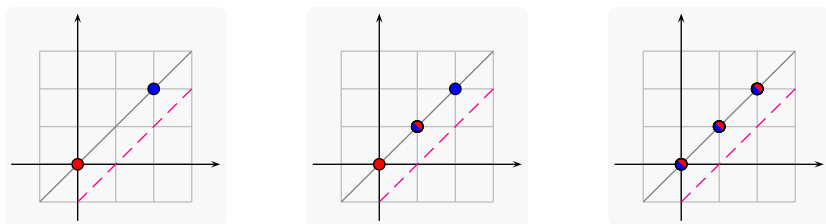
\begin{center}
\psset{unit=0.5cm}
\begin{pspicture}(-1.95,-1.7)(3.95,4.2)
\pspolygon[linecolor=white,fillstyle=solid,fillcolor=backgroundgray,linearc=0.3](-1.95,-1.7)(3.95,-1.7)(3.95,4.2)(-1.95,4.2)

\psset{linewidth=0.5pt,linecolor=gray,linestyle=solid,fillstyle=none}
\psline[linecolor=lightgray](-1,-1)(-1,3)
\psline[linecolor=lightgray](0,-1)(0,3)
\psline[linecolor=lightgray](1,-1)(1,3)
\psline[linecolor=lightgray](2,-1)(2,3)
\psline[linecolor=lightgray](3,-1)(3,3)
\psline[linecolor=lightgray](-1,-1)(3,-1)
\psline[linecolor=lightgray](-1,0)(3,0)
\psline[linecolor=lightgray](-1,1)(3,1)
\psline[linecolor=lightgray](-1,2)(3,2)
\psline[linecolor=lightgray](-1,3)(3,3)
\psline[linecolor=black,arrows=->](-1.75,0)(3.75,0)
\psline[linecolor=black,arrows=->](0,-1.5)(0,4)
\psline(-1,-1)(3,3)
\psline[linecolor=magenta,linestyle=dashed](0,-1)(3,2)
\psset{linewidth=1pt,linecolor=black,linestyle=solid,fillstyle=none}
\pscircle[linecolor=black,linewidth=0.5pt,fillstyle=solid,fillcolor=red](0,0){0.15}
\pscircle[linecolor=black,linewidth=0.5pt,fillstyle=solid,fillcolor=blue](2,2){0.15}
\end{pspicture}\hfil
\psset{unit=0.5cm}
\begin{pspicture}(-1.95,-1.7)(3.95,4.2)
\pspolygon[linecolor=white,fillstyle=solid,fillcolor=backgroundgray,linearc=0.3](-1.95,-1.7)(3.95,-1.7)(3.95,4.2)(-1.95,4.2)

\psset{linewidth=0.5pt,linecolor=gray,linestyle=solid,fillstyle=none}
\psline[linecolor=lightgray](-1,-1)(-1,3)
\psline[linecolor=lightgray](0,-1)(0,3)
\psline[linecolor=lightgray](1,-1)(1,3)
\psline[linecolor=lightgray](2,-1)(2,3)
\psline[linecolor=lightgray](3,-1)(3,3)
\psline[linecolor=lightgray](-1,-1)(3,-1)
\psline[linecolor=lightgray](-1,0)(3,0)
\psline[linecolor=lightgray](-1,1)(3,1)
\psline[linecolor=lightgray](-1,2)(3,2)
\psline[linecolor=lightgray](-1,3)(3,3)
\psline[linecolor=black,arrows=->](-1.75,0)(3.75,0)
\psline[linecolor=black,arrows=->](0,-1.5)(0,4)
\psline(-1,-1)(3,3)
\psline[linecolor=magenta,linestyle=dashed](0,-1)(3,2)
\psset{linewidth=1pt,linecolor=black,linestyle=solid,fillstyle=none}
\pscircle[linecolor=black,linewidth=0.5pt,fillstyle=solid,fillcolor=red](0,0){0.15}
\pscircle[linecolor=black,linewidth=0.5pt,fillstyle=solid,fillcolor=red](1,1){0.15}
\pscircle[linecolor=black,linewidth=0.5pt,fillstyle=solid,fillcolor=blue](2,2){0.15}
\psarc[linecolor=black,linewidth=0.5pt,fillstyle=solid,fillcolor=red](1,1){0.15}{-45.0}{135.0}
\psarc[linecolor=black,linewidth=0.5pt,fillstyle=solid,fillcolor=blue](1,1){0.15}{135.0}{315.0}
\pscircle[linecolor=black,linewidth=0.5pt](1,1){0.15}
\psarc[linecolor=black,linewidth=0.5pt,fillstyle=solid,fillcolor=red](1,1){0.15}{-45.0}{135.0}
\psarc[linecolor=black,linewidth=0.5pt,fillstyle=solid,fillcolor=blue](1,1){0.15}{135.0}{315.0}
\pscircle[linecolor=black,linewidth=0.5pt](1,1){0.15}
\end{pspicture}\hfil
\psset{unit=0.5cm}
\begin{pspicture}(-1.95,-1.7)(3.95,4.2)
\pspolygon[linecolor=white,fillstyle=solid,fillcolor=backgroundgray,linearc=0.3](-1.95,-1.7)(3.95,-1.7)(3.95,4.2)(-1.95,4.2)

\psset{linewidth=0.5pt,linecolor=gray,linestyle=solid,fillstyle=none}
\psline[linecolor=lightgray](-1,-1)(-1,3)
\psline[linecolor=lightgray](0,-1)(0,3)
\psline[linecolor=lightgray](1,-1)(1,3)
\psline[linecolor=lightgray](2,-1)(2,3)
\psline[linecolor=lightgray](3,-1)(3,3)
\psline[linecolor=lightgray](-1,-1)(3,-1)
\psline[linecolor=lightgray](-1,0)(3,0)
\psline[linecolor=lightgray](-1,1)(3,1)
\psline[linecolor=lightgray](-1,2)(3,2)
\psline[linecolor=lightgray](-1,3)(3,3)
\psline[linecolor=black,arrows=->](-1.75,0)(3.75,0)
\psline[linecolor=black,arrows=->](0,-1.5)(0,4)
\psline(-1,-1)(3,3)
\psline[linecolor=magenta,linestyle=dashed](0,-1)(3,2)
\psset{linewidth=1pt,linecolor=black,linestyle=solid,fillstyle=none}
\pscircle[linecolor=black,linewidth=0.5pt,fillstyle=solid,fillcolor=red](0,0){0.15}
\pscircle[linecolor=black,linewidth=0.5pt,fillstyle=solid,fillcolor=red](1,1){0.15}
\pscircle[linecolor=black,linewidth=0.5pt,fillstyle=solid,fillcolor=red](2,2){0.15}
\psarc[linecolor=black,linewidth=0.5pt,fillstyle=solid,fillcolor=red](2,2){0.15}{-45.0}{135.0}
\psarc[linecolor=black,linewidth=0.5pt,fillstyle=solid,fillcolor=blue](2,2){0.15}{135.0}{315.0}
\pscircle[linecolor=black,linewidth=0.5pt](2,2){0.15}
\psarc[linecolor=black,linewidth=0.5pt,fillstyle=solid,fillcolor=red](1,1){0.15}{-45.0}{135.0}
\psarc[linecolor=black,linewidth=0.5pt,fillstyle=solid,fillcolor=blue](1,1){0.15}{135.0}{315.0}
\pscircle[linecolor=black,linewidth=0.5pt](1,1){0.15}
\psarc[linecolor=black,linewidth=0.5pt,fillstyle=solid,fillcolor=red](0,0){0.15}{-45.0}{135.0}
\psarc[linecolor=black,linewidth=0.5pt,fillstyle=solid,fillcolor=blue](0,0){0.15}{135.0}{315.0}
\pscircle[linecolor=black,linewidth=0.5pt](0,0){0.15}
\psarc[linecolor=black,linewidth=0.5pt,fillstyle=solid,fillcolor=red](1,1){0.15}{-45.0}{135.0}
\psarc[linecolor=black,linewidth=0.5pt,fillstyle=solid,fillcolor=blue](1,1){0.15}{135.0}{315.0}
\pscircle[linecolor=black,linewidth=0.5pt](1,1){0.15}
\psarc[linecolor=black,linewidth=0.5pt,fillstyle=solid,fillcolor=red](0,0){0.15}{-45.0}{135.0}
\psarc[linecolor=black,linewidth=0.5pt,fillstyle=solid,fillcolor=blue](0,0){0.15}{135.0}{315.0}
\pscircle[linecolor=black,linewidth=0.5pt](0,0){0.15}
\psarc[linecolor=black,linewidth=0.5pt,fillstyle=solid,fillcolor=red](2,2){0.15}{-45.0}{135.0}
\psarc[linecolor=black,linewidth=0.5pt,fillstyle=solid,fillcolor=blue](2,2){0.15}{135.0}{315.0}
\pscircle[linecolor=black,linewidth=0.5pt](2,2){0.15}
\end{pspicture}\end{center}

In order to show \eqref{eq:cigler-odd} for $n>0$, observe that in this specific
situation there are also ``enforced \EM{terminal}
segments'' of horizontal steps to the right for the  lattice paths, since no terminal point $B_j$
can be reached by a vertical upward step. The following pictures illustrate the situation
for $m=3$ and $n=1,2$, implying
\bit
\item $M=2-k-m=-2$ and $N=4,5$ for the \lhs\ of \eqref{eq:cigler-odd}, shown in the left
	pictures below,
\item $M=1-k+m =3$ and $N=1,2$ for the \rhs\ of \eqref{eq:cigler-odd}, shown in the right
	pictures below:
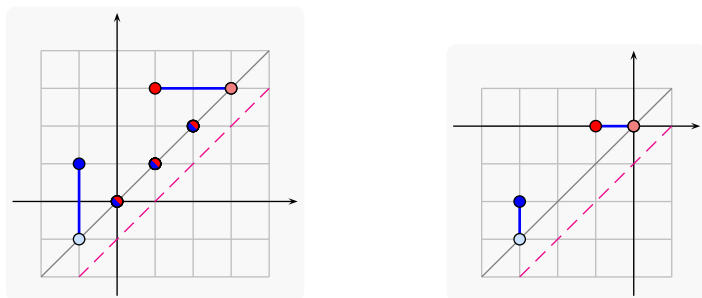
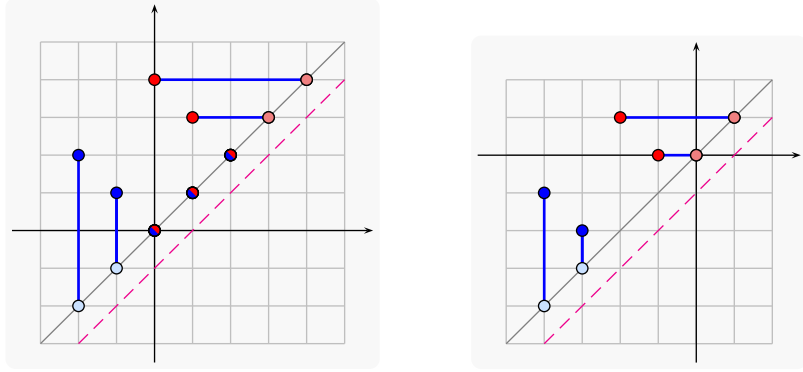
\eit
\begin{center}
\psset{unit=0.5cm}
\begin{pspicture}(-2.95,-2.7)(4.95,5.2)
\pspolygon[linecolor=white,fillstyle=solid,fillcolor=backgroundgray,linearc=0.3](-2.95,-2.7)(4.95,-2.7)(4.95,5.2)(-2.95,5.2)

\psset{linewidth=0.5pt,linecolor=gray,linestyle=solid,fillstyle=none}
\psline[linecolor=lightgray](-2,-2)(-2,4)
\psline[linecolor=lightgray](-1,-2)(-1,4)
\psline[linecolor=lightgray](0,-2)(0,4)
\psline[linecolor=lightgray](1,-2)(1,4)
\psline[linecolor=lightgray](2,-2)(2,4)
\psline[linecolor=lightgray](3,-2)(3,4)
\psline[linecolor=lightgray](4,-2)(4,4)
\psline[linecolor=lightgray](-2,-2)(4,-2)
\psline[linecolor=lightgray](-2,-1)(4,-1)
\psline[linecolor=lightgray](-2,0)(4,0)
\psline[linecolor=lightgray](-2,1)(4,1)
\psline[linecolor=lightgray](-2,2)(4,2)
\psline[linecolor=lightgray](-2,3)(4,3)
\psline[linecolor=lightgray](-2,4)(4,4)
\psline[linecolor=black,arrows=->](-2.75,0)(4.75,0)
\psline[linecolor=black,arrows=->](0,-2.5)(0,5)
\psline(-2,-2)(4,4)
\psline[linecolor=magenta,linestyle=dashed](-1,-2)(4,3)
\psset{linewidth=1pt,linecolor=black,linestyle=solid,fillstyle=none}
\psline[linearc=0.2,linecolor=blue](1,3)(2,3)(3,3)
\psline[linecolor=blue](-1,-1)(-1,0)(-1,1)
\pscircle[linecolor=black,linewidth=0.5pt,fillstyle=solid,fillcolor=red](0,0){0.15}
\pscircle[linecolor=black,linewidth=0.5pt,fillstyle=solid,fillcolor=red](1,1){0.15}
\pscircle[linecolor=black,linewidth=0.5pt,fillstyle=solid,fillcolor=red](2,2){0.15}
\pscircle[linecolor=black,linewidth=0.5pt,fillstyle=solid,fillcolor=red](3,3){0.15}
\psarc[linecolor=black,linewidth=0.5pt,fillstyle=solid,fillcolor=red](2,2){0.15}{-45.0}{135.0}
\psarc[linecolor=black,linewidth=0.5pt,fillstyle=solid,fillcolor=blue](2,2){0.15}{135.0}{315.0}
\pscircle[linecolor=black,linewidth=0.5pt](2,2){0.15}
\psarc[linecolor=black,linewidth=0.5pt,fillstyle=solid,fillcolor=red](1,1){0.15}{-45.0}{135.0}
\psarc[linecolor=black,linewidth=0.5pt,fillstyle=solid,fillcolor=blue](1,1){0.15}{135.0}{315.0}
\pscircle[linecolor=black,linewidth=0.5pt](1,1){0.15}
\psarc[linecolor=black,linewidth=0.5pt,fillstyle=solid,fillcolor=red](0,0){0.15}{-45.0}{135.0}
\psarc[linecolor=black,linewidth=0.5pt,fillstyle=solid,fillcolor=blue](0,0){0.15}{135.0}{315.0}
\pscircle[linecolor=black,linewidth=0.5pt](0,0){0.15}
\pscircle[linecolor=black,linewidth=0.5pt,fillstyle=solid,fillcolor=blue](-1,1){0.15}
\pscircle[linecolor=black,linewidth=0.5pt,fillstyle=solid,fillcolor=mfzartstahlblau](-1,-1){0.15}
\psarc[linecolor=black,linewidth=0.5pt,fillstyle=solid,fillcolor=red](1,1){0.15}{-45.0}{135.0}
\psarc[linecolor=black,linewidth=0.5pt,fillstyle=solid,fillcolor=blue](1,1){0.15}{135.0}{315.0}
\pscircle[linecolor=black,linewidth=0.5pt](1,1){0.15}
\psarc[linecolor=black,linewidth=0.5pt,fillstyle=solid,fillcolor=red](0,0){0.15}{-45.0}{135.0}
\psarc[linecolor=black,linewidth=0.5pt,fillstyle=solid,fillcolor=blue](0,0){0.15}{135.0}{315.0}
\pscircle[linecolor=black,linewidth=0.5pt](0,0){0.15}
\psarc[linecolor=black,linewidth=0.5pt,fillstyle=solid,fillcolor=red](2,2){0.15}{-45.0}{135.0}
\psarc[linecolor=black,linewidth=0.5pt,fillstyle=solid,fillcolor=blue](2,2){0.15}{135.0}{315.0}
\pscircle[linecolor=black,linewidth=0.5pt](2,2){0.15}
\pscircle[linecolor=black,linewidth=0.5pt,fillstyle=solid,fillcolor=red](1,3){0.15}
\pscircle[linecolor=black,linewidth=0.5pt,fillstyle=solid,fillcolor=lightcoral](3,3){0.15}
\end{pspicture}\hfil
\psset{unit=0.5cm}
\begin{pspicture}(-4.95,-4.7)(1.95,2.2)
\pspolygon[linecolor=white,fillstyle=solid,fillcolor=backgroundgray,linearc=0.3](-4.95,-4.7)(1.95,-4.7)(1.95,2.2)(-4.95,2.2)

\psset{linewidth=0.5pt,linecolor=gray,linestyle=solid,fillstyle=none}
\psline[linecolor=lightgray](-4,-4)(-4,1)
\psline[linecolor=lightgray](-3,-4)(-3,1)
\psline[linecolor=lightgray](-2,-4)(-2,1)
\psline[linecolor=lightgray](-1,-4)(-1,1)
\psline[linecolor=lightgray](0,-4)(0,1)
\psline[linecolor=lightgray](1,-4)(1,1)
\psline[linecolor=lightgray](-4,-4)(1,-4)
\psline[linecolor=lightgray](-4,-3)(1,-3)
\psline[linecolor=lightgray](-4,-2)(1,-2)
\psline[linecolor=lightgray](-4,-1)(1,-1)
\psline[linecolor=lightgray](-4,0)(1,0)
\psline[linecolor=lightgray](-4,1)(1,1)
\psline[linecolor=black,arrows=->](-4.75,0)(1.75,0)
\psline[linecolor=black,arrows=->](0,-4.5)(0,2)
\psline(-4,-4)(1,1)
\psline[linecolor=magenta,linestyle=dashed](-3,-4)(1,0)
\psset{linewidth=1pt,linecolor=black,linestyle=solid,fillstyle=none}
\psline[linearc=0.2,linecolor=blue](-1,0)(0,0)
\psline[linecolor=blue](-3,-3)(-3,-2)
\pscircle[linecolor=black,linewidth=0.5pt,fillstyle=solid,fillcolor=red](0,0){0.15}
\pscircle[linecolor=black,linewidth=0.5pt,fillstyle=solid,fillcolor=blue](-3,-2){0.15}
\pscircle[linecolor=black,linewidth=0.5pt,fillstyle=solid,fillcolor=mfzartstahlblau](-3,-3){0.15}
\pscircle[linecolor=black,linewidth=0.5pt,fillstyle=solid,fillcolor=red](-1,0){0.15}
\pscircle[linecolor=black,linewidth=0.5pt,fillstyle=solid,fillcolor=lightcoral](0,0){0.15}
\end{pspicture}\end{center}
\begin{center}
\psset{unit=0.5cm}
\begin{pspicture}(-3.95,-3.7)(5.95,6.2)
\pspolygon[linecolor=white,fillstyle=solid,fillcolor=backgroundgray,linearc=0.3](-3.95,-3.7)(5.95,-3.7)(5.95,6.2)(-3.95,6.2)

\psset{linewidth=0.5pt,linecolor=gray,linestyle=solid,fillstyle=none}
\psline[linecolor=lightgray](-3,-3)(-3,5)
\psline[linecolor=lightgray](-2,-3)(-2,5)
\psline[linecolor=lightgray](-1,-3)(-1,5)
\psline[linecolor=lightgray](0,-3)(0,5)
\psline[linecolor=lightgray](1,-3)(1,5)
\psline[linecolor=lightgray](2,-3)(2,5)
\psline[linecolor=lightgray](3,-3)(3,5)
\psline[linecolor=lightgray](4,-3)(4,5)
\psline[linecolor=lightgray](5,-3)(5,5)
\psline[linecolor=lightgray](-3,-3)(5,-3)
\psline[linecolor=lightgray](-3,-2)(5,-2)
\psline[linecolor=lightgray](-3,-1)(5,-1)
\psline[linecolor=lightgray](-3,0)(5,0)
\psline[linecolor=lightgray](-3,1)(5,1)
\psline[linecolor=lightgray](-3,2)(5,2)
\psline[linecolor=lightgray](-3,3)(5,3)
\psline[linecolor=lightgray](-3,4)(5,4)
\psline[linecolor=lightgray](-3,5)(5,5)
\psline[linecolor=black,arrows=->](-3.75,0)(5.75,0)
\psline[linecolor=black,arrows=->](0,-3.5)(0,6)
\psline(-3,-3)(5,5)
\psline[linecolor=magenta,linestyle=dashed](-2,-3)(5,4)
\psset{linewidth=1pt,linecolor=black,linestyle=solid,fillstyle=none}
\psline[linearc=0.2,linecolor=blue](1,3)(2,3)(3,3)
\psline[linecolor=blue](-1,-1)(-1,0)(-1,1)
\psline[linearc=0.2,linecolor=blue](0,4)(1,4)(2,4)(3,4)(4,4)
\psline[linecolor=blue](-1,-1)(-1,0)(-1,1)
\psline[linecolor=blue](-2,-2)(-2,-1)(-2,0)(-2,1)(-2,2)
\pscircle[linecolor=black,linewidth=0.5pt,fillstyle=solid,fillcolor=red](0,0){0.15}
\pscircle[linecolor=black,linewidth=0.5pt,fillstyle=solid,fillcolor=red](1,1){0.15}
\pscircle[linecolor=black,linewidth=0.5pt,fillstyle=solid,fillcolor=red](2,2){0.15}
\pscircle[linecolor=black,linewidth=0.5pt,fillstyle=solid,fillcolor=red](3,3){0.15}
\pscircle[linecolor=black,linewidth=0.5pt,fillstyle=solid,fillcolor=red](4,4){0.15}
\psarc[linecolor=black,linewidth=0.5pt,fillstyle=solid,fillcolor=red](2,2){0.15}{-45.0}{135.0}
\psarc[linecolor=black,linewidth=0.5pt,fillstyle=solid,fillcolor=blue](2,2){0.15}{135.0}{315.0}
\pscircle[linecolor=black,linewidth=0.5pt](2,2){0.15}
\psarc[linecolor=black,linewidth=0.5pt,fillstyle=solid,fillcolor=red](1,1){0.15}{-45.0}{135.0}
\psarc[linecolor=black,linewidth=0.5pt,fillstyle=solid,fillcolor=blue](1,1){0.15}{135.0}{315.0}
\pscircle[linecolor=black,linewidth=0.5pt](1,1){0.15}
\psarc[linecolor=black,linewidth=0.5pt,fillstyle=solid,fillcolor=red](0,0){0.15}{-45.0}{135.0}
\psarc[linecolor=black,linewidth=0.5pt,fillstyle=solid,fillcolor=blue](0,0){0.15}{135.0}{315.0}
\pscircle[linecolor=black,linewidth=0.5pt](0,0){0.15}
\pscircle[linecolor=black,linewidth=0.5pt,fillstyle=solid,fillcolor=blue](-1,1){0.15}
\pscircle[linecolor=black,linewidth=0.5pt,fillstyle=solid,fillcolor=blue](-2,2){0.15}
\pscircle[linecolor=black,linewidth=0.5pt,fillstyle=solid,fillcolor=mfzartstahlblau](-1,-1){0.15}
\pscircle[linecolor=black,linewidth=0.5pt,fillstyle=solid,fillcolor=mfzartstahlblau](-2,-2){0.15}
\psarc[linecolor=black,linewidth=0.5pt,fillstyle=solid,fillcolor=red](1,1){0.15}{-45.0}{135.0}
\psarc[linecolor=black,linewidth=0.5pt,fillstyle=solid,fillcolor=blue](1,1){0.15}{135.0}{315.0}
\pscircle[linecolor=black,linewidth=0.5pt](1,1){0.15}
\psarc[linecolor=black,linewidth=0.5pt,fillstyle=solid,fillcolor=red](0,0){0.15}{-45.0}{135.0}
\psarc[linecolor=black,linewidth=0.5pt,fillstyle=solid,fillcolor=blue](0,0){0.15}{135.0}{315.0}
\pscircle[linecolor=black,linewidth=0.5pt](0,0){0.15}
\psarc[linecolor=black,linewidth=0.5pt,fillstyle=solid,fillcolor=red](2,2){0.15}{-45.0}{135.0}
\psarc[linecolor=black,linewidth=0.5pt,fillstyle=solid,fillcolor=blue](2,2){0.15}{135.0}{315.0}
\pscircle[linecolor=black,linewidth=0.5pt](2,2){0.15}
\pscircle[linecolor=black,linewidth=0.5pt,fillstyle=solid,fillcolor=red](1,3){0.15}
\pscircle[linecolor=black,linewidth=0.5pt,fillstyle=solid,fillcolor=lightcoral](3,3){0.15}
\pscircle[linecolor=black,linewidth=0.5pt,fillstyle=solid,fillcolor=red](0,4){0.15}
\pscircle[linecolor=black,linewidth=0.5pt,fillstyle=solid,fillcolor=lightcoral](4,4){0.15}
\end{pspicture}\hfil
\psset{unit=0.5cm}
\begin{pspicture}(-5.95,-5.7)(2.95,3.2)
\pspolygon[linecolor=white,fillstyle=solid,fillcolor=backgroundgray,linearc=0.3](-5.95,-5.7)(2.95,-5.7)(2.95,3.2)(-5.95,3.2)

\psset{linewidth=0.5pt,linecolor=gray,linestyle=solid,fillstyle=none}
\psline[linecolor=lightgray](-5,-5)(-5,2)
\psline[linecolor=lightgray](-4,-5)(-4,2)
\psline[linecolor=lightgray](-3,-5)(-3,2)
\psline[linecolor=lightgray](-2,-5)(-2,2)
\psline[linecolor=lightgray](-1,-5)(-1,2)
\psline[linecolor=lightgray](0,-5)(0,2)
\psline[linecolor=lightgray](1,-5)(1,2)
\psline[linecolor=lightgray](2,-5)(2,2)
\psline[linecolor=lightgray](-5,-5)(2,-5)
\psline[linecolor=lightgray](-5,-4)(2,-4)
\psline[linecolor=lightgray](-5,-3)(2,-3)
\psline[linecolor=lightgray](-5,-2)(2,-2)
\psline[linecolor=lightgray](-5,-1)(2,-1)
\psline[linecolor=lightgray](-5,0)(2,0)
\psline[linecolor=lightgray](-5,1)(2,1)
\psline[linecolor=lightgray](-5,2)(2,2)
\psline[linecolor=black,arrows=->](-5.75,0)(2.75,0)
\psline[linecolor=black,arrows=->](0,-5.5)(0,3)
\psline(-5,-5)(2,2)
\psline[linecolor=magenta,linestyle=dashed](-4,-5)(2,1)
\psset{linewidth=1pt,linecolor=black,linestyle=solid,fillstyle=none}
\psline[linearc=0.2,linecolor=blue](-1,0)(0,0)
\psline[linecolor=blue](-3,-3)(-3,-2)
\psline[linearc=0.2,linecolor=blue](-2,1)(-1,1)(0,1)(1,1)
\psline[linecolor=blue](-3,-3)(-3,-2)
\psline[linecolor=blue](-4,-4)(-4,-3)(-4,-2)(-4,-1)
\pscircle[linecolor=black,linewidth=0.5pt,fillstyle=solid,fillcolor=red](0,0){0.15}
\pscircle[linecolor=black,linewidth=0.5pt,fillstyle=solid,fillcolor=red](1,1){0.15}
\pscircle[linecolor=black,linewidth=0.5pt,fillstyle=solid,fillcolor=blue](-3,-2){0.15}
\pscircle[linecolor=black,linewidth=0.5pt,fillstyle=solid,fillcolor=blue](-4,-1){0.15}
\pscircle[linecolor=black,linewidth=0.5pt,fillstyle=solid,fillcolor=mfzartstahlblau](-3,-3){0.15}
\pscircle[linecolor=black,linewidth=0.5pt,fillstyle=solid,fillcolor=mfzartstahlblau](-4,-4){0.15}
\pscircle[linecolor=black,linewidth=0.5pt,fillstyle=solid,fillcolor=red](-1,0){0.15}
\pscircle[linecolor=black,linewidth=0.5pt,fillstyle=solid,fillcolor=lightcoral](0,0){0.15}
\pscircle[linecolor=black,linewidth=0.5pt,fillstyle=solid,fillcolor=red](-2,1){0.15}
\pscircle[linecolor=black,linewidth=0.5pt,fillstyle=solid,fillcolor=lightcoral](1,1){0.15}
\end{pspicture}\end{center}
Now observe that the \EM{enforced} initial and terminal points corresponding to
the \lhs\ of \eqref{eq:cigler-odd} are simply obtained by a \EM{translation} of the
\EM{enforced} initial and terminal points corresponding to the \rhs\ (recall
the considerations and pictures in section~\ref{sec:even-odd}), and that the
$m$ \EM{enforced two--faced points}
\bit
\item correspond to the \EM{same} translation of the \EM{forbidden line}
\item and contribute a factor $\pas{-1}^{\binom{m}2}$ to the sign
of the corresponding permutation (which has fixed points $m+1,m+2,\dots,N-1$).
\eit

These observations conclude the proof of the \EM{odd identity}
\eqref{eq:cigler-odd} for the special case $k=1$.

The proof of the \EM{even identity} \eqref{eq:cigler-odd} with
parameters $K=2k = 2$ and
\bit
\item $M=-m$ and $N=n+m+1$ for the \lhs,
\item $M=m$ and $N=n$ for the \rhs,
\eit
is completely analogous to the odd case, see the following pictures
for $m=2$ (i.e., $M=1-k-m=-2$ and $M=1-k+m=2$, respectively, for the \lhs\ and
the \rhs):
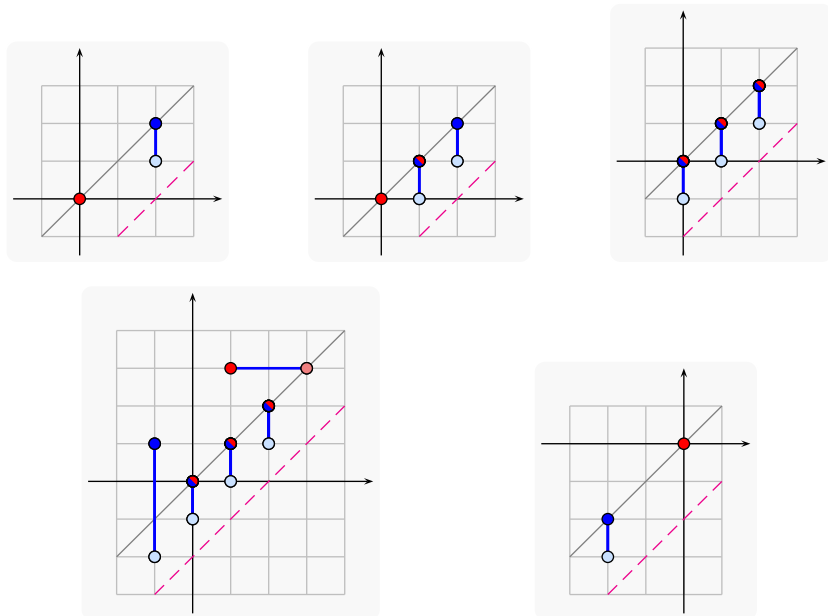
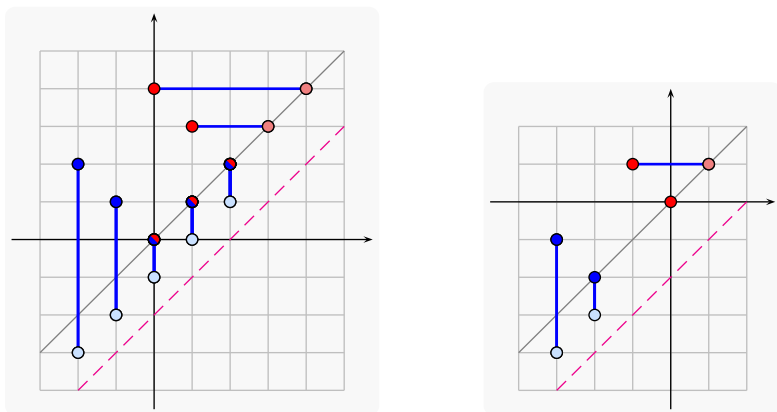
\begin{center}
\psset{unit=0.5cm}
\begin{pspicture}(-1.95,-1.7)(3.95,4.2)
\pspolygon[linecolor=white,fillstyle=solid,fillcolor=backgroundgray,linearc=0.3](-1.95,-1.7)(3.95,-1.7)(3.95,4.2)(-1.95,4.2)

\psset{linewidth=0.5pt,linecolor=gray,linestyle=solid,fillstyle=none}
\psline[linecolor=lightgray](-1,-1)(-1,3)
\psline[linecolor=lightgray](0,-1)(0,3)
\psline[linecolor=lightgray](1,-1)(1,3)
\psline[linecolor=lightgray](2,-1)(2,3)
\psline[linecolor=lightgray](3,-1)(3,3)
\psline[linecolor=lightgray](-1,-1)(3,-1)
\psline[linecolor=lightgray](-1,0)(3,0)
\psline[linecolor=lightgray](-1,1)(3,1)
\psline[linecolor=lightgray](-1,2)(3,2)
\psline[linecolor=lightgray](-1,3)(3,3)
\psline[linecolor=black,arrows=->](-1.75,0)(3.75,0)
\psline[linecolor=black,arrows=->](0,-1.5)(0,4)
\psline(-1,-1)(3,3)
\psline[linecolor=magenta,linestyle=dashed](1,-1)(3,1)
\psset{linewidth=1pt,linecolor=black,linestyle=solid,fillstyle=none}
\psline[linearc=0.2,linecolor=blue](2,1)(2,2)
\psline[linecolor=blue](2,1)(2,2)
\pscircle[linecolor=black,linewidth=0.5pt,fillstyle=solid,fillcolor=red](0,0){0.15}
\pscircle[linecolor=black,linewidth=0.5pt,fillstyle=solid,fillcolor=blue](2,2){0.15}
\pscircle[linecolor=black,linewidth=0.5pt,fillstyle=solid,fillcolor=mfzartstahlblau](2,1){0.15}
\pscircle[linecolor=black,linewidth=0.5pt,fillstyle=solid,fillcolor=mfzartstahlblau](2,1){0.15}
\pscircle[linecolor=black,linewidth=0.5pt,fillstyle=solid,fillcolor=blue](2,2){0.15}
\end{pspicture}\hfil
\psset{unit=0.5cm}
\begin{pspicture}(-1.95,-1.7)(3.95,4.2)
\pspolygon[linecolor=white,fillstyle=solid,fillcolor=backgroundgray,linearc=0.3](-1.95,-1.7)(3.95,-1.7)(3.95,4.2)(-1.95,4.2)

\psset{linewidth=0.5pt,linecolor=gray,linestyle=solid,fillstyle=none}
\psline[linecolor=lightgray](-1,-1)(-1,3)
\psline[linecolor=lightgray](0,-1)(0,3)
\psline[linecolor=lightgray](1,-1)(1,3)
\psline[linecolor=lightgray](2,-1)(2,3)
\psline[linecolor=lightgray](3,-1)(3,3)
\psline[linecolor=lightgray](-1,-1)(3,-1)
\psline[linecolor=lightgray](-1,0)(3,0)
\psline[linecolor=lightgray](-1,1)(3,1)
\psline[linecolor=lightgray](-1,2)(3,2)
\psline[linecolor=lightgray](-1,3)(3,3)
\psline[linecolor=black,arrows=->](-1.75,0)(3.75,0)
\psline[linecolor=black,arrows=->](0,-1.5)(0,4)
\psline(-1,-1)(3,3)
\psline[linecolor=magenta,linestyle=dashed](1,-1)(3,1)
\psset{linewidth=1pt,linecolor=black,linestyle=solid,fillstyle=none}
\psline[linearc=0.2,linecolor=blue](2,1)(2,2)
\psline[linecolor=blue](2,1)(2,2)
\psline[linearc=0.2,linecolor=blue](1,0)(1,1)
\psline[linecolor=blue](2,1)(2,2)
\psline[linecolor=blue](1,0)(1,1)
\pscircle[linecolor=black,linewidth=0.5pt,fillstyle=solid,fillcolor=red](0,0){0.15}
\pscircle[linecolor=black,linewidth=0.5pt,fillstyle=solid,fillcolor=red](1,1){0.15}
\pscircle[linecolor=black,linewidth=0.5pt,fillstyle=solid,fillcolor=blue](2,2){0.15}
\psarc[linecolor=black,linewidth=0.5pt,fillstyle=solid,fillcolor=red](1,1){0.15}{-45.0}{135.0}
\psarc[linecolor=black,linewidth=0.5pt,fillstyle=solid,fillcolor=blue](1,1){0.15}{135.0}{315.0}
\pscircle[linecolor=black,linewidth=0.5pt](1,1){0.15}
\pscircle[linecolor=black,linewidth=0.5pt,fillstyle=solid,fillcolor=mfzartstahlblau](2,1){0.15}
\pscircle[linecolor=black,linewidth=0.5pt,fillstyle=solid,fillcolor=mfzartstahlblau](1,0){0.15}
\pscircle[linecolor=black,linewidth=0.5pt,fillstyle=solid,fillcolor=mfzartstahlblau](2,1){0.15}
\pscircle[linecolor=black,linewidth=0.5pt,fillstyle=solid,fillcolor=blue](2,2){0.15}
\pscircle[linecolor=black,linewidth=0.5pt,fillstyle=solid,fillcolor=mfzartstahlblau](1,0){0.15}
\psarc[linecolor=black,linewidth=0.5pt,fillstyle=solid,fillcolor=red](1,1){0.15}{-45.0}{135.0}
\psarc[linecolor=black,linewidth=0.5pt,fillstyle=solid,fillcolor=blue](1,1){0.15}{135.0}{315.0}
\pscircle[linecolor=black,linewidth=0.5pt](1,1){0.15}
\end{pspicture}\hfil
\psset{unit=0.5cm}
\begin{pspicture}(-1.95,-2.7)(3.95,4.2)
\pspolygon[linecolor=white,fillstyle=solid,fillcolor=backgroundgray,linearc=0.3](-1.95,-2.7)(3.95,-2.7)(3.95,4.2)(-1.95,4.2)

\psset{linewidth=0.5pt,linecolor=gray,linestyle=solid,fillstyle=none}
\psline[linecolor=lightgray](-1,-2)(-1,3)
\psline[linecolor=lightgray](0,-2)(0,3)
\psline[linecolor=lightgray](1,-2)(1,3)
\psline[linecolor=lightgray](2,-2)(2,3)
\psline[linecolor=lightgray](3,-2)(3,3)
\psline[linecolor=lightgray](-1,-2)(3,-2)
\psline[linecolor=lightgray](-1,-1)(3,-1)
\psline[linecolor=lightgray](-1,0)(3,0)
\psline[linecolor=lightgray](-1,1)(3,1)
\psline[linecolor=lightgray](-1,2)(3,2)
\psline[linecolor=lightgray](-1,3)(3,3)
\psline[linecolor=black,arrows=->](-1.75,0)(3.75,0)
\psline[linecolor=black,arrows=->](0,-2.5)(0,4)
\psline(-1,-1)(3,3)
\psline[linecolor=magenta,linestyle=dashed](0,-2)(3,1)
\psset{linewidth=1pt,linecolor=black,linestyle=solid,fillstyle=none}
\psline[linearc=0.2,linecolor=blue](2,1)(2,2)
\psline[linecolor=blue](2,1)(2,2)
\psline[linearc=0.2,linecolor=blue](1,0)(1,1)
\psline[linecolor=blue](2,1)(2,2)
\psline[linecolor=blue](1,0)(1,1)
\psline[linearc=0.2,linecolor=blue](0,-1)(0,0)
\psline[linecolor=blue](2,1)(2,2)
\psline[linecolor=blue](1,0)(1,1)
\psline[linecolor=blue](0,-1)(0,0)
\pscircle[linecolor=black,linewidth=0.5pt,fillstyle=solid,fillcolor=red](0,0){0.15}
\pscircle[linecolor=black,linewidth=0.5pt,fillstyle=solid,fillcolor=red](1,1){0.15}
\pscircle[linecolor=black,linewidth=0.5pt,fillstyle=solid,fillcolor=red](2,2){0.15}
\psarc[linecolor=black,linewidth=0.5pt,fillstyle=solid,fillcolor=red](2,2){0.15}{-45.0}{135.0}
\psarc[linecolor=black,linewidth=0.5pt,fillstyle=solid,fillcolor=blue](2,2){0.15}{135.0}{315.0}
\pscircle[linecolor=black,linewidth=0.5pt](2,2){0.15}
\psarc[linecolor=black,linewidth=0.5pt,fillstyle=solid,fillcolor=red](1,1){0.15}{-45.0}{135.0}
\psarc[linecolor=black,linewidth=0.5pt,fillstyle=solid,fillcolor=blue](1,1){0.15}{135.0}{315.0}
\pscircle[linecolor=black,linewidth=0.5pt](1,1){0.15}
\psarc[linecolor=black,linewidth=0.5pt,fillstyle=solid,fillcolor=red](0,0){0.15}{-45.0}{135.0}
\psarc[linecolor=black,linewidth=0.5pt,fillstyle=solid,fillcolor=blue](0,0){0.15}{135.0}{315.0}
\pscircle[linecolor=black,linewidth=0.5pt](0,0){0.15}
\pscircle[linecolor=black,linewidth=0.5pt,fillstyle=solid,fillcolor=mfzartstahlblau](2,1){0.15}
\pscircle[linecolor=black,linewidth=0.5pt,fillstyle=solid,fillcolor=mfzartstahlblau](1,0){0.15}
\pscircle[linecolor=black,linewidth=0.5pt,fillstyle=solid,fillcolor=mfzartstahlblau](0,-1){0.15}
\pscircle[linecolor=black,linewidth=0.5pt,fillstyle=solid,fillcolor=mfzartstahlblau](2,1){0.15}
\pscircle[linecolor=black,linewidth=0.5pt,fillstyle=solid,fillcolor=blue](2,2){0.15}
\pscircle[linecolor=black,linewidth=0.5pt,fillstyle=solid,fillcolor=mfzartstahlblau](1,0){0.15}
\psarc[linecolor=black,linewidth=0.5pt,fillstyle=solid,fillcolor=red](1,1){0.15}{-45.0}{135.0}
\psarc[linecolor=black,linewidth=0.5pt,fillstyle=solid,fillcolor=blue](1,1){0.15}{135.0}{315.0}
\pscircle[linecolor=black,linewidth=0.5pt](1,1){0.15}
\pscircle[linecolor=black,linewidth=0.5pt,fillstyle=solid,fillcolor=mfzartstahlblau](0,-1){0.15}
\psarc[linecolor=black,linewidth=0.5pt,fillstyle=solid,fillcolor=red](0,0){0.15}{-45.0}{135.0}
\psarc[linecolor=black,linewidth=0.5pt,fillstyle=solid,fillcolor=blue](0,0){0.15}{135.0}{315.0}
\pscircle[linecolor=black,linewidth=0.5pt](0,0){0.15}
\psarc[linecolor=black,linewidth=0.5pt,fillstyle=solid,fillcolor=red](2,2){0.15}{-45.0}{135.0}
\psarc[linecolor=black,linewidth=0.5pt,fillstyle=solid,fillcolor=blue](2,2){0.15}{135.0}{315.0}
\pscircle[linecolor=black,linewidth=0.5pt](2,2){0.15}
\end{pspicture}\end{center}
\begin{center}
\psset{unit=0.5cm}
\begin{pspicture}(-2.95,-3.7)(4.95,5.2)
\pspolygon[linecolor=white,fillstyle=solid,fillcolor=backgroundgray,linearc=0.3](-2.95,-3.7)(4.95,-3.7)(4.95,5.2)(-2.95,5.2)

\psset{linewidth=0.5pt,linecolor=gray,linestyle=solid,fillstyle=none}
\psline[linecolor=lightgray](-2,-3)(-2,4)
\psline[linecolor=lightgray](-1,-3)(-1,4)
\psline[linecolor=lightgray](0,-3)(0,4)
\psline[linecolor=lightgray](1,-3)(1,4)
\psline[linecolor=lightgray](2,-3)(2,4)
\psline[linecolor=lightgray](3,-3)(3,4)
\psline[linecolor=lightgray](4,-3)(4,4)
\psline[linecolor=lightgray](-2,-3)(4,-3)
\psline[linecolor=lightgray](-2,-2)(4,-2)
\psline[linecolor=lightgray](-2,-1)(4,-1)
\psline[linecolor=lightgray](-2,0)(4,0)
\psline[linecolor=lightgray](-2,1)(4,1)
\psline[linecolor=lightgray](-2,2)(4,2)
\psline[linecolor=lightgray](-2,3)(4,3)
\psline[linecolor=lightgray](-2,4)(4,4)
\psline[linecolor=black,arrows=->](-2.75,0)(4.75,0)
\psline[linecolor=black,arrows=->](0,-3.5)(0,5)
\psline(-2,-2)(4,4)
\psline[linecolor=magenta,linestyle=dashed](-1,-3)(4,2)
\psset{linewidth=1pt,linecolor=black,linestyle=solid,fillstyle=none}
\psline[linearc=0.2,linecolor=blue](2,1)(2,2)
\psline[linecolor=blue](2,1)(2,2)
\psline[linearc=0.2,linecolor=blue](1,0)(1,1)
\psline[linecolor=blue](2,1)(2,2)
\psline[linecolor=blue](1,0)(1,1)
\psline[linearc=0.2,linecolor=blue](0,-1)(0,0)
\psline[linecolor=blue](2,1)(2,2)
\psline[linecolor=blue](1,0)(1,1)
\psline[linecolor=blue](0,-1)(0,0)
\psline[linearc=0.2,linecolor=blue](-1,-2)(-1,-1)(-1,0)(-1,1)
\psline[linearc=0.2,linecolor=blue](1,3)(2,3)(3,3)
\psline[linecolor=blue](2,1)(2,2)
\psline[linecolor=blue](1,0)(1,1)
\psline[linecolor=blue](0,-1)(0,0)
\psline[linecolor=blue](-1,-2)(-1,-1)(-1,0)(-1,1)
\pscircle[linecolor=black,linewidth=0.5pt,fillstyle=solid,fillcolor=red](0,0){0.15}
\pscircle[linecolor=black,linewidth=0.5pt,fillstyle=solid,fillcolor=red](1,1){0.15}
\pscircle[linecolor=black,linewidth=0.5pt,fillstyle=solid,fillcolor=red](2,2){0.15}
\pscircle[linecolor=black,linewidth=0.5pt,fillstyle=solid,fillcolor=red](3,3){0.15}
\psarc[linecolor=black,linewidth=0.5pt,fillstyle=solid,fillcolor=red](2,2){0.15}{-45.0}{135.0}
\psarc[linecolor=black,linewidth=0.5pt,fillstyle=solid,fillcolor=blue](2,2){0.15}{135.0}{315.0}
\pscircle[linecolor=black,linewidth=0.5pt](2,2){0.15}
\psarc[linecolor=black,linewidth=0.5pt,fillstyle=solid,fillcolor=red](1,1){0.15}{-45.0}{135.0}
\psarc[linecolor=black,linewidth=0.5pt,fillstyle=solid,fillcolor=blue](1,1){0.15}{135.0}{315.0}
\pscircle[linecolor=black,linewidth=0.5pt](1,1){0.15}
\psarc[linecolor=black,linewidth=0.5pt,fillstyle=solid,fillcolor=red](0,0){0.15}{-45.0}{135.0}
\psarc[linecolor=black,linewidth=0.5pt,fillstyle=solid,fillcolor=blue](0,0){0.15}{135.0}{315.0}
\pscircle[linecolor=black,linewidth=0.5pt](0,0){0.15}
\pscircle[linecolor=black,linewidth=0.5pt,fillstyle=solid,fillcolor=blue](-1,1){0.15}
\pscircle[linecolor=black,linewidth=0.5pt,fillstyle=solid,fillcolor=mfzartstahlblau](2,1){0.15}
\pscircle[linecolor=black,linewidth=0.5pt,fillstyle=solid,fillcolor=mfzartstahlblau](1,0){0.15}
\pscircle[linecolor=black,linewidth=0.5pt,fillstyle=solid,fillcolor=mfzartstahlblau](0,-1){0.15}
\pscircle[linecolor=black,linewidth=0.5pt,fillstyle=solid,fillcolor=mfzartstahlblau](-1,-2){0.15}
\pscircle[linecolor=black,linewidth=0.5pt,fillstyle=solid,fillcolor=mfzartstahlblau](2,1){0.15}
\pscircle[linecolor=black,linewidth=0.5pt,fillstyle=solid,fillcolor=blue](2,2){0.15}
\pscircle[linecolor=black,linewidth=0.5pt,fillstyle=solid,fillcolor=mfzartstahlblau](1,0){0.15}
\psarc[linecolor=black,linewidth=0.5pt,fillstyle=solid,fillcolor=red](1,1){0.15}{-45.0}{135.0}
\psarc[linecolor=black,linewidth=0.5pt,fillstyle=solid,fillcolor=blue](1,1){0.15}{135.0}{315.0}
\pscircle[linecolor=black,linewidth=0.5pt](1,1){0.15}
\pscircle[linecolor=black,linewidth=0.5pt,fillstyle=solid,fillcolor=mfzartstahlblau](0,-1){0.15}
\psarc[linecolor=black,linewidth=0.5pt,fillstyle=solid,fillcolor=red](0,0){0.15}{-45.0}{135.0}
\psarc[linecolor=black,linewidth=0.5pt,fillstyle=solid,fillcolor=blue](0,0){0.15}{135.0}{315.0}
\pscircle[linecolor=black,linewidth=0.5pt](0,0){0.15}
\psarc[linecolor=black,linewidth=0.5pt,fillstyle=solid,fillcolor=red](2,2){0.15}{-45.0}{135.0}
\psarc[linecolor=black,linewidth=0.5pt,fillstyle=solid,fillcolor=blue](2,2){0.15}{135.0}{315.0}
\pscircle[linecolor=black,linewidth=0.5pt](2,2){0.15}
\pscircle[linecolor=black,linewidth=0.5pt,fillstyle=solid,fillcolor=red](1,3){0.15}
\pscircle[linecolor=black,linewidth=0.5pt,fillstyle=solid,fillcolor=lightcoral](3,3){0.15}
\pscircle[linecolor=black,linewidth=0.5pt,fillstyle=solid,fillcolor=mfzartstahlblau](-1,-2){0.15}
\pscircle[linecolor=black,linewidth=0.5pt,fillstyle=solid,fillcolor=blue](-1,1){0.15}
\end{pspicture}\hfil
\psset{unit=0.5cm}
\begin{pspicture}(-3.95,-4.7)(1.95,2.2)
\pspolygon[linecolor=white,fillstyle=solid,fillcolor=backgroundgray,linearc=0.3](-3.95,-4.7)(1.95,-4.7)(1.95,2.2)(-3.95,2.2)

\psset{linewidth=0.5pt,linecolor=gray,linestyle=solid,fillstyle=none}
\psline[linecolor=lightgray](-3,-4)(-3,1)
\psline[linecolor=lightgray](-2,-4)(-2,1)
\psline[linecolor=lightgray](-1,-4)(-1,1)
\psline[linecolor=lightgray](0,-4)(0,1)
\psline[linecolor=lightgray](1,-4)(1,1)
\psline[linecolor=lightgray](-3,-4)(1,-4)
\psline[linecolor=lightgray](-3,-3)(1,-3)
\psline[linecolor=lightgray](-3,-2)(1,-2)
\psline[linecolor=lightgray](-3,-1)(1,-1)
\psline[linecolor=lightgray](-3,0)(1,0)
\psline[linecolor=lightgray](-3,1)(1,1)
\psline[linecolor=black,arrows=->](-3.75,0)(1.75,0)
\psline[linecolor=black,arrows=->](0,-4.5)(0,2)
\psline(-3,-3)(1,1)
\psline[linecolor=magenta,linestyle=dashed](-2,-4)(1,-1)
\psset{linewidth=1pt,linecolor=black,linestyle=solid,fillstyle=none}
\psline[linearc=0.2,linecolor=blue](-2,-3)(-2,-2)
\psline[linecolor=blue](-2,-3)(-2,-2)
\pscircle[linecolor=black,linewidth=0.5pt,fillstyle=solid,fillcolor=red](0,0){0.15}
\pscircle[linecolor=black,linewidth=0.5pt,fillstyle=solid,fillcolor=blue](-2,-2){0.15}
\pscircle[linecolor=black,linewidth=0.5pt,fillstyle=solid,fillcolor=mfzartstahlblau](-2,-3){0.15}
\end{pspicture}\end{center}
\begin{center}
\psset{unit=0.5cm}
\begin{pspicture}(-3.95,-4.7)(5.95,6.2)
\pspolygon[linecolor=white,fillstyle=solid,fillcolor=backgroundgray,linearc=0.3](-3.95,-4.7)(5.95,-4.7)(5.95,6.2)(-3.95,6.2)

\psset{linewidth=0.5pt,linecolor=gray,linestyle=solid,fillstyle=none}
\psline[linecolor=lightgray](-3,-4)(-3,5)
\psline[linecolor=lightgray](-2,-4)(-2,5)
\psline[linecolor=lightgray](-1,-4)(-1,5)
\psline[linecolor=lightgray](0,-4)(0,5)
\psline[linecolor=lightgray](1,-4)(1,5)
\psline[linecolor=lightgray](2,-4)(2,5)
\psline[linecolor=lightgray](3,-4)(3,5)
\psline[linecolor=lightgray](4,-4)(4,5)
\psline[linecolor=lightgray](5,-4)(5,5)
\psline[linecolor=lightgray](-3,-4)(5,-4)
\psline[linecolor=lightgray](-3,-3)(5,-3)
\psline[linecolor=lightgray](-3,-2)(5,-2)
\psline[linecolor=lightgray](-3,-1)(5,-1)
\psline[linecolor=lightgray](-3,0)(5,0)
\psline[linecolor=lightgray](-3,1)(5,1)
\psline[linecolor=lightgray](-3,2)(5,2)
\psline[linecolor=lightgray](-3,3)(5,3)
\psline[linecolor=lightgray](-3,4)(5,4)
\psline[linecolor=lightgray](-3,5)(5,5)
\psline[linecolor=black,arrows=->](-3.75,0)(5.75,0)
\psline[linecolor=black,arrows=->](0,-4.5)(0,6)
\psline(-3,-3)(5,5)
\psline[linecolor=magenta,linestyle=dashed](-2,-4)(5,3)
\psset{linewidth=1pt,linecolor=black,linestyle=solid,fillstyle=none}
\psline[linearc=0.2,linecolor=blue](2,1)(2,2)
\psline[linecolor=blue](2,1)(2,2)
\psline[linearc=0.2,linecolor=blue](1,0)(1,1)
\psline[linecolor=blue](2,1)(2,2)
\psline[linecolor=blue](1,0)(1,1)
\psline[linearc=0.2,linecolor=blue](0,-1)(0,0)
\psline[linecolor=blue](2,1)(2,2)
\psline[linecolor=blue](1,0)(1,1)
\psline[linecolor=blue](0,-1)(0,0)
\psline[linearc=0.2,linecolor=blue](-1,-2)(-1,-1)(-1,0)(-1,1)
\psline[linearc=0.2,linecolor=blue](1,3)(2,3)(3,3)
\psline[linecolor=blue](2,1)(2,2)
\psline[linecolor=blue](1,0)(1,1)
\psline[linecolor=blue](0,-1)(0,0)
\psline[linecolor=blue](-1,-2)(-1,-1)(-1,0)(-1,1)
\psline[linearc=0.2,linecolor=blue](-2,-3)(-2,-2)(-2,-1)(-2,0)(-2,1)(-2,2)
\psline[linearc=0.2,linecolor=blue](0,4)(1,4)(2,4)(3,4)(4,4)
\psline[linecolor=blue](2,1)(2,2)
\psline[linecolor=blue](1,0)(1,1)
\psline[linecolor=blue](0,-1)(0,0)
\psline[linecolor=blue](-1,-2)(-1,-1)(-1,0)(-1,1)
\psline[linecolor=blue](-2,-3)(-2,-2)(-2,-1)(-2,0)(-2,1)(-2,2)
\pscircle[linecolor=black,linewidth=0.5pt,fillstyle=solid,fillcolor=red](0,0){0.15}
\pscircle[linecolor=black,linewidth=0.5pt,fillstyle=solid,fillcolor=red](1,1){0.15}
\pscircle[linecolor=black,linewidth=0.5pt,fillstyle=solid,fillcolor=red](2,2){0.15}
\pscircle[linecolor=black,linewidth=0.5pt,fillstyle=solid,fillcolor=red](3,3){0.15}
\pscircle[linecolor=black,linewidth=0.5pt,fillstyle=solid,fillcolor=red](4,4){0.15}
\psarc[linecolor=black,linewidth=0.5pt,fillstyle=solid,fillcolor=red](2,2){0.15}{-45.0}{135.0}
\psarc[linecolor=black,linewidth=0.5pt,fillstyle=solid,fillcolor=blue](2,2){0.15}{135.0}{315.0}
\pscircle[linecolor=black,linewidth=0.5pt](2,2){0.15}
\psarc[linecolor=black,linewidth=0.5pt,fillstyle=solid,fillcolor=red](1,1){0.15}{-45.0}{135.0}
\psarc[linecolor=black,linewidth=0.5pt,fillstyle=solid,fillcolor=blue](1,1){0.15}{135.0}{315.0}
\pscircle[linecolor=black,linewidth=0.5pt](1,1){0.15}
\psarc[linecolor=black,linewidth=0.5pt,fillstyle=solid,fillcolor=red](0,0){0.15}{-45.0}{135.0}
\psarc[linecolor=black,linewidth=0.5pt,fillstyle=solid,fillcolor=blue](0,0){0.15}{135.0}{315.0}
\pscircle[linecolor=black,linewidth=0.5pt](0,0){0.15}
\pscircle[linecolor=black,linewidth=0.5pt,fillstyle=solid,fillcolor=blue](-1,1){0.15}
\pscircle[linecolor=black,linewidth=0.5pt,fillstyle=solid,fillcolor=blue](-2,2){0.15}
\pscircle[linecolor=black,linewidth=0.5pt,fillstyle=solid,fillcolor=mfzartstahlblau](2,1){0.15}
\pscircle[linecolor=black,linewidth=0.5pt,fillstyle=solid,fillcolor=mfzartstahlblau](1,0){0.15}
\pscircle[linecolor=black,linewidth=0.5pt,fillstyle=solid,fillcolor=mfzartstahlblau](0,-1){0.15}
\pscircle[linecolor=black,linewidth=0.5pt,fillstyle=solid,fillcolor=mfzartstahlblau](-1,-2){0.15}
\pscircle[linecolor=black,linewidth=0.5pt,fillstyle=solid,fillcolor=mfzartstahlblau](-2,-3){0.15}
\pscircle[linecolor=black,linewidth=0.5pt,fillstyle=solid,fillcolor=mfzartstahlblau](2,1){0.15}
\pscircle[linecolor=black,linewidth=0.5pt,fillstyle=solid,fillcolor=blue](2,2){0.15}
\pscircle[linecolor=black,linewidth=0.5pt,fillstyle=solid,fillcolor=mfzartstahlblau](1,0){0.15}
\psarc[linecolor=black,linewidth=0.5pt,fillstyle=solid,fillcolor=red](1,1){0.15}{-45.0}{135.0}
\psarc[linecolor=black,linewidth=0.5pt,fillstyle=solid,fillcolor=blue](1,1){0.15}{135.0}{315.0}
\pscircle[linecolor=black,linewidth=0.5pt](1,1){0.15}
\pscircle[linecolor=black,linewidth=0.5pt,fillstyle=solid,fillcolor=mfzartstahlblau](0,-1){0.15}
\psarc[linecolor=black,linewidth=0.5pt,fillstyle=solid,fillcolor=red](0,0){0.15}{-45.0}{135.0}
\psarc[linecolor=black,linewidth=0.5pt,fillstyle=solid,fillcolor=blue](0,0){0.15}{135.0}{315.0}
\pscircle[linecolor=black,linewidth=0.5pt](0,0){0.15}
\psarc[linecolor=black,linewidth=0.5pt,fillstyle=solid,fillcolor=red](2,2){0.15}{-45.0}{135.0}
\psarc[linecolor=black,linewidth=0.5pt,fillstyle=solid,fillcolor=blue](2,2){0.15}{135.0}{315.0}
\pscircle[linecolor=black,linewidth=0.5pt](2,2){0.15}
\pscircle[linecolor=black,linewidth=0.5pt,fillstyle=solid,fillcolor=red](1,3){0.15}
\pscircle[linecolor=black,linewidth=0.5pt,fillstyle=solid,fillcolor=lightcoral](3,3){0.15}
\pscircle[linecolor=black,linewidth=0.5pt,fillstyle=solid,fillcolor=mfzartstahlblau](-1,-2){0.15}
\pscircle[linecolor=black,linewidth=0.5pt,fillstyle=solid,fillcolor=blue](-1,1){0.15}
\pscircle[linecolor=black,linewidth=0.5pt,fillstyle=solid,fillcolor=red](0,4){0.15}
\pscircle[linecolor=black,linewidth=0.5pt,fillstyle=solid,fillcolor=lightcoral](4,4){0.15}
\pscircle[linecolor=black,linewidth=0.5pt,fillstyle=solid,fillcolor=mfzartstahlblau](-2,-3){0.15}
\pscircle[linecolor=black,linewidth=0.5pt,fillstyle=solid,fillcolor=blue](-2,2){0.15}
\end{pspicture}\hfil
\psset{unit=0.5cm}
\begin{pspicture}(-4.95,-5.7)(2.95,3.2)
\pspolygon[linecolor=white,fillstyle=solid,fillcolor=backgroundgray,linearc=0.3](-4.95,-5.7)(2.95,-5.7)(2.95,3.2)(-4.95,3.2)

\psset{linewidth=0.5pt,linecolor=gray,linestyle=solid,fillstyle=none}
\psline[linecolor=lightgray](-4,-5)(-4,2)
\psline[linecolor=lightgray](-3,-5)(-3,2)
\psline[linecolor=lightgray](-2,-5)(-2,2)
\psline[linecolor=lightgray](-1,-5)(-1,2)
\psline[linecolor=lightgray](0,-5)(0,2)
\psline[linecolor=lightgray](1,-5)(1,2)
\psline[linecolor=lightgray](2,-5)(2,2)
\psline[linecolor=lightgray](-4,-5)(2,-5)
\psline[linecolor=lightgray](-4,-4)(2,-4)
\psline[linecolor=lightgray](-4,-3)(2,-3)
\psline[linecolor=lightgray](-4,-2)(2,-2)
\psline[linecolor=lightgray](-4,-1)(2,-1)
\psline[linecolor=lightgray](-4,0)(2,0)
\psline[linecolor=lightgray](-4,1)(2,1)
\psline[linecolor=lightgray](-4,2)(2,2)
\psline[linecolor=black,arrows=->](-4.75,0)(2.75,0)
\psline[linecolor=black,arrows=->](0,-5.5)(0,3)
\psline(-4,-4)(2,2)
\psline[linecolor=magenta,linestyle=dashed](-3,-5)(2,0)
\psset{linewidth=1pt,linecolor=black,linestyle=solid,fillstyle=none}
\psline[linearc=0.2,linecolor=blue](-2,-3)(-2,-2)
\psline[linecolor=blue](-2,-3)(-2,-2)
\psline[linearc=0.2,linecolor=blue](-3,-4)(-3,-3)(-3,-2)(-3,-1)
\psline[linearc=0.2,linecolor=blue](-1,1)(0,1)(1,1)
\psline[linecolor=blue](-2,-3)(-2,-2)
\psline[linecolor=blue](-3,-4)(-3,-3)(-3,-2)(-3,-1)
\pscircle[linecolor=black,linewidth=0.5pt,fillstyle=solid,fillcolor=red](0,0){0.15}
\pscircle[linecolor=black,linewidth=0.5pt,fillstyle=solid,fillcolor=red](1,1){0.15}
\pscircle[linecolor=black,linewidth=0.5pt,fillstyle=solid,fillcolor=blue](-2,-2){0.15}
\pscircle[linecolor=black,linewidth=0.5pt,fillstyle=solid,fillcolor=blue](-3,-1){0.15}
\pscircle[linecolor=black,linewidth=0.5pt,fillstyle=solid,fillcolor=mfzartstahlblau](-2,-3){0.15}
\pscircle[linecolor=black,linewidth=0.5pt,fillstyle=solid,fillcolor=mfzartstahlblau](-3,-4){0.15}
\pscircle[linecolor=black,linewidth=0.5pt,fillstyle=solid,fillcolor=red](-1,1){0.15}
\pscircle[linecolor=black,linewidth=0.5pt,fillstyle=solid,fillcolor=lightcoral](1,1){0.15}
\pscircle[linecolor=black,linewidth=0.5pt,fillstyle=solid,fillcolor=mfzartstahlblau](-3,-4){0.15}
\pscircle[linecolor=black,linewidth=0.5pt,fillstyle=solid,fillcolor=blue](-3,-1){0.15}
\end{pspicture}\end{center}
The same reasoning as for the odd identity concludes the proof.
\end{proof}

\secA{Bijective proof of Cigler's Conjecture for $k>1$}
\label{sec:proof}
The rest of this paper is devoted to a bijective proof of Ciglers's Conjecture for $k>1$.
While all single arguments are elementary, the chain of arguments
is long and complicated: we will make extensive use of \EM{pictures}
to explain them, hoping that the saying ``A picture is worth a thousand words'' proves to be true.

\secB{Encoding permutations of survivors}
\label{sec:rhs-sign}
Let us call a permutation $\pi$ \EM{admissible} if there are
survivors ($N$--tuples of \nilp s)
connecting initial point $A_i$ to terminal point $B_{\pi\of i}$ for
$i=0,1,\dots,N-1$.

Note that for every survivor, the lattice path reaching terminal
point $B_i$ may have length $0$ (if $B_i$ is also an enforced
initial point: this can only happen for the \lhs; see
subsection~\ref{sec:even-odd}), but otherwise must end
\bit
\item either with an upwards step \EM{from below}
\item or with a rightwards step \EM{from the left};
\eit
see the following pictures:
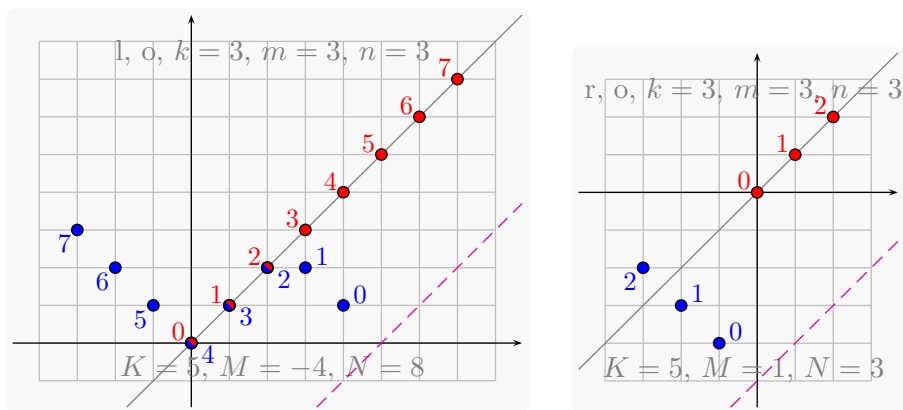
\begin{center}
\psset{unit=0.5cm}
\begin{pspicture}(-4.9,-1.9)(8.9,8.9)
\pspolygon[linecolor=white,fillstyle=solid,fillcolor=backgroundgray,linearc=0.3](-4.9,-1.9)(8.9,-1.9)(8.9,8.9)(-4.9,8.9)

\psset{linewidth=0.5pt,linecolor=gray,linestyle=solid,fillstyle=none}
\psline[linecolor=lightgray](-4,-1)(-4,8)
\psline[linecolor=lightgray](-3,-1)(-3,8)
\psline[linecolor=lightgray](-2,-1)(-2,8)
\psline[linecolor=lightgray](-1,-1)(-1,8)
\psline[linecolor=lightgray](0,-1)(0,8)
\psline[linecolor=lightgray](1,-1)(1,8)
\psline[linecolor=lightgray](2,-1)(2,8)
\psline[linecolor=lightgray](3,-1)(3,8)
\psline[linecolor=lightgray](4,-1)(4,8)
\psline[linecolor=lightgray](5,-1)(5,8)
\psline[linecolor=lightgray](6,-1)(6,8)
\psline[linecolor=lightgray](7,-1)(7,8)
\psline[linecolor=lightgray](8,-1)(8,8)
\psline[linecolor=lightgray](-4,-1)(8,-1)
\psline[linecolor=lightgray](-4,0)(8,0)
\psline[linecolor=lightgray](-4,1)(8,1)
\psline[linecolor=lightgray](-4,2)(8,2)
\psline[linecolor=lightgray](-4,3)(8,3)
\psline[linecolor=lightgray](-4,4)(8,4)
\psline[linecolor=lightgray](-4,5)(8,5)
\psline[linecolor=lightgray](-4,6)(8,6)
\psline[linecolor=lightgray](-4,7)(8,7)
\psline[linecolor=lightgray](-4,8)(8,8)
\psline[linecolor=black,arrows=->](-4.7,0)(8.7,0)
\psline[linecolor=black,arrows=->](0,-1.7)(0,8.7)
\rput(2,7.7){ {\small\gray l, o, $k=3$, $m=3$, $n=3$}}
\rput(2,-0.7){ {\small\gray $K=5$, $M=-4$, $N=8$}}
\psset{linewidth=0.5pt,linecolor=gray,linestyle=solid,fillstyle=none}
\psline(-1.7,-1.7)(8.7,8.7)
\psline[linecolor=magenta,linestyle=dashed](3.3,-1.7)(8.7,3.7)
\pscircle[linecolor=black,linewidth=0.5pt,fillstyle=solid,fillcolor=red](0,0){0.15}
\rput(-0.48,0.3){ {\footnotesize\red $0$}}
\pscircle[linecolor=black,linewidth=0.5pt,fillstyle=solid,fillcolor=red](1,1){0.15}
\rput(0.52,1.3){ {\footnotesize\red $1$}}
\pscircle[linecolor=black,linewidth=0.5pt,fillstyle=solid,fillcolor=red](2,2){0.15}
\rput(1.52,2.3){ {\footnotesize\red $2$}}
\pscircle[linecolor=black,linewidth=0.5pt,fillstyle=solid,fillcolor=red](3,3){0.15}
\rput(2.52,3.3){ {\footnotesize\red $3$}}
\pscircle[linecolor=black,linewidth=0.5pt,fillstyle=solid,fillcolor=red](4,4){0.15}
\rput(3.52,4.3){ {\footnotesize\red $4$}}
\pscircle[linecolor=black,linewidth=0.5pt,fillstyle=solid,fillcolor=red](5,5){0.15}
\rput(4.52,5.3){ {\footnotesize\red $5$}}
\pscircle[linecolor=black,linewidth=0.5pt,fillstyle=solid,fillcolor=red](6,6){0.15}
\rput(5.52,6.3){ {\footnotesize\red $6$}}
\pscircle[linecolor=black,linewidth=0.5pt,fillstyle=solid,fillcolor=red](7,7){0.15}
\rput(6.52,7.3){ {\footnotesize\red $7$}}
\psarc[linecolor=black,linewidth=0.5pt,fillstyle=solid,fillcolor=red](2,2){0.15}{-45.0}{135.0}
\psarc[linecolor=black,linewidth=0.5pt,fillstyle=solid,fillcolor=blue](2,2){0.15}{135.0}{315.0}
\pscircle[linecolor=black,linewidth=0.5pt](2,2){0.15}
\rput(2.3,1.7){ {\footnotesize\blue $2$}}
\psarc[linecolor=black,linewidth=0.5pt,fillstyle=solid,fillcolor=red](1,1){0.15}{-45.0}{135.0}
\psarc[linecolor=black,linewidth=0.5pt,fillstyle=solid,fillcolor=blue](1,1){0.15}{135.0}{315.0}
\pscircle[linecolor=black,linewidth=0.5pt](1,1){0.15}
\rput(1.3,0.7){ {\footnotesize\blue $3$}}
\psarc[linecolor=black,linewidth=0.5pt,fillstyle=solid,fillcolor=red](0,0){0.15}{-45.0}{135.0}
\psarc[linecolor=black,linewidth=0.5pt,fillstyle=solid,fillcolor=blue](0,0){0.15}{135.0}{315.0}
\pscircle[linecolor=black,linewidth=0.5pt](0,0){0.15}
\rput(0.3,-0.3){ {\footnotesize\blue $4$}}
\pscircle[linecolor=black,linewidth=0.5pt,fillstyle=solid,fillcolor=blue](4,1){0.15}
\rput(4.3,1.3){ {\footnotesize\blue $0$}}
\pscircle[linecolor=black,linewidth=0.5pt,fillstyle=solid,fillcolor=blue](3,2){0.15}
\rput(3.3,2.3){ {\footnotesize\blue $1$}}
\pscircle[linecolor=black,linewidth=0.5pt,fillstyle=solid,fillcolor=blue](-1,1){0.15}
\rput(-1.48,0.64){ {\footnotesize\blue $5$}}
\pscircle[linecolor=black,linewidth=0.5pt,fillstyle=solid,fillcolor=blue](-2,2){0.15}
\rput(-2.48,1.64){ {\footnotesize\blue $6$}}
\pscircle[linecolor=black,linewidth=0.5pt,fillstyle=solid,fillcolor=blue](-3,3){0.15}
\rput(-3.48,2.64){ {\footnotesize\blue $7$}}
\end{pspicture}\hfil
\psset{unit=0.5cm}
\begin{pspicture}(-4.9,-5.9)(3.9,3.9)
\pspolygon[linecolor=white,fillstyle=solid,fillcolor=backgroundgray,linearc=0.3](-4.9,-5.9)(3.9,-5.9)(3.9,3.9)(-4.9,3.9)

\psset{linewidth=0.5pt,linecolor=gray,linestyle=solid,fillstyle=none}
\psline[linecolor=lightgray](-4,-5)(-4,3)
\psline[linecolor=lightgray](-3,-5)(-3,3)
\psline[linecolor=lightgray](-2,-5)(-2,3)
\psline[linecolor=lightgray](-1,-5)(-1,3)
\psline[linecolor=lightgray](0,-5)(0,3)
\psline[linecolor=lightgray](1,-5)(1,3)
\psline[linecolor=lightgray](2,-5)(2,3)
\psline[linecolor=lightgray](3,-5)(3,3)
\psline[linecolor=lightgray](-4,-5)(3,-5)
\psline[linecolor=lightgray](-4,-4)(3,-4)
\psline[linecolor=lightgray](-4,-3)(3,-3)
\psline[linecolor=lightgray](-4,-2)(3,-2)
\psline[linecolor=lightgray](-4,-1)(3,-1)
\psline[linecolor=lightgray](-4,0)(3,0)
\psline[linecolor=lightgray](-4,1)(3,1)
\psline[linecolor=lightgray](-4,2)(3,2)
\psline[linecolor=lightgray](-4,3)(3,3)
\psline[linecolor=black,arrows=->](-4.7,0)(3.7,0)
\psline[linecolor=black,arrows=->](0,-5.7)(0,3.7)
\rput(-0.5,2.7){ {\small\gray r, o, $k=3$, $m=3$, $n=3$}}
\rput(-0.5,-4.7){ {\small\gray $K=5$, $M=1$, $N=3$}}
\psset{linewidth=0.5pt,linecolor=gray,linestyle=solid,fillstyle=none}
\psline(-4.7,-4.7)(3.7,3.7)
\psline[linecolor=magenta,linestyle=dashed](-0.7,-5.7)(3.7,-1.3)
\pscircle[linecolor=black,linewidth=0.5pt,fillstyle=solid,fillcolor=red](0,0){0.15}
\rput(-0.48,0.3){ {\footnotesize\red $0$}}
\pscircle[linecolor=black,linewidth=0.5pt,fillstyle=solid,fillcolor=red](1,1){0.15}
\rput(0.52,1.3){ {\footnotesize\red $1$}}
\pscircle[linecolor=black,linewidth=0.5pt,fillstyle=solid,fillcolor=red](2,2){0.15}
\rput(1.52,2.3){ {\footnotesize\red $2$}}
\pscircle[linecolor=black,linewidth=0.5pt,fillstyle=solid,fillcolor=blue](-1,-4){0.15}
\rput(-0.7,-3.7){ {\footnotesize\blue $0$}}
\pscircle[linecolor=black,linewidth=0.5pt,fillstyle=solid,fillcolor=blue](-2,-3){0.15}
\rput(-1.7,-2.7){ {\footnotesize\blue $1$}}
\pscircle[linecolor=black,linewidth=0.5pt,fillstyle=solid,fillcolor=blue](-3,-2){0.15}
\rput(-3.48,-2.36){ {\footnotesize\blue $2$}}
\end{pspicture}\end{center}
Considering the \lhs, note that all lattice paths starting
at enforced initial points
\bit
\item \EM{below} the diagonal must end with an upwards step,
\item \EM{on} the diagonal have length $0$,
\item \EM{above} the diagonal must end with a rightwards step.
\eit
This means that \EM{all} admissible permutations for the \lhs\ must map 
\bit
\item indices $k-1,k,k+1,\dots,k+\pas{m+\Iverson{\text{even situation}}}-2$
\item to indices $0,1,\dots,\pas{m+\Iverson{\text{even situation}}}-1$
\eit
\EM{in reverse order} (in the above picture: $2\mapsto2$, $3\mapsto1$ and 
$4\mapsto0$). We may disregard this ``enforced subpermutation'' and
renumber the initial and terminal points as indicated in the following
picture: 
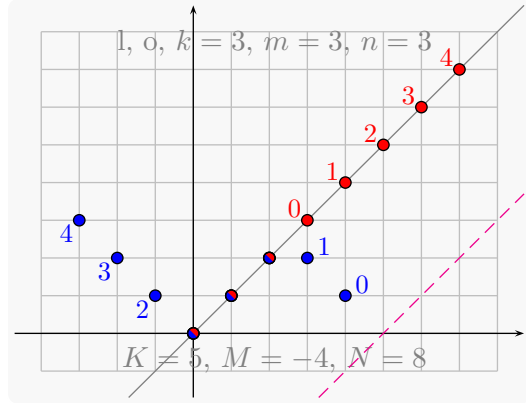
\cps{points_numbered_with_terminals_5_-4_8_renumbered}

With this renumbering for the \lhs\ (and the original numbering
for the \rhs), we assign to $\pi$ the $01$--code (omitting terminal
point $B_0$)
$$
c\of\pi = \pas{\Iverson{\text{Terminal point $B_i$ is reached from the left}}}_{i\geq 1}.
$$

By considering \EM{enforced terminal segments}, as illustrated in the
following pictures, is easy see that the mapping $\pi\mapsto c\of\pi$ is \EM{injective}:
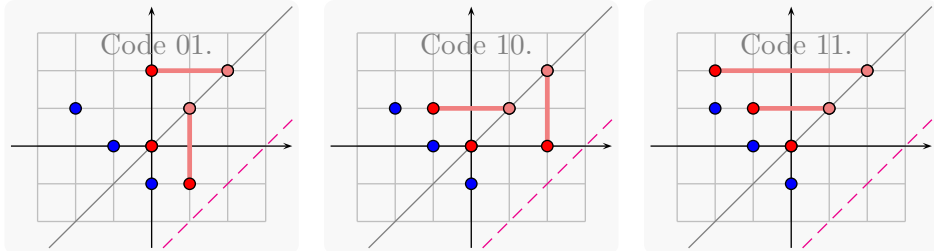
\begin{center}
\psset{unit=0.5cm}
\begin{pspicture}(-3.9,-2.9)(3.9,3.9)
\pspolygon[linecolor=white,fillstyle=solid,fillcolor=backgroundgray,linearc=0.3](-3.9,-2.9)(3.9,-2.9)(3.9,3.9)(-3.9,3.9)

\psset{linewidth=0.5pt,linecolor=gray,linestyle=solid,fillstyle=none}
\psline[linecolor=lightgray](-3,-2)(-3,3)
\psline[linecolor=lightgray](-2,-2)(-2,3)
\psline[linecolor=lightgray](-1,-2)(-1,3)
\psline[linecolor=lightgray](0,-2)(0,3)
\psline[linecolor=lightgray](1,-2)(1,3)
\psline[linecolor=lightgray](2,-2)(2,3)
\psline[linecolor=lightgray](3,-2)(3,3)
\psline[linecolor=lightgray](-3,-2)(3,-2)
\psline[linecolor=lightgray](-3,-1)(3,-1)
\psline[linecolor=lightgray](-3,0)(3,0)
\psline[linecolor=lightgray](-3,1)(3,1)
\psline[linecolor=lightgray](-3,2)(3,2)
\psline[linecolor=lightgray](-3,3)(3,3)
\psline[linecolor=black,arrows=->](-3.7,0)(3.7,0)
\psline[linecolor=black,arrows=->](0,-2.7)(0,3.7)
\rput(0,2.7){ {\small\gray Code 01.}}
\psset{linewidth=0.5pt,linecolor=gray,linestyle=solid,fillstyle=none}
\psline(-2.7,-2.7)(3.7,3.7)
\psline[linecolor=magenta,linestyle=dashed](0.3,-2.7)(3.7,0.7)
\pscircle[linecolor=black,linewidth=0.5pt,fillstyle=solid,fillcolor=red](0,0){0.15}
\pscircle[linecolor=black,linewidth=0.5pt,fillstyle=solid,fillcolor=red](1,1){0.15}
\pscircle[linecolor=black,linewidth=0.5pt,fillstyle=solid,fillcolor=red](2,2){0.15}
\pscircle[linecolor=black,linewidth=0.5pt,fillstyle=solid,fillcolor=blue](0,-1){0.15}
\pscircle[linecolor=black,linewidth=0.5pt,fillstyle=solid,fillcolor=blue](-1,0){0.15}
\pscircle[linecolor=black,linewidth=0.5pt,fillstyle=solid,fillcolor=blue](-2,1){0.15}
\pscircle[linecolor=black,linewidth=0.5pt,fillstyle=solid,fillcolor=red](0,0){0.15}
\psline[linewidth=0.125,linearc=0.2,arrows=c-c,linecolor=lightcoral](1,1)(1,-1)
\pscircle[linecolor=black,linewidth=0.5pt,fillstyle=solid,fillcolor=lightcoral](1,1){0.15}
\pscircle[linecolor=black,linewidth=0.5pt,fillstyle=solid,fillcolor=red](1,-1){0.15}
\psline[linewidth=0.125,linearc=0.2,arrows=c-c,linecolor=lightcoral](2,2)(0,2)
\pscircle[linecolor=black,linewidth=0.5pt,fillstyle=solid,fillcolor=lightcoral](2,2){0.15}
\pscircle[linecolor=black,linewidth=0.5pt,fillstyle=solid,fillcolor=red](0,2){0.15}
\end{pspicture}\hfil
\psset{unit=0.5cm}
\begin{pspicture}(-3.9,-2.9)(3.9,3.9)
\pspolygon[linecolor=white,fillstyle=solid,fillcolor=backgroundgray,linearc=0.3](-3.9,-2.9)(3.9,-2.9)(3.9,3.9)(-3.9,3.9)

\psset{linewidth=0.5pt,linecolor=gray,linestyle=solid,fillstyle=none}
\psline[linecolor=lightgray](-3,-2)(-3,3)
\psline[linecolor=lightgray](-2,-2)(-2,3)
\psline[linecolor=lightgray](-1,-2)(-1,3)
\psline[linecolor=lightgray](0,-2)(0,3)
\psline[linecolor=lightgray](1,-2)(1,3)
\psline[linecolor=lightgray](2,-2)(2,3)
\psline[linecolor=lightgray](3,-2)(3,3)
\psline[linecolor=lightgray](-3,-2)(3,-2)
\psline[linecolor=lightgray](-3,-1)(3,-1)
\psline[linecolor=lightgray](-3,0)(3,0)
\psline[linecolor=lightgray](-3,1)(3,1)
\psline[linecolor=lightgray](-3,2)(3,2)
\psline[linecolor=lightgray](-3,3)(3,3)
\psline[linecolor=black,arrows=->](-3.7,0)(3.7,0)
\psline[linecolor=black,arrows=->](0,-2.7)(0,3.7)
\rput(0,2.7){ {\small\gray Code 10.}}
\psset{linewidth=0.5pt,linecolor=gray,linestyle=solid,fillstyle=none}
\psline(-2.7,-2.7)(3.7,3.7)
\psline[linecolor=magenta,linestyle=dashed](0.3,-2.7)(3.7,0.7)
\pscircle[linecolor=black,linewidth=0.5pt,fillstyle=solid,fillcolor=red](0,0){0.15}
\pscircle[linecolor=black,linewidth=0.5pt,fillstyle=solid,fillcolor=red](1,1){0.15}
\pscircle[linecolor=black,linewidth=0.5pt,fillstyle=solid,fillcolor=red](2,2){0.15}
\pscircle[linecolor=black,linewidth=0.5pt,fillstyle=solid,fillcolor=blue](0,-1){0.15}
\pscircle[linecolor=black,linewidth=0.5pt,fillstyle=solid,fillcolor=blue](-1,0){0.15}
\pscircle[linecolor=black,linewidth=0.5pt,fillstyle=solid,fillcolor=blue](-2,1){0.15}
\pscircle[linecolor=black,linewidth=0.5pt,fillstyle=solid,fillcolor=red](0,0){0.15}
\psline[linewidth=0.125,linearc=0.2,arrows=c-c,linecolor=lightcoral](1,1)(-1,1)
\pscircle[linecolor=black,linewidth=0.5pt,fillstyle=solid,fillcolor=lightcoral](1,1){0.15}
\pscircle[linecolor=black,linewidth=0.5pt,fillstyle=solid,fillcolor=red](-1,1){0.15}
\psline[linewidth=0.125,linearc=0.2,arrows=c-c,linecolor=lightcoral](2,2)(2,0)
\pscircle[linecolor=black,linewidth=0.5pt,fillstyle=solid,fillcolor=lightcoral](2,2){0.15}
\pscircle[linecolor=black,linewidth=0.5pt,fillstyle=solid,fillcolor=red](2,0){0.15}
\end{pspicture}\hfil
\psset{unit=0.5cm}
\begin{pspicture}(-3.9,-2.9)(3.9,3.9)
\pspolygon[linecolor=white,fillstyle=solid,fillcolor=backgroundgray,linearc=0.3](-3.9,-2.9)(3.9,-2.9)(3.9,3.9)(-3.9,3.9)

\psset{linewidth=0.5pt,linecolor=gray,linestyle=solid,fillstyle=none}
\psline[linecolor=lightgray](-3,-2)(-3,3)
\psline[linecolor=lightgray](-2,-2)(-2,3)
\psline[linecolor=lightgray](-1,-2)(-1,3)
\psline[linecolor=lightgray](0,-2)(0,3)
\psline[linecolor=lightgray](1,-2)(1,3)
\psline[linecolor=lightgray](2,-2)(2,3)
\psline[linecolor=lightgray](3,-2)(3,3)
\psline[linecolor=lightgray](-3,-2)(3,-2)
\psline[linecolor=lightgray](-3,-1)(3,-1)
\psline[linecolor=lightgray](-3,0)(3,0)
\psline[linecolor=lightgray](-3,1)(3,1)
\psline[linecolor=lightgray](-3,2)(3,2)
\psline[linecolor=lightgray](-3,3)(3,3)
\psline[linecolor=black,arrows=->](-3.7,0)(3.7,0)
\psline[linecolor=black,arrows=->](0,-2.7)(0,3.7)
\rput(0,2.7){ {\small\gray Code 11.}}
\psset{linewidth=0.5pt,linecolor=gray,linestyle=solid,fillstyle=none}
\psline(-2.7,-2.7)(3.7,3.7)
\psline[linecolor=magenta,linestyle=dashed](0.3,-2.7)(3.7,0.7)
\pscircle[linecolor=black,linewidth=0.5pt,fillstyle=solid,fillcolor=red](0,0){0.15}
\pscircle[linecolor=black,linewidth=0.5pt,fillstyle=solid,fillcolor=red](1,1){0.15}
\pscircle[linecolor=black,linewidth=0.5pt,fillstyle=solid,fillcolor=red](2,2){0.15}
\pscircle[linecolor=black,linewidth=0.5pt,fillstyle=solid,fillcolor=blue](0,-1){0.15}
\pscircle[linecolor=black,linewidth=0.5pt,fillstyle=solid,fillcolor=blue](-1,0){0.15}
\pscircle[linecolor=black,linewidth=0.5pt,fillstyle=solid,fillcolor=blue](-2,1){0.15}
\pscircle[linecolor=black,linewidth=0.5pt,fillstyle=solid,fillcolor=red](0,0){0.15}
\psline[linewidth=0.125,linearc=0.2,linecolor=lightcoral](1,1)(-1,1)
\pscircle[linecolor=black,linewidth=0.5pt,fillstyle=solid,fillcolor=lightcoral](1,1){0.15}
\pscircle[linecolor=black,linewidth=0.5pt,fillstyle=solid,fillcolor=red](-1,1){0.15}
\psline[linewidth=0.125,linearc=0.2,linecolor=lightcoral](2,2)(-2,2)
\pscircle[linecolor=black,linewidth=0.5pt,fillstyle=solid,fillcolor=lightcoral](2,2){0.15}
\pscircle[linecolor=black,linewidth=0.5pt,fillstyle=solid,fillcolor=red](-2,2){0.15}
\end{pspicture}\end{center}
(Note that no admissible permutation is mapped to
$01$--code $00$ in the situation depicted above, since the paths must not intersect the forbidden line.) We call a $01$--code \EM{admissible} if it is assigned
to an admissible permutation, so (by definition) there is a \EM{bijection} between
the set of admissible permutations and the set of admissible codes.

\secB{Relations between signs of admissible permutations and their codes}
We want to express the sign of an admissible permutation $\pi$ in terms
of its corresponding $01$--code $c\of\pi$: as \EM{base cases}, we consider the $01$--codes
\bit
\item consisting of $k-1$ zeros, followed by $n$ ones, for the
	\lhs,
\item consisting of $l\geq0$ zeros, followed by $n-l$ ones, for the \rhs\ 
	(note that $l=1$ means that the permutation has precisely one
	inversion).
\eit
The permutations corresponding to these base cases are visualized in
the following pictures:
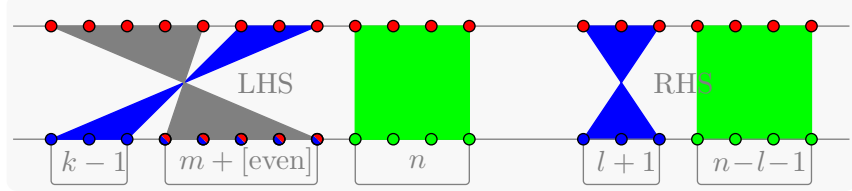
\cps{base_cases}
Clearly, the signs of these permutations are
\bit
\item $\pas{-1}^{\binom{m+\Iverson{\text{even}}+k-1}2}$ for the
	\lhs\ (since the number of inversions of this permutation is ${\binom{m+\Iverson{\text{even}}+k-1}2}$),
\item $\pas{-1}^{\binom{l+1}2}$ for the
	\rhs\ (since the permutation has ${\binom{l+1}2}$
	inversions (note that this is true also for $l=0$).
\eit
For a $01$--code $c=\pas{c_i}_{i\geq 1}$, an \EM{inversion} of $c$
is a pair of indices $\pas{i<j}$ such that $c_i > c_j$. Denote by
$\inv\of c$ the \EM{number} of inversions of $c$.

If we number the zeros in some $01$--code $c$ from left to right
(starting with $1$) and assume that the $i$--th zero appears at
index $z_i$ ($z_i\geq 1$ by definition) in $c$,
then the number of inversions of $c$ equals the sum
$$
\inv\of c  = \sum_{i=1}^{\#\pas{\text{zeros in $c$}}} \pas{z_i-i}.
$$
It is easy to see: a permutation $\pi$ with code $c\of\pi$
\bit
\item 
for the \lhs\ has
\begin{equation}
\label{eq:sign-for-lhs}
{\binom{m+\Iverson{\text{even}}}2} + \inv\of{c\of\pi}
\end{equation}
inversions,
\item  for the \rhs\ has
\begin{equation}
\label{eq:sign-for-rhs}
\binom{l+1}2 + \inv\of{c\of\pi} = \sum_{i=1}^l {z_i}
\end{equation}
inversions, if $c$ contains exactly $l$ zeros.
\eit
(We will use this facts at the end of our proof.)

Assume that in some $01$--code $c=\pas{c_k}_{k\geq 1}$ we have for
$1\leq i < j$
\bit
\item $c_i=1$
\item and $c_j=0$,
\eit
and let $l$ be the number of components $c_k=0$, for $i<k<j$:
Then \EM{swapping} $c_i$ and $c_j$
\bit
\item decreases the number of inversions involving $c_i$ by $l+1$,
\item increases the number of inversions by $1$ for each of the $j-i-1-l$
	ones \EM{between} $i$ and $j$,
\eit
so the change in sign effected by such swap is $\pas{-1}^{\abs{j-i}}$, and the
same holds true for $c_i=0$ and $c_j=1$: we note this simple observation
for later reference.

\secB{Rewrite the identities as equations for sums over $01$--codes}
\def\gfcode#1{{\mathbf{gf}\of{#1}}}

If we consider the set of all admissible codes for the \lhs\ and
for the \rhs\ in equations \eqref{eq:cigler-even} and \eqref{eq:cigler-odd}
and write $\gfcode{c}$ for the  generating function of survivors
corresponding to $c$ (i.e., with permutation $\pi$ of terminal points corresponding to $c$),
then we may
write \eqref{eq:cigler-even} and \eqref{eq:cigler-odd} in a uniform way:
\begin{equation}
\label{eq:cigler-as-sums}
\sum_{C\text{ for LHS}} \signum\pi \gfcode{C}
= 
\sum_{c\text{ for RHS}} \signum\pi \gfcode{c}.
\end{equation}
In general, neither the summation ranges nor the summands of this
identity coincide: our proof will involve
``reducing'' the (larger) sum for the \lhs\ to the sum of the \rhs,
by certain cancellations resembling the \LGV--involution.

\secB{Plan for our bijective proof}
The plan for our bijective proof is simple: 
\begin{enumerate}
\item Find \EM{another} sign--reversing involution $\psi$ on the set of 
all survivors (i.e., \nilp s) corresponding to the \lhs.
\item Show that there is a bijection $\xi$
\bit
\item from the set of the fixed points of $\psi$ (i.e., the \EM{survivors
of the second cancellation} corresponding to $\psi$) for the \lhs
\item to the set of survivors corresponding to the \rhs
\item which changes the sign by a \EM{constant factor} $f$ (i.e., $\signum{\xi\of o} = f\cdot \signum o$), namely
	\bit
	\item $f = \pas{-1}^{\binom{m+k}2}$ in the even case,
	\item $f = \pas{-1}^{\binom{m+k-1}2}$ in the odd case.
	\eit
\eit
\end{enumerate}
As it will turn out, the bijection $\xi$ ``respects the summands in
\eqref{eq:cigler-as-sums}'' in the following sense: 
it implies a bijection between
\bit
\item the summation range (i.e., $01$--codes $C$)
	of \EM{survivors of the second cancellation} (effected by $\psi$)
	appearing in the \lhs\ of \eqref{eq:cigler-as-sums}
\item and the summation range (i.e., $01$--codes $c=\xi\of C$) for the \rhs\ of \eqref{eq:cigler-as-sums},
\eit
and it maps
\bit
\item \EM{survivors of the second cancellation} corresponding
to $01$--code $C$ in the \lhs
\item to survivors corresponding to $01$--code $c=\xi\of C$ in the \rhs.
\eit


\secB{Step 1: Construct the sign--reversing involution $\psi$}
We shall construct the sign--reversing involution $\psi$ on the set of
survivors for the \lhs\ by combining the \LGV--idea
with simple reflections on the diagonal line $d = \pas{y=x}$.
\secC{Reflection for the \lhs: the folded situation}
If we \EM{reflect} on the diagonal $d$ \EM{all} \nilp s for the \lhs\
which lie \EM{above} $d$, then we call the result of such reflection the \EM{folded
situation}. In the following pictures, \EM{reflected} paths 
will be coloured \EM{green}, as well as their reflected initial points
(and all paths and initial points \EM{below} $d$, which are \EM{not reflected},
will be coloured blue, as before).

We illustrate this for the odd identity with  parameters $k=3$, $m=2$ and $n=2$
(i.e., $K=2k-1=5$, $M=2-k-m=-3$, and $N=n+m+k-1 = 6$).
``In terms of initial, two--faced and terminal points'', the transition from the \EM{original}
situation to the \EM{folded} situation looks like this:
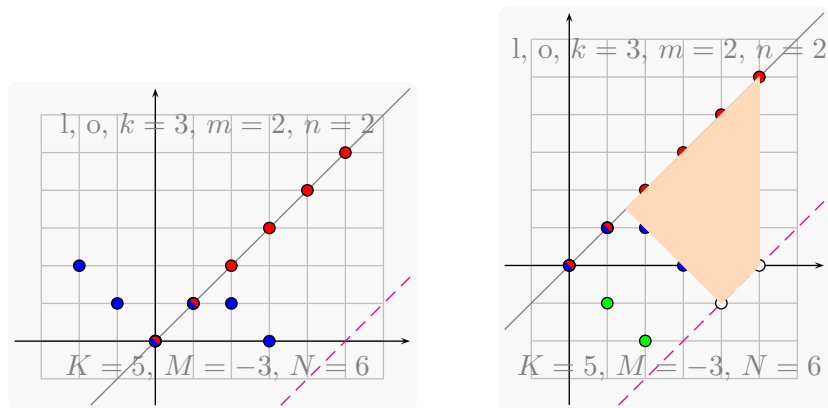
\begin{center}
\psset{unit=0.5cm}
\begin{pspicture}(-3.9,-1.9)(6.9,6.9)
\pspolygon[linecolor=white,fillstyle=solid,fillcolor=backgroundgray,linearc=0.3](-3.9,-1.9)(6.9,-1.9)(6.9,6.9)(-3.9,6.9)

\psset{linewidth=0.5pt,linecolor=gray,linestyle=solid,fillstyle=none}
\psline[linecolor=lightgray](-3,-1)(-3,6)
\psline[linecolor=lightgray](-2,-1)(-2,6)
\psline[linecolor=lightgray](-1,-1)(-1,6)
\psline[linecolor=lightgray](0,-1)(0,6)
\psline[linecolor=lightgray](1,-1)(1,6)
\psline[linecolor=lightgray](2,-1)(2,6)
\psline[linecolor=lightgray](3,-1)(3,6)
\psline[linecolor=lightgray](4,-1)(4,6)
\psline[linecolor=lightgray](5,-1)(5,6)
\psline[linecolor=lightgray](6,-1)(6,6)
\psline[linecolor=lightgray](-3,-1)(6,-1)
\psline[linecolor=lightgray](-3,0)(6,0)
\psline[linecolor=lightgray](-3,1)(6,1)
\psline[linecolor=lightgray](-3,2)(6,2)
\psline[linecolor=lightgray](-3,3)(6,3)
\psline[linecolor=lightgray](-3,4)(6,4)
\psline[linecolor=lightgray](-3,5)(6,5)
\psline[linecolor=lightgray](-3,6)(6,6)
\psline[linecolor=black,arrows=->](-3.7,0)(6.7,0)
\psline[linecolor=black,arrows=->](0,-1.7)(0,6.7)
\rput(1.5,5.7){ {\small\gray l, o, $k=3$, $m=2$, $n=2$}}
\rput(1.5,-0.7){ {\small\gray $K=5$, $M=-3$, $N=6$}}
\psset{linewidth=0.5pt,linecolor=gray,linestyle=solid,fillstyle=none}
\psline(-1.7,-1.7)(6.7,6.7)
\psline[linecolor=magenta,linestyle=dashed](3.3,-1.7)(6.7,1.7)
\pscircle[linecolor=black,linewidth=0.5pt,fillstyle=solid,fillcolor=red](0,0){0.15}
\pscircle[linecolor=black,linewidth=0.5pt,fillstyle=solid,fillcolor=red](1,1){0.15}
\pscircle[linecolor=black,linewidth=0.5pt,fillstyle=solid,fillcolor=red](2,2){0.15}
\pscircle[linecolor=black,linewidth=0.5pt,fillstyle=solid,fillcolor=red](3,3){0.15}
\pscircle[linecolor=black,linewidth=0.5pt,fillstyle=solid,fillcolor=red](4,4){0.15}
\pscircle[linecolor=black,linewidth=0.5pt,fillstyle=solid,fillcolor=red](5,5){0.15}
\psarc[linecolor=black,linewidth=0.5pt,fillstyle=solid,fillcolor=red](1,1){0.15}{-45.0}{135.0}
\psarc[linecolor=black,linewidth=0.5pt,fillstyle=solid,fillcolor=blue](1,1){0.15}{135.0}{315.0}
\pscircle[linecolor=black,linewidth=0.5pt](1,1){0.15}
\psarc[linecolor=black,linewidth=0.5pt,fillstyle=solid,fillcolor=red](0,0){0.15}{-45.0}{135.0}
\psarc[linecolor=black,linewidth=0.5pt,fillstyle=solid,fillcolor=blue](0,0){0.15}{135.0}{315.0}
\pscircle[linecolor=black,linewidth=0.5pt](0,0){0.15}
\pscircle[linecolor=black,linewidth=0.5pt,fillstyle=solid,fillcolor=blue](3,0){0.15}
\pscircle[linecolor=black,linewidth=0.5pt,fillstyle=solid,fillcolor=blue](2,1){0.15}
\pscircle[linecolor=black,linewidth=0.5pt,fillstyle=solid,fillcolor=blue](-1,1){0.15}
\pscircle[linecolor=black,linewidth=0.5pt,fillstyle=solid,fillcolor=blue](-2,2){0.15}
\end{pspicture}\hfil
\psset{unit=0.5cm}
\begin{pspicture}(-1.9,-3.9)(6.9,6.9)
\pspolygon[linecolor=white,fillstyle=solid,fillcolor=backgroundgray,linearc=0.3](-1.9,-3.9)(6.9,-3.9)(6.9,6.9)(-1.9,6.9)

\psset{linewidth=0.5pt,linecolor=gray,linestyle=solid,fillstyle=none}
\psline[linecolor=lightgray](-1,-3)(-1,6)
\psline[linecolor=lightgray](0,-3)(0,6)
\psline[linecolor=lightgray](1,-3)(1,6)
\psline[linecolor=lightgray](2,-3)(2,6)
\psline[linecolor=lightgray](3,-3)(3,6)
\psline[linecolor=lightgray](4,-3)(4,6)
\psline[linecolor=lightgray](5,-3)(5,6)
\psline[linecolor=lightgray](6,-3)(6,6)
\psline[linecolor=lightgray](-1,-3)(6,-3)
\psline[linecolor=lightgray](-1,-2)(6,-2)
\psline[linecolor=lightgray](-1,-1)(6,-1)
\psline[linecolor=lightgray](-1,0)(6,0)
\psline[linecolor=lightgray](-1,1)(6,1)
\psline[linecolor=lightgray](-1,2)(6,2)
\psline[linecolor=lightgray](-1,3)(6,3)
\psline[linecolor=lightgray](-1,4)(6,4)
\psline[linecolor=lightgray](-1,5)(6,5)
\psline[linecolor=lightgray](-1,6)(6,6)
\psline[linecolor=black,arrows=->](-1.7,0)(6.7,0)
\psline[linecolor=black,arrows=->](0,-3.7)(0,6.7)
\rput(2.5,5.7){ {\small\gray l, o, $k=3$, $m=2$, $n=2$}}
\rput(2.5,-2.7){ {\small\gray $K=5$, $M=-3$, $N=6$}}
\psset{linewidth=0.5pt,linecolor=gray,linestyle=solid,fillstyle=none}
\psline(-1.7,-1.7)(6.7,6.7)
\psline[linecolor=magenta,linestyle=dashed](1.3,-3.7)(6.7,1.7)
\pscircle[linecolor=black,linewidth=0.5pt,fillstyle=solid,fillcolor=red](0,0){0.15}
\pscircle[linecolor=black,linewidth=0.5pt,fillstyle=solid,fillcolor=red](1,1){0.15}
\pscircle[linecolor=black,linewidth=0.5pt,fillstyle=solid,fillcolor=red](2,2){0.15}
\pscircle[linecolor=black,linewidth=0.5pt,fillstyle=solid,fillcolor=red](3,3){0.15}
\pscircle[linecolor=black,linewidth=0.5pt,fillstyle=solid,fillcolor=red](4,4){0.15}
\pscircle[linecolor=black,linewidth=0.5pt,fillstyle=solid,fillcolor=red](5,5){0.15}
\psarc[linecolor=black,linewidth=0.5pt,fillstyle=solid,fillcolor=red](1,1){0.15}{-45.0}{135.0}
\psarc[linecolor=black,linewidth=0.5pt,fillstyle=solid,fillcolor=blue](1,1){0.15}{135.0}{315.0}
\pscircle[linecolor=black,linewidth=0.5pt](1,1){0.15}
\psarc[linecolor=black,linewidth=0.5pt,fillstyle=solid,fillcolor=red](0,0){0.15}{-45.0}{135.0}
\psarc[linecolor=black,linewidth=0.5pt,fillstyle=solid,fillcolor=blue](0,0){0.15}{135.0}{315.0}
\pscircle[linecolor=black,linewidth=0.5pt](0,0){0.15}
\pscircle[linecolor=black,linewidth=0.5pt,fillstyle=solid,fillcolor=blue](3,0){0.15}
\pscircle[linecolor=black,linewidth=0.5pt,fillstyle=solid,fillcolor=blue](2,1){0.15}
\pscircle[linecolor=black,linewidth=0.5pt,fillstyle=solid,fillcolor=green](1,-1){0.15}
\pscircle[linecolor=black,linewidth=0.5pt,fillstyle=solid,fillcolor=green](2,-2){0.15}
\pscircle[linecolor=black,linewidth=0.5pt,fillstyle=solid,fillcolor=white](4,-1){0.15}
\pscircle[linecolor=black,linewidth=0.5pt,fillstyle=solid,fillcolor=white](5,0){0.15}
\pspolygon[fillstyle=solid,linecolor=peachpuff,fillcolor=peachpuff,opacity=0.4,strokeopacity=0.6](4,-1)(5,0)(5,5)(1.5,1.5)
\end{pspicture}\end{center}
Note that in the right picture above (showing the folded situation), a certain
trapezoid is marked in colour, and
additional white points on the forbidden line points are shown. We call this the \EM{essential
region}: its significance and the meaning of the white points will soon become clear.

\secC{Reflection for the \rhs: The reflected situation}
If we \EM{reflect} on $d$ the forbidden line
$y=x-K$ and \EM{all} \nilp s for the \rhs, then we call the result of such reflection the \EM{reflected
situation}. 

Again, we illustrate this for the odd identity with  parameters $k=3$, $m=2$ and $n=2$
(i.e., $K=2k-1=5$, $M=1-k+m=0$, and $N=n = 2$).
``In terms of initial and terminal points'', the transition from the \EM{original}
situation to the \EM{reflected} situation looks like this:
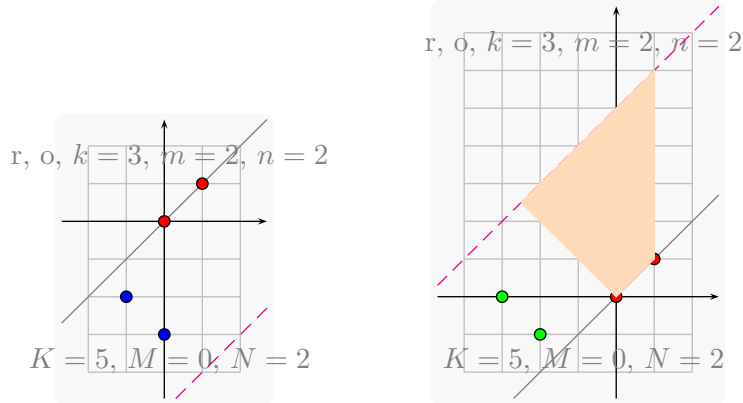
\begin{center}
\psset{unit=0.5cm}
\begin{pspicture}(-2.9,-4.9)(2.9,2.9)
\pspolygon[linecolor=white,fillstyle=solid,fillcolor=backgroundgray,linearc=0.3](-2.9,-4.9)(2.9,-4.9)(2.9,2.9)(-2.9,2.9)

\psset{linewidth=0.5pt,linecolor=gray,linestyle=solid,fillstyle=none}
\psline[linecolor=lightgray](-2,-4)(-2,2)
\psline[linecolor=lightgray](-1,-4)(-1,2)
\psline[linecolor=lightgray](0,-4)(0,2)
\psline[linecolor=lightgray](1,-4)(1,2)
\psline[linecolor=lightgray](2,-4)(2,2)
\psline[linecolor=lightgray](-2,-4)(2,-4)
\psline[linecolor=lightgray](-2,-3)(2,-3)
\psline[linecolor=lightgray](-2,-2)(2,-2)
\psline[linecolor=lightgray](-2,-1)(2,-1)
\psline[linecolor=lightgray](-2,0)(2,0)
\psline[linecolor=lightgray](-2,1)(2,1)
\psline[linecolor=lightgray](-2,2)(2,2)
\psline[linecolor=black,arrows=->](-2.7,0)(2.7,0)
\psline[linecolor=black,arrows=->](0,-4.7)(0,2.7)
\rput(0,1.7){ {\small\gray r, o, $k=3$, $m=2$, $n=2$}}
\rput(0,-3.7){ {\small\gray $K=5$, $M=0$, $N=2$}}
\psset{linewidth=0.5pt,linecolor=gray,linestyle=solid,fillstyle=none}
\psline(-2.7,-2.7)(2.7,2.7)
\psline[linecolor=magenta,linestyle=dashed](0.3,-4.7)(2.7,-2.3)
\pscircle[linecolor=black,linewidth=0.5pt,fillstyle=solid,fillcolor=red](0,0){0.15}
\pscircle[linecolor=black,linewidth=0.5pt,fillstyle=solid,fillcolor=red](1,1){0.15}
\pscircle[linecolor=black,linewidth=0.5pt,fillstyle=solid,fillcolor=blue](0,-3){0.15}
\pscircle[linecolor=black,linewidth=0.5pt,fillstyle=solid,fillcolor=blue](-1,-2){0.15}
\end{pspicture}\hfil
\psset{unit=0.5cm}
\begin{pspicture}(-4.9,-2.9)(2.9,7.9)
\pspolygon[linecolor=white,fillstyle=solid,fillcolor=backgroundgray,linearc=0.3](-4.9,-2.9)(2.9,-2.9)(2.9,7.9)(-4.9,7.9)

\psset{linewidth=0.5pt,linecolor=gray,linestyle=solid,fillstyle=none}
\psline[linecolor=lightgray](-4,-2)(-4,7)
\psline[linecolor=lightgray](-3,-2)(-3,7)
\psline[linecolor=lightgray](-2,-2)(-2,7)
\psline[linecolor=lightgray](-1,-2)(-1,7)
\psline[linecolor=lightgray](0,-2)(0,7)
\psline[linecolor=lightgray](1,-2)(1,7)
\psline[linecolor=lightgray](2,-2)(2,7)
\psline[linecolor=lightgray](-4,-2)(2,-2)
\psline[linecolor=lightgray](-4,-1)(2,-1)
\psline[linecolor=lightgray](-4,0)(2,0)
\psline[linecolor=lightgray](-4,1)(2,1)
\psline[linecolor=lightgray](-4,2)(2,2)
\psline[linecolor=lightgray](-4,3)(2,3)
\psline[linecolor=lightgray](-4,4)(2,4)
\psline[linecolor=lightgray](-4,5)(2,5)
\psline[linecolor=lightgray](-4,6)(2,6)
\psline[linecolor=lightgray](-4,7)(2,7)
\psline[linecolor=black,arrows=->](-4.7,0)(2.7,0)
\psline[linecolor=black,arrows=->](0,-2.7)(0,7.7)
\rput(-1,6.7){ {\small\gray r, o, $k=3$, $m=2$, $n=2$}}
\rput(-1,-1.7){ {\small\gray $K=5$, $M=0$, $N=2$}}
\psset{linewidth=0.5pt,linecolor=gray,linestyle=solid,fillstyle=none}
\psline(-2.7,-2.7)(2.7,2.7)
\psline[linecolor=magenta,linestyle=dashed](-4.7,0.3)(2.7,7.7)
\pscircle[linecolor=black,linewidth=0.5pt,fillstyle=solid,fillcolor=red](0,0){0.15}
\pscircle[linecolor=black,linewidth=0.5pt,fillstyle=solid,fillcolor=red](1,1){0.15}
\pscircle[linecolor=black,linewidth=0.5pt,fillstyle=solid,fillcolor=green](-3,0){0.15}
\pscircle[linecolor=black,linewidth=0.5pt,fillstyle=solid,fillcolor=green](-2,-1){0.15}
\pspolygon[fillstyle=solid,linecolor=peachpuff,fillcolor=peachpuff,opacity=0.4,strokeopacity=0.6](0,0)(1,1)(1,6)(-2.5,2.5)
\end{pspicture}\end{center}
Note that in the right picture above (showing the reflected situation), a certain
trapezoid is marked in colour again, which we also call the \EM{essential region}
(for the reflected situation).

\secC{Essential regions}
It is easy to see (recall the considerations and illustrations
in section~\ref{sec:even-odd}) that there is a \EM{translation} 
which transforms
\bit
\item the essential regions for the folded situation
\item to the essential region
for the reflected situation,
\eit
which maps bijectively 
\bit
\item green initial points to green initial points
\item and white additional points to terminal points,
\eit
and which \EM{swaps the roles} of the diagonal and the forbidden line
for the folded situation and the reflected situation (see the essential
regions in the above pictures).

This congruence of essential regions is crucial for the second step in our proof:
%
%
Clearly, the reflections explained above give \EM{sign--preserving
bijections} (since the permutation $\pi$ is not affected by the reflections)
from the ``original situation'' to the folded or reflected
situation, respectively, so instead of defining a sign--preserving bijection
between
\bit
\item the survivors for the \lhs\ which also ``survive'' the cancellation effected
by the \EM{second} sign--reversing involution $\psi$ (yet to be defined)
\item and the survivors for the \rhs,
\eit
we will define the sign--preserving bijection $\xi$ between
\bit
\item the ``twofold'' survivors for the \lhs\ \EM{in the folded situation}
\item and the survivors for the \rhs\ \EM{in the reflected situation}.
\eit
The construction of this bijection $\xi$ will involve \EM{only} the essential regions, leaving segments of green
paths \EM{outside} the essential regions \EM{unchanged}.

\secC{Folded overlays}
In the following illustration, the left picture shows \nilp s (for the \lhs),
and the right picture shows the \EM{folded} situation,
obtained by the reflection of all paths above the diagonal $d=\pas{y=x}$
on $d$:
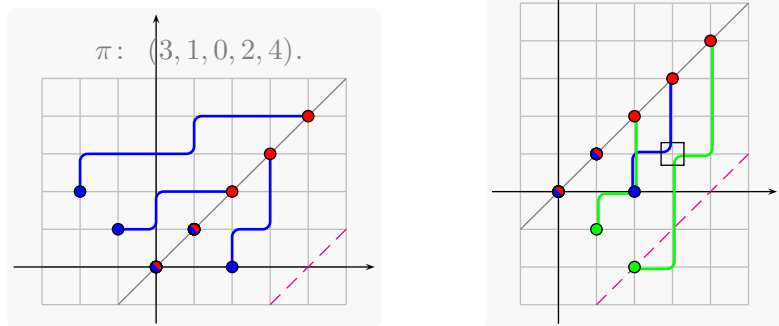
\begin{center}
\psset{unit=0.5cm}
\begin{pspicture}(-3.95,-1.7)(5.95,6.9)
\pspolygon[linecolor=white,fillstyle=solid,fillcolor=backgroundgray,linearc=0.3](-3.95,-1.7)(5.95,-1.7)(5.95,6.9)(-3.95,6.9)

\psset{linewidth=0.5pt,linecolor=gray,linestyle=solid,fillstyle=none}
\psline[linecolor=lightgray](-3,-1)(-3,5)
\psline[linecolor=lightgray](-2,-1)(-2,5)
\psline[linecolor=lightgray](-1,-1)(-1,5)
\psline[linecolor=lightgray](0,-1)(0,5)
\psline[linecolor=lightgray](1,-1)(1,5)
\psline[linecolor=lightgray](2,-1)(2,5)
\psline[linecolor=lightgray](3,-1)(3,5)
\psline[linecolor=lightgray](4,-1)(4,5)
\psline[linecolor=lightgray](5,-1)(5,5)
\psline[linecolor=lightgray](-3,-1)(5,-1)
\psline[linecolor=lightgray](-3,0)(5,0)
\psline[linecolor=lightgray](-3,1)(5,1)
\psline[linecolor=lightgray](-3,2)(5,2)
\psline[linecolor=lightgray](-3,3)(5,3)
\psline[linecolor=lightgray](-3,4)(5,4)
\psline[linecolor=lightgray](-3,5)(5,5)
\psline[linecolor=black,arrows=->](-3.75,0)(5.75,0)
\psline[linecolor=black,arrows=->](0,-1.5)(0,6.7)
\psline(-1,-1)(5,5)
\psline[linecolor=magenta,linestyle=dashed](3,-1)(5,1)
\psset{linewidth=1pt,linecolor=black,linestyle=solid,fillstyle=none}
\psline[linearc=0.2,linecolor=blue](2,0)(2,1)(3,1)(3,2)(3,3)
\psline[linearc=0.2,linecolor=blue](-2,2)(-2,3)(-1,3)(0,3)(1,3)(1,4)(2,4)(3,4)(4,4)
\psline[linearc=0.2,linecolor=blue](-1,1)(0,1)(0,2)(1,2)(2,2)
\pscircle[linecolor=black,linewidth=0.5pt,fillstyle=solid,fillcolor=red](0,0){0.15}
\pscircle[linecolor=black,linewidth=0.5pt,fillstyle=solid,fillcolor=red](1,1){0.15}
\pscircle[linecolor=black,linewidth=0.5pt,fillstyle=solid,fillcolor=red](2,2){0.15}
\pscircle[linecolor=black,linewidth=0.5pt,fillstyle=solid,fillcolor=red](3,3){0.15}
\pscircle[linecolor=black,linewidth=0.5pt,fillstyle=solid,fillcolor=red](4,4){0.15}
\pscircle[linecolor=black,linewidth=0.5pt,fillstyle=solid,fillcolor=blue](2,0){0.15}
\psarc[linecolor=black,linewidth=0.5pt,fillstyle=solid,fillcolor=red](1,1){0.15}{-45.0}{135.0}
\psarc[linecolor=black,linewidth=0.5pt,fillstyle=solid,fillcolor=blue](1,1){0.15}{135.0}{315.0}
\pscircle[linecolor=black,linewidth=0.5pt](1,1){0.15}
\psarc[linecolor=black,linewidth=0.5pt,fillstyle=solid,fillcolor=red](0,0){0.15}{-45.0}{135.0}
\psarc[linecolor=black,linewidth=0.5pt,fillstyle=solid,fillcolor=blue](0,0){0.15}{135.0}{315.0}
\pscircle[linecolor=black,linewidth=0.5pt](0,0){0.15}
\pscircle[linecolor=black,linewidth=0.5pt,fillstyle=solid,fillcolor=blue](-1,1){0.15}
\pscircle[linecolor=black,linewidth=0.5pt,fillstyle=solid,fillcolor=blue](-2,2){0.15}
\psarc[linecolor=black,linewidth=0.5pt,fillstyle=solid,fillcolor=red](0,0){0.15}{-45.0}{135.0}
\psarc[linecolor=black,linewidth=0.5pt,fillstyle=solid,fillcolor=blue](0,0){0.15}{135.0}{315.0}
\pscircle[linecolor=black,linewidth=0.5pt](0,0){0.15}
\psarc[linecolor=black,linewidth=0.5pt,fillstyle=solid,fillcolor=red](1,1){0.15}{-45.0}{135.0}
\psarc[linecolor=black,linewidth=0.5pt,fillstyle=solid,fillcolor=blue](1,1){0.15}{135.0}{315.0}
\pscircle[linecolor=black,linewidth=0.5pt](1,1){0.15}
\rput(1,5.7){ {\small\gray $\pi\fdg\pas{3,1,0,2,4}$.}}
\end{pspicture}\hfil
\psset{unit=0.5cm}
\begin{pspicture}(-1.95,-3.7)(5.95,6.2)
\pspolygon[linecolor=white,fillstyle=solid,fillcolor=backgroundgray,linearc=0.3](-1.95,-3.7)(5.95,-3.7)(5.95,6.2)(-1.95,6.2)

\psset{linewidth=0.5pt,linecolor=gray,linestyle=solid,fillstyle=none}
\psline[linecolor=lightgray](-1,-3)(-1,5)
\psline[linecolor=lightgray](0,-3)(0,5)
\psline[linecolor=lightgray](1,-3)(1,5)
\psline[linecolor=lightgray](2,-3)(2,5)
\psline[linecolor=lightgray](3,-3)(3,5)
\psline[linecolor=lightgray](4,-3)(4,5)
\psline[linecolor=lightgray](5,-3)(5,5)
\psline[linecolor=lightgray](-1,-3)(5,-3)
\psline[linecolor=lightgray](-1,-2)(5,-2)
\psline[linecolor=lightgray](-1,-1)(5,-1)
\psline[linecolor=lightgray](-1,0)(5,0)
\psline[linecolor=lightgray](-1,1)(5,1)
\psline[linecolor=lightgray](-1,2)(5,2)
\psline[linecolor=lightgray](-1,3)(5,3)
\psline[linecolor=lightgray](-1,4)(5,4)
\psline[linecolor=lightgray](-1,5)(5,5)
\psline[linecolor=black,arrows=->](-1.75,0)(5.75,0)
\psline[linecolor=black,arrows=->](0,-3.5)(0,6)
\psline(-1,-1)(5,5)
\psline[linecolor=magenta,linestyle=dashed](1,-3)(5,1)
\psset{linewidth=1pt,linecolor=black,linestyle=solid,fillstyle=none}
\psline[linearc=0.2,linecolor=blue](1.95,0.05)(1.95,1.05)(2.95,1.05)(2.95,2.05)(2.95,3.05)
\psline[linearc=0.2,linecolor=green](2.05,-2.05)(3.05,-2.05)(3.05,-1.05)(3.05,-0.05)(3.05,0.95)(4.05,0.95)(4.05,1.95)(4.05,2.95)(4.05,3.95)
\psline[linearc=0.2,linecolor=green](1.05,-1.05)(1.05,-0.05)(2.05,-0.05)(2.05,0.95)(2.05,1.95)
\pscircle[linecolor=black,linewidth=0.5pt,fillstyle=solid,fillcolor=red](0,0){0.15}
\pscircle[linecolor=black,linewidth=0.5pt,fillstyle=solid,fillcolor=red](1,1){0.15}
\pscircle[linecolor=black,linewidth=0.5pt,fillstyle=solid,fillcolor=red](2,2){0.15}
\pscircle[linecolor=black,linewidth=0.5pt,fillstyle=solid,fillcolor=red](3,3){0.15}
\pscircle[linecolor=black,linewidth=0.5pt,fillstyle=solid,fillcolor=red](4,4){0.15}
\pscircle[linecolor=black,linewidth=0.5pt,fillstyle=solid,fillcolor=blue](2,0){0.15}
\psarc[linecolor=black,linewidth=0.5pt,fillstyle=solid,fillcolor=red](1,1){0.15}{-45.0}{135.0}
\psarc[linecolor=black,linewidth=0.5pt,fillstyle=solid,fillcolor=blue](1,1){0.15}{135.0}{315.0}
\pscircle[linecolor=black,linewidth=0.5pt](1,1){0.15}
\psarc[linecolor=black,linewidth=0.5pt,fillstyle=solid,fillcolor=red](0,0){0.15}{-45.0}{135.0}
\psarc[linecolor=black,linewidth=0.5pt,fillstyle=solid,fillcolor=blue](0,0){0.15}{135.0}{315.0}
\pscircle[linecolor=black,linewidth=0.5pt](0,0){0.15}
\pscircle[linecolor=black,linewidth=0.5pt,fillstyle=solid,fillcolor=green](1,-1){0.15}
\pscircle[linecolor=black,linewidth=0.5pt,fillstyle=solid,fillcolor=green](2,-2){0.15}
\psarc[linecolor=black,linewidth=0.5pt,fillstyle=solid,fillcolor=red](0,0){0.15}{-45.0}{135.0}
\psarc[linecolor=black,linewidth=0.5pt,fillstyle=solid,fillcolor=blue](0,0){0.15}{135.0}{315.0}
\pscircle[linecolor=black,linewidth=0.5pt](0,0){0.15}
\psarc[linecolor=black,linewidth=0.5pt,fillstyle=solid,fillcolor=red](1,1){0.15}{-45.0}{135.0}
\psarc[linecolor=black,linewidth=0.5pt,fillstyle=solid,fillcolor=blue](1,1){0.15}{135.0}{315.0}
\pscircle[linecolor=black,linewidth=0.5pt](1,1){0.15}
\pspolygon[linecolor=black,linewidth=0.5pt](2.7,0.7)(3.3,0.7)(3.3,1.3)(2.7,1.3)
\end{pspicture}\end{center}
Observe that by this reflection, we
get an \EM{overlay} of green and blue paths, such that
\bit
\item all paths (blue or green) are restricted to the range $y\leq x$,
\item blue paths \EM{additionally} are restricted to the range $y>x-K$,
\item blue and green paths are \EM{nonintersecting}, i.e., there is no
	point of intersection of two blue paths or two green paths.
\eit
We shall call such overlay of green and blue paths a \EM{folded overlay}.
(As already mentioned, folded overlays are in \EM{sign--preserving bijection}
with the ``original'' \nilp s).

Note that a folded overlay appears as a directed subgraph of the lattice $\Z^2$
with green and blue edges, which
\bit
\item have the \EM{normal} direction (rightwards or upwards),
\item but might be traversed in the \EM{reversed} direction (leftwards or downwards),
\eit
and where double edges (which necessarily must be of opposite colours blue and green)
are allowed. We call vertices in this graph which are incident \EM{only} with a green (or blue,
respectively) edge a green (or blue, respectively) point, and we call points which
are incident with green \EM{and} blue edges \EM{bicoloured}. 

\secC{Points of bicoloured intersection and the \LGV--idea for folded overlays}
An \EM{overlay of green and blue paths} may have bicoloured points where some
\EM{blue} path intersects
some \EM{green} path: in the above right picture, there are three such \EM{points
of bicoloured intersection}, the \EM{maximal} (in lexicographic order) of them
is indicated by a small square. Clearly, 
\bit
\item every \EM{terminal} point and every initial point
	which is \EM{not bicoloured}
	is incident with \EM{precisely one} step which is either green or blue,
\item every \EM{bicoloured} point is incident with
	\bit
	\item precisely two blue and precisely two green steps \EM{if}
	it is not an initial point,
	\item precisely one blue and precisely one green step \EM{if}
	it is an initial point,
	\eit
\item  every other point is incident with \EM{precisely two} steps
	of the \EM{same} colour (either green or blue).
\eit
In the above pictures, we immediately see how a sign--reversing
involution resembling the \LGV--method might work --- just look at the pictures below,
where terminal segments of paths after the maximal point of bicoloured intersection
are \EM{exchanged}:
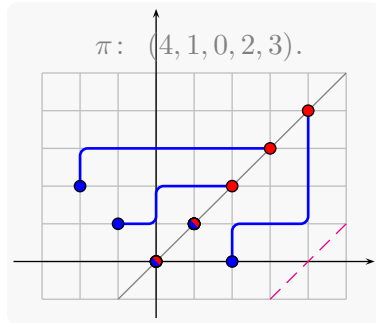
\begin{center}
\psset{unit=0.5cm}
\begin{pspicture}(-3.95,-1.7)(5.95,6.9)
\pspolygon[linecolor=white,fillstyle=solid,fillcolor=backgroundgray,linearc=0.3](-3.95,-1.7)(5.95,-1.7)(5.95,6.9)(-3.95,6.9)

\psset{linewidth=0.5pt,linecolor=gray,linestyle=solid,fillstyle=none}
\psline[linecolor=lightgray](-3,-1)(-3,5)
\psline[linecolor=lightgray](-2,-1)(-2,5)
\psline[linecolor=lightgray](-1,-1)(-1,5)
\psline[linecolor=lightgray](0,-1)(0,5)
\psline[linecolor=lightgray](1,-1)(1,5)
\psline[linecolor=lightgray](2,-1)(2,5)
\psline[linecolor=lightgray](3,-1)(3,5)
\psline[linecolor=lightgray](4,-1)(4,5)
\psline[linecolor=lightgray](5,-1)(5,5)
\psline[linecolor=lightgray](-3,-1)(5,-1)
\psline[linecolor=lightgray](-3,0)(5,0)
\psline[linecolor=lightgray](-3,1)(5,1)
\psline[linecolor=lightgray](-3,2)(5,2)
\psline[linecolor=lightgray](-3,3)(5,3)
\psline[linecolor=lightgray](-3,4)(5,4)
\psline[linecolor=lightgray](-3,5)(5,5)
\psline[linecolor=black,arrows=->](-3.75,0)(5.75,0)
\psline[linecolor=black,arrows=->](0,-1.5)(0,6.7)
\psline(-1,-1)(5,5)
\psline[linecolor=magenta,linestyle=dashed](3,-1)(5,1)
\psset{linewidth=1pt,linecolor=black,linestyle=solid,fillstyle=none}
\psline[linearc=0.2,linecolor=blue](2,0)(2,1)(3,1)(4,1)(4,2)(4,3)(4,4)
\psline[linearc=0.2,linecolor=blue](-2,2)(-2,3)(-1,3)(0,3)(1,3)(2,3)(3,3)
\psline[linearc=0.2,linecolor=blue](-1,1)(0,1)(0,2)(1,2)(2,2)
\pscircle[linecolor=black,linewidth=0.5pt,fillstyle=solid,fillcolor=red](0,0){0.15}
\pscircle[linecolor=black,linewidth=0.5pt,fillstyle=solid,fillcolor=red](1,1){0.15}
\pscircle[linecolor=black,linewidth=0.5pt,fillstyle=solid,fillcolor=red](2,2){0.15}
\pscircle[linecolor=black,linewidth=0.5pt,fillstyle=solid,fillcolor=red](3,3){0.15}
\pscircle[linecolor=black,linewidth=0.5pt,fillstyle=solid,fillcolor=red](4,4){0.15}
\pscircle[linecolor=black,linewidth=0.5pt,fillstyle=solid,fillcolor=blue](2,0){0.15}
\psarc[linecolor=black,linewidth=0.5pt,fillstyle=solid,fillcolor=red](1,1){0.15}{-45.0}{135.0}
\psarc[linecolor=black,linewidth=0.5pt,fillstyle=solid,fillcolor=blue](1,1){0.15}{135.0}{315.0}
\pscircle[linecolor=black,linewidth=0.5pt](1,1){0.15}
\psarc[linecolor=black,linewidth=0.5pt,fillstyle=solid,fillcolor=red](0,0){0.15}{-45.0}{135.0}
\psarc[linecolor=black,linewidth=0.5pt,fillstyle=solid,fillcolor=blue](0,0){0.15}{135.0}{315.0}
\pscircle[linecolor=black,linewidth=0.5pt](0,0){0.15}
\pscircle[linecolor=black,linewidth=0.5pt,fillstyle=solid,fillcolor=blue](-1,1){0.15}
\pscircle[linecolor=black,linewidth=0.5pt,fillstyle=solid,fillcolor=blue](-2,2){0.15}
\psarc[linecolor=black,linewidth=0.5pt,fillstyle=solid,fillcolor=red](0,0){0.15}{-45.0}{135.0}
\psarc[linecolor=black,linewidth=0.5pt,fillstyle=solid,fillcolor=blue](0,0){0.15}{135.0}{315.0}
\pscircle[linecolor=black,linewidth=0.5pt](0,0){0.15}
\psarc[linecolor=black,linewidth=0.5pt,fillstyle=solid,fillcolor=red](1,1){0.15}{-45.0}{135.0}
\psarc[linecolor=black,linewidth=0.5pt,fillstyle=solid,fillcolor=blue](1,1){0.15}{135.0}{315.0}
\pscircle[linecolor=black,linewidth=0.5pt](1,1){0.15}
\rput(1,5.7){ {\small\gray $\pi\fdg\pas{4,1,0,2,3}$.}}
\end{pspicture}\hfil
\psset{unit=0.5cm}
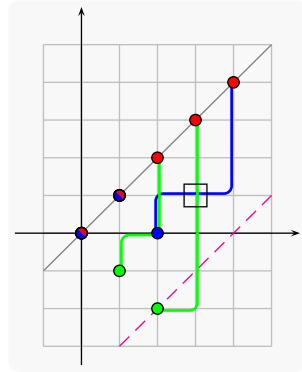
\begin{pspicture}(-1.95,-3.7)(5.95,6.2)
\pspolygon[linecolor=white,fillstyle=solid,fillcolor=backgroundgray,linearc=0.3](-1.95,-3.7)(5.95,-3.7)(5.95,6.2)(-1.95,6.2)

\psset{linewidth=0.5pt,linecolor=gray,linestyle=solid,fillstyle=none}
\psline[linecolor=lightgray](-1,-3)(-1,5)
\psline[linecolor=lightgray](0,-3)(0,5)
\psline[linecolor=lightgray](1,-3)(1,5)
\psline[linecolor=lightgray](2,-3)(2,5)
\psline[linecolor=lightgray](3,-3)(3,5)
\psline[linecolor=lightgray](4,-3)(4,5)
\psline[linecolor=lightgray](5,-3)(5,5)
\psline[linecolor=lightgray](-1,-3)(5,-3)
\psline[linecolor=lightgray](-1,-2)(5,-2)
\psline[linecolor=lightgray](-1,-1)(5,-1)
\psline[linecolor=lightgray](-1,0)(5,0)
\psline[linecolor=lightgray](-1,1)(5,1)
\psline[linecolor=lightgray](-1,2)(5,2)
\psline[linecolor=lightgray](-1,3)(5,3)
\psline[linecolor=lightgray](-1,4)(5,4)
\psline[linecolor=lightgray](-1,5)(5,5)
\psline[linecolor=black,arrows=->](-1.75,0)(5.75,0)
\psline[linecolor=black,arrows=->](0,-3.5)(0,6)
\psline(-1,-1)(5,5)
\psline[linecolor=magenta,linestyle=dashed](1,-3)(5,1)
\psset{linewidth=1pt,linecolor=black,linestyle=solid,fillstyle=none}
\psline[linearc=0.2,linecolor=blue](1.95,0.05)(1.95,1.05)(2.95,1.05)(3.95,1.05)(3.95,2.05)(3.95,3.05)(3.95,4.05)
\psline[linearc=0.2,linecolor=green](2.05,-2.05)(3.05,-2.05)(3.05,-1.05)(3.05,-0.05)(3.05,0.95)(3.05,1.95)(3.05,2.95)
\psline[linearc=0.2,linecolor=green](1.05,-1.05)(1.05,-0.05)(2.05,-0.05)(2.05,0.95)(2.05,1.95)
\pscircle[linecolor=black,linewidth=0.5pt,fillstyle=solid,fillcolor=red](0,0){0.15}
\pscircle[linecolor=black,linewidth=0.5pt,fillstyle=solid,fillcolor=red](1,1){0.15}
\pscircle[linecolor=black,linewidth=0.5pt,fillstyle=solid,fillcolor=red](2,2){0.15}
\pscircle[linecolor=black,linewidth=0.5pt,fillstyle=solid,fillcolor=red](3,3){0.15}
\pscircle[linecolor=black,linewidth=0.5pt,fillstyle=solid,fillcolor=red](4,4){0.15}
\pscircle[linecolor=black,linewidth=0.5pt,fillstyle=solid,fillcolor=blue](2,0){0.15}
\psarc[linecolor=black,linewidth=0.5pt,fillstyle=solid,fillcolor=red](1,1){0.15}{-45.0}{135.0}
\psarc[linecolor=black,linewidth=0.5pt,fillstyle=solid,fillcolor=blue](1,1){0.15}{135.0}{315.0}
\pscircle[linecolor=black,linewidth=0.5pt](1,1){0.15}
\psarc[linecolor=black,linewidth=0.5pt,fillstyle=solid,fillcolor=red](0,0){0.15}{-45.0}{135.0}
\psarc[linecolor=black,linewidth=0.5pt,fillstyle=solid,fillcolor=blue](0,0){0.15}{135.0}{315.0}
\pscircle[linecolor=black,linewidth=0.5pt](0,0){0.15}
\pscircle[linecolor=black,linewidth=0.5pt,fillstyle=solid,fillcolor=green](1,-1){0.15}
\pscircle[linecolor=black,linewidth=0.5pt,fillstyle=solid,fillcolor=green](2,-2){0.15}
\psarc[linecolor=black,linewidth=0.5pt,fillstyle=solid,fillcolor=red](0,0){0.15}{-45.0}{135.0}
\psarc[linecolor=black,linewidth=0.5pt,fillstyle=solid,fillcolor=blue](0,0){0.15}{135.0}{315.0}
\pscircle[linecolor=black,linewidth=0.5pt](0,0){0.15}
\psarc[linecolor=black,linewidth=0.5pt,fillstyle=solid,fillcolor=red](1,1){0.15}{-45.0}{135.0}
\psarc[linecolor=black,linewidth=0.5pt,fillstyle=solid,fillcolor=blue](1,1){0.15}{135.0}{315.0}
\pscircle[linecolor=black,linewidth=0.5pt](1,1){0.15}
\pspolygon[linecolor=black,linewidth=0.5pt](2.7,0.7)(3.3,0.7)(3.3,1.3)(2.7,1.3)
\end{pspicture}\end{center}
However, we need a slightly more complicated construction.

\secC{The significance of essential regions}
Clearly, a \EM{point of bicoloured intersection} 
in some
\EM{folded overlay} can only be located
\bit
\item \EM{weakly above} the forbidden line (since there are no blue steps strictly below)
\item and \EM{weakly above} the line containing the blue initial points (since there are no blue steps strictly below)
\item and \EM{weakly below} the line $\pas{y=x}$ (since there are no points above)
\item and \EM{weakly to the left} of the line $x = M-1$ (since there are no points to the right).
\eit
Note that the region just described is precisely the \EM{essential region}
for the folded situation: in the following, we will mostly focus our attention on this essential region.

\secC{Bicoloured connections and the sign--reversing involution}
For each terminal point $B_j$ in a \EM{folded overlay}, we construct the \EM{bicoloured path}
starting at $B_j$ as follows: we start with 
\bit
\item $P=B_j$,
\item the \EM{reversed} direction as \EM{current direction}
\item and the \EM{unique} colour
of the path ending in $B_j$ as \EM{current colour}.
\eit
As long as this is possible, we traverse the edge \EM{incident with $P$}
\bit
\item of the current colour
\item in the current direction,
\eit
to its other endpoint $Q$ and set $P=Q$; and if $Q$ is a point of bicoloured intersection,
we \EM{swap the current colour (blue/green) and the current direction (normal/reversed)}.

It is easy to see that this construction must end in an initial or terminal point $Q$
where we cannot continue (i.e., where there is no incident edge of the current colour and
current direction): we call such bicoloured path
from $B_j$ to $Q$ a \EM{bicoloured connection of $B_j$ and $Q$}.

The following pictures illustrate this construction: in the left picture, we start at the
green terminal point labeled $Q$, and in the right picture, we start at the
green terminal point labeled $P$.
\cps{bicoloured_connection_example}
(Observe that there is also a shorter bicoloured connection,
indicated in the pictures by the crosshatched area.)

Now we are in the position to define the involution $\psi$ on the set of all
\EM{overlays of green and blue paths} in the folded situation of the \lhs:
Let $o$ be an overlay of green and blue paths: we call a
\EM{bicoloured connection} $c$ in $o$
\bit
\item which \EM{connects two different} terminal points $B_{a} \neq B_{b}$,
\item and whose \EM{green} segments \EM{never intersect the forbidden line}
\eit
an \EM{involutive bicoloured connection}. (Note that $B_{a}$ and $B_{b}$
must have opposite colours, and that an involutive bicoloured connection
never leaves the essential region.)

A moment's thought shows that \EM{swapping the colours} (green to blue and vice versa) of
\EM{all} edges belonging to an involutive bicoloured connection gives another overlay of
green and blue paths $o^\prime$:
\bit
\item By construction, the blue paths \EM{and} the green paths in $o^\prime$ are \EM{nonintersecting},
\item and the blue paths in $o^\prime$ do not intersect the forbidden line.
\eit
(The right picture above is obtained by this swapping of colours, and vice versa.)

So for all folded overlays $o$
\bit
\item which do not contain any involutive bicoloured connection, we simply set $\psi\of o = o$,
\item which contain involutive bicoloured connections, we choose
the maximal terminal point (in lexicographic order) for which an involutive bicoloured connection $c$
exists and define $\psi\of o$ as the  folded overlay $o^\prime$ obtained by the swapping of
colours in $c$ (as described above).
\eit

In the following pictures,
the possible courses of bicoloured connections through a point of
bicoloured intersection are indicated by small arrows:
\bit
\item They show ``essentially all'' possible situations at \EM{points of intersections} in the network of blue and green
paths (modulo swap of colours blue and green, or directions horizontal and vertical) 
\item and make clear that each
\EM{bicoloured} edge would form a (trivial) bicoloured connection of the
vertices it is incident with (so bicoloured edges never belong to
bicoloured connections of initial or terminal points).
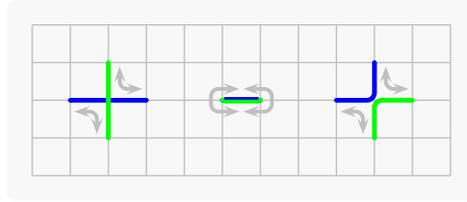
\eit
\cps{the_turns}
A moment's thought shows that swapping colours in a bicoloured connection 
\bit
\item does not change any \EM{other} bicoloured connection,
\item does not destroy any existing bicoloured connection,
\item and does not introduce new bicoloured connections:
\eit
So the mapping $\psi$ is, in fact, an \EM{involution}.

Moreover, bicoloured connections may intersect, but (by construction) can \EM{never cross}
(neither ``itself'' nor ``one another''). 
Therefore, a bicoloured connection of \EM{terminal points} $B_a$ and $B_b$
forms an ``impenetrable barrier'' for other bicoloured connections
(see the above pictures, where the area enclosed by this ``barrier''
is coloured). But this implies that the number of terminal points \EM{between}
	$B_a$ and $B_b$ must be \EM{even}, so for $o\neq\psi\of o$, we have 
	$\signum{\psi\of o} = -\signum o$, since
	the swapping of colours amounts to swapping 
	
	\bit
	\item a zero at position $a$
	\item and a one at position $b$
	\eit
	(or vice versa) with \EM{odd distance} $\abs{a-b}$
	in the $01$--code corresponding to $o$ (see the considerations in section
	\ref{sec:rhs-sign}).
	 Therefore, the mapping $\psi$ is, in fact,
	\EM{sign--reversing}.

So the determinant
on the \lhs\ in \eqref{eq:cigler-even} and \eqref{eq:cigler-odd} appears
as the generating function
of the fixed points of $\psi$, who ``survive'' the cancellation effected by $\psi$. Clearly,
these fixed points are the
folded overlays which \EM{do not} contain any involutive bicoloured connection:
We shall
call them
\EM{folded survivors}.

The next step in our proof is the exhibition of a bijection $\xi$ between
the
\bit
\item the folded survivors in the \lhs\ (in the folded situation)
\item and \EM{all} (``normal'') survivors in the \rhs\ (in the reflected situation)
\eit
of identities \eqref{eq:cigler-even} and \eqref{eq:cigler-odd}.


\secB{Step 2: Exhibit the bijection $\xi$}
The bijection $\xi$ we shall exhibit in this section between
\bit
\item \EM{folded survivors} (i.e, \nilp s from the \EM{folded situation} of the
    \lhs\ which ``survive'' the
	cancellation by the sign--preserving involution $\psi$),
\item and (``normal'') survivors (i.e., \nilp s) from the \EM{reflected} \rhs
\eit
will 
leave \EM{unchanged} all segments of paths \EM{outside the essential
regions} and employ a ``natural transformation'' from
\bit
\item the restriction to the essential region of folded survivors, which have
	\bit
	\item terminal points on the diagonal $y=x$
	\item and the forbidden line for blue paths $y=x-K$,
	\eit
\item and the restriction to the essential region of (``normal'') survivors, which have
	\bit
	\item terminal points on the diagonal (corresponding by a translation to the forbidden
		line)
	\item and the \EM{reflected} forbidden line $y=x+K$ (corresponding by a translation to
		the diagonal).
	\eit
\eit
In order to conceive this ``natural transformation'' inside the essential
region, we will show that a folded overlay is a \EM{folded survivor} if
and only if it has a certain \EM{simple structure}.
This will be the most strenuous part of this section: once we have achieved
this, the ``natural transformation'' will basically be obtained by
``inspection of pictures''.

\secC{Examples of folded survivors}
By definition, a \EM{folded survivor} is a folded overlay which \EM{does not} have 
an \EM{involutive bicoloured connection}, i.e.,
\bit
\item either there is \EM{no point of intersection at all}
\item or each bicoloured connection of terminal points contains some
	green segment which intersects the forbidden line $y=x-k$.
\eit
The following pictures illustrate this:
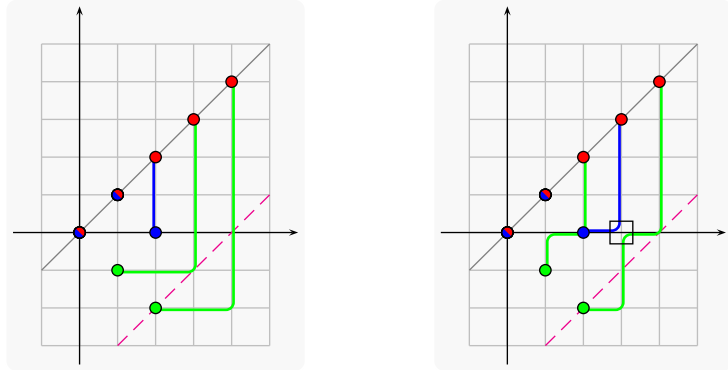
\begin{center}
\psset{unit=0.5cm}
\begin{pspicture}(-1.95,-3.7)(5.95,6.2)
\pspolygon[linecolor=white,fillstyle=solid,fillcolor=backgroundgray,linearc=0.3](-1.95,-3.7)(5.95,-3.7)(5.95,6.2)(-1.95,6.2)

\psset{linewidth=0.5pt,linecolor=gray,linestyle=solid,fillstyle=none}
\psline[linecolor=lightgray](-1,-3)(-1,5)
\psline[linecolor=lightgray](0,-3)(0,5)
\psline[linecolor=lightgray](1,-3)(1,5)
\psline[linecolor=lightgray](2,-3)(2,5)
\psline[linecolor=lightgray](3,-3)(3,5)
\psline[linecolor=lightgray](4,-3)(4,5)
\psline[linecolor=lightgray](5,-3)(5,5)
\psline[linecolor=lightgray](-1,-3)(5,-3)
\psline[linecolor=lightgray](-1,-2)(5,-2)
\psline[linecolor=lightgray](-1,-1)(5,-1)
\psline[linecolor=lightgray](-1,0)(5,0)
\psline[linecolor=lightgray](-1,1)(5,1)
\psline[linecolor=lightgray](-1,2)(5,2)
\psline[linecolor=lightgray](-1,3)(5,3)
\psline[linecolor=lightgray](-1,4)(5,4)
\psline[linecolor=lightgray](-1,5)(5,5)
\psline[linecolor=black,arrows=->](-1.75,0)(5.75,0)
\psline[linecolor=black,arrows=->](0,-3.5)(0,6)
\psline(-1,-1)(5,5)
\psline[linecolor=magenta,linestyle=dashed](1,-3)(5,1)
\psset{linewidth=1pt,linecolor=black,linestyle=solid,fillstyle=none}
\psline[linearc=0.2,linecolor=blue](1.95,0.05)(1.95,1.05)(1.95,2.05)
\psline[linearc=0.2,linecolor=green](2.05,-2.05)(3.05,-2.05)(4.05,-2.05)(4.05,-1.05)(4.05,-0.05)(4.05,0.95)(4.05,1.95)(4.05,2.95)(4.05,3.95)
\psline[linearc=0.2,linecolor=green](1.05,-1.05)(2.05,-1.05)(3.05,-1.05)(3.05,-0.05)(3.05,0.95)(3.05,1.95)(3.05,2.95)
\pscircle[linecolor=black,linewidth=0.5pt,fillstyle=solid,fillcolor=red](0,0){0.15}
\pscircle[linecolor=black,linewidth=0.5pt,fillstyle=solid,fillcolor=red](1,1){0.15}
\pscircle[linecolor=black,linewidth=0.5pt,fillstyle=solid,fillcolor=red](2,2){0.15}
\pscircle[linecolor=black,linewidth=0.5pt,fillstyle=solid,fillcolor=red](3,3){0.15}
\pscircle[linecolor=black,linewidth=0.5pt,fillstyle=solid,fillcolor=red](4,4){0.15}
\pscircle[linecolor=black,linewidth=0.5pt,fillstyle=solid,fillcolor=blue](2,0){0.15}
\psarc[linecolor=black,linewidth=0.5pt,fillstyle=solid,fillcolor=red](1,1){0.15}{-45.0}{135.0}
\psarc[linecolor=black,linewidth=0.5pt,fillstyle=solid,fillcolor=blue](1,1){0.15}{135.0}{315.0}
\pscircle[linecolor=black,linewidth=0.5pt](1,1){0.15}
\psarc[linecolor=black,linewidth=0.5pt,fillstyle=solid,fillcolor=red](0,0){0.15}{-45.0}{135.0}
\psarc[linecolor=black,linewidth=0.5pt,fillstyle=solid,fillcolor=blue](0,0){0.15}{135.0}{315.0}
\pscircle[linecolor=black,linewidth=0.5pt](0,0){0.15}
\pscircle[linecolor=black,linewidth=0.5pt,fillstyle=solid,fillcolor=green](1,-1){0.15}
\pscircle[linecolor=black,linewidth=0.5pt,fillstyle=solid,fillcolor=green](2,-2){0.15}
\psarc[linecolor=black,linewidth=0.5pt,fillstyle=solid,fillcolor=red](0,0){0.15}{-45.0}{135.0}
\psarc[linecolor=black,linewidth=0.5pt,fillstyle=solid,fillcolor=blue](0,0){0.15}{135.0}{315.0}
\pscircle[linecolor=black,linewidth=0.5pt](0,0){0.15}
\psarc[linecolor=black,linewidth=0.5pt,fillstyle=solid,fillcolor=red](1,1){0.15}{-45.0}{135.0}
\psarc[linecolor=black,linewidth=0.5pt,fillstyle=solid,fillcolor=blue](1,1){0.15}{135.0}{315.0}
\pscircle[linecolor=black,linewidth=0.5pt](1,1){0.15}
\end{pspicture}\hfil
\psset{unit=0.5cm}
\begin{pspicture}(-1.95,-3.7)(5.95,6.2)
\pspolygon[linecolor=white,fillstyle=solid,fillcolor=backgroundgray,linearc=0.3](-1.95,-3.7)(5.95,-3.7)(5.95,6.2)(-1.95,6.2)

\psset{linewidth=0.5pt,linecolor=gray,linestyle=solid,fillstyle=none}
\psline[linecolor=lightgray](-1,-3)(-1,5)
\psline[linecolor=lightgray](0,-3)(0,5)
\psline[linecolor=lightgray](1,-3)(1,5)
\psline[linecolor=lightgray](2,-3)(2,5)
\psline[linecolor=lightgray](3,-3)(3,5)
\psline[linecolor=lightgray](4,-3)(4,5)
\psline[linecolor=lightgray](5,-3)(5,5)
\psline[linecolor=lightgray](-1,-3)(5,-3)
\psline[linecolor=lightgray](-1,-2)(5,-2)
\psline[linecolor=lightgray](-1,-1)(5,-1)
\psline[linecolor=lightgray](-1,0)(5,0)
\psline[linecolor=lightgray](-1,1)(5,1)
\psline[linecolor=lightgray](-1,2)(5,2)
\psline[linecolor=lightgray](-1,3)(5,3)
\psline[linecolor=lightgray](-1,4)(5,4)
\psline[linecolor=lightgray](-1,5)(5,5)
\psline[linecolor=black,arrows=->](-1.75,0)(5.75,0)
\psline[linecolor=black,arrows=->](0,-3.5)(0,6)
\psline(-1,-1)(5,5)
\psline[linecolor=magenta,linestyle=dashed](1,-3)(5,1)
\psset{linewidth=1pt,linecolor=black,linestyle=solid,fillstyle=none}
\psline[linearc=0.2,linecolor=blue](1.95,0.05)(2.95,0.05)(2.95,1.05)(2.95,2.05)(2.95,3.05)
\psline[linearc=0.2,linecolor=green](2.05,-2.05)(3.05,-2.05)(3.05,-1.05)(3.05,-0.05)(4.05,-0.05)(4.05,0.95)(4.05,1.95)(4.05,2.95)(4.05,3.95)
\psline[linearc=0.2,linecolor=green](1.05,-1.05)(1.05,-0.05)(2.05,-0.05)(2.05,0.95)(2.05,1.95)
\pscircle[linecolor=black,linewidth=0.5pt,fillstyle=solid,fillcolor=red](0,0){0.15}
\pscircle[linecolor=black,linewidth=0.5pt,fillstyle=solid,fillcolor=red](1,1){0.15}
\pscircle[linecolor=black,linewidth=0.5pt,fillstyle=solid,fillcolor=red](2,2){0.15}
\pscircle[linecolor=black,linewidth=0.5pt,fillstyle=solid,fillcolor=red](3,3){0.15}
\pscircle[linecolor=black,linewidth=0.5pt,fillstyle=solid,fillcolor=red](4,4){0.15}
\pscircle[linecolor=black,linewidth=0.5pt,fillstyle=solid,fillcolor=blue](2,0){0.15}
\psarc[linecolor=black,linewidth=0.5pt,fillstyle=solid,fillcolor=red](1,1){0.15}{-45.0}{135.0}
\psarc[linecolor=black,linewidth=0.5pt,fillstyle=solid,fillcolor=blue](1,1){0.15}{135.0}{315.0}
\pscircle[linecolor=black,linewidth=0.5pt](1,1){0.15}
\psarc[linecolor=black,linewidth=0.5pt,fillstyle=solid,fillcolor=red](0,0){0.15}{-45.0}{135.0}
\psarc[linecolor=black,linewidth=0.5pt,fillstyle=solid,fillcolor=blue](0,0){0.15}{135.0}{315.0}
\pscircle[linecolor=black,linewidth=0.5pt](0,0){0.15}
\pscircle[linecolor=black,linewidth=0.5pt,fillstyle=solid,fillcolor=green](1,-1){0.15}
\pscircle[linecolor=black,linewidth=0.5pt,fillstyle=solid,fillcolor=green](2,-2){0.15}
\psarc[linecolor=black,linewidth=0.5pt,fillstyle=solid,fillcolor=red](0,0){0.15}{-45.0}{135.0}
\psarc[linecolor=black,linewidth=0.5pt,fillstyle=solid,fillcolor=blue](0,0){0.15}{135.0}{315.0}
\pscircle[linecolor=black,linewidth=0.5pt](0,0){0.15}
\psarc[linecolor=black,linewidth=0.5pt,fillstyle=solid,fillcolor=red](1,1){0.15}{-45.0}{135.0}
\psarc[linecolor=black,linewidth=0.5pt,fillstyle=solid,fillcolor=blue](1,1){0.15}{135.0}{315.0}
\pscircle[linecolor=black,linewidth=0.5pt](1,1){0.15}
\pspolygon[linecolor=black,linewidth=0.5pt](2.7,-0.3)(3.3,-0.3)(3.3,0.3)(2.7,0.3)
\end{pspicture}\end{center}


It is hard to conceive the structure of folded survivors in small examples, so we shall
consider larger ones: in order to limit the size
of the corresponding pictures, we observe that for any folded overlay $\mathcal F$
there are constants $z\in\Z$ such that $\mathcal F$ does not have points
of bicoloured intersection $\pas{x,y}$ with $y>z$, so we will ``cut away''
the uninteresting part above such level $z$ in the following
graphical illustrations. 

The following picture shows a folded survivor 
in the odd case with parameters $k=7$, $m=4$ and $n=9$ (which gives
$K=13$, $M=-9$ and $N=19$ for the \lhs), where we have
``cut away'' the part above level $5$:
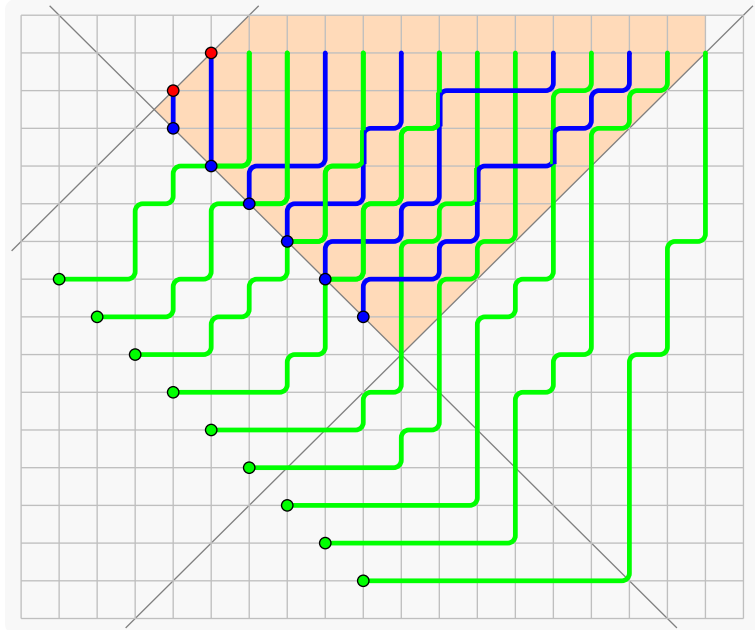
\cps{myex0_13_-9_19}
Here, the four two--faced points and only two of the terminal points are shown
(the others are ``cut away''), and the \EM{essential region}
(which, too is cut) is indicated as coloured area.

It is easy to verify ``by inspection'' that this folded
overlay has no involutive bicoloured connection. We shall
use this folded survivor as a running example in this section.


Since the transformation we want to present only involves the part
\EM{inside the essential region} (where points of bicoloured
intersections might be located), we shall omit the
green segments of paths \EM{outside} the essential region
in the following pictures.

\secC{The planar graph corresponding to a folded overlay}
Let $\mathcal F$ be some folded overlay: if we 
\bit
\item consider only points and steps
\EM{inside the essential region},
\item disregard the orientation (but not the colour!)
of all edges, 
\item and ``merge'' all
double (undirected) edges thus obtained into
\EM{bicoloured} single edges (having colours blue \EM{and} green; we shall
call the edges of \EM{unique} colour \EM{unicoloured}),
\eit
 then we obtain a simple
\EM{planar} graph. 
Note that all ``interior'' vertices of this graph have degree $2$, $3$ or $4$,
and ``loose ends'' (vertices of degree $1$) can only appear on the ``boundary'',
i.e., they are
\bit
\item either terminal points
\item or initial blue points
\item or green points on the forbidden line.
\eit 

We can ``close such loose ends'' by introducing \EM{artificial slanted edges}
connecting neighbouring ``loose ends'', with the following ``colouring rule'':
\bit
\item slanted edges on the line $\pas{y=x}$ have the same colour as their left point,
\item slanted edges on the line $\pas{y=-x-2M-K+2}$ (this is the line containing the blue intital points)
	are coloured green,
\item slanted edges on the forbidden line $\pas{y=x-K}$ are coloured blue;
\eit
see the following picture for an illustration, where we also added
blue and green terminal points (shifted downwards, to limit the height
of the picture) together with the entries in the corresponding 
$01$--code: recall that, by definition,
\bit
\item blue terminal points are reached \EM{from below},
\item green terminal points are reached \EM{from the left}
\eit
in the original (unfolded) situation (the rightmost isolated green
point on the forbidden
line indicates the $9$--th green path, all whose steps are completely
``cut away'') in this picture:
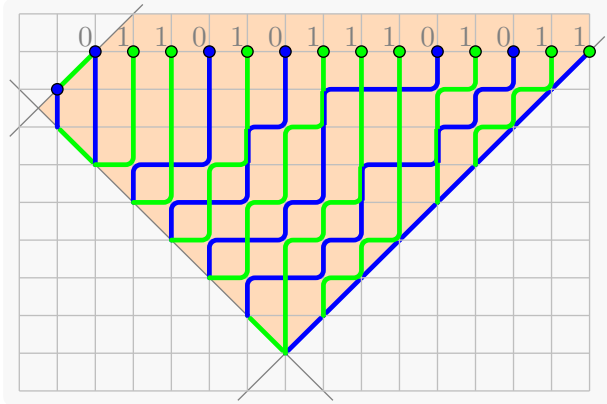
\cps{myex1a_13_-9_19}
The only reason for introducing these \EM{artificial slanted edges} is to avoid
the tedious distinction between ``partially open'' and ``closed areas'', which
now uniformly appear as \EM{finite faces} in the \EM{simple planar graph} thus
obtained: 
We call it the \EM{planar graph
corresponding to $\mathcal F$} and denote it by $\overline{\mathcal F}$ 


\secC{Faces of the planar graph corresponding to a folded overlay}
\label{sec:faces}
It is obvious that $\overline{\mathcal F}$ is \EM{connected}, and that
the \EM{boundary} $\partial f$ of each finite face $f$ of $\overline{\mathcal F}$ is a \EM{closed curve}
in $\R^2$, which we may traverse in counterclockwise orientation. But
it is \EM{not so obvious} that this closed curve must correspond to a \EM{circle}
in the graph--theoretical sense (i.e., passes through every point at most once).
 
When traversing (in counterclockwise orientation) the boundary $\partial f$ of
a finite face $f$ of of some \EM{general} planar subgraph $G$ (without vertices
of degree $1$) of the lattice $\Z^2$, there are edges
\bit
\item which are traversed rightwards or upwards: call those the \EM{black edges},
\item which are traversed leftwards or downwards: call those the \EM{red edges}.
\eit
Now observe that, in general, there can be precisely eight \EM{types of kinks} in
$\partial f$, as shown
in the following picture (where the crosshatched areas indicate the finite face $f$
whose boundary contains the kink):
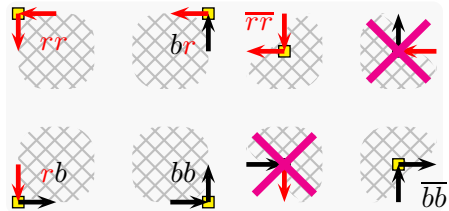
\def\rra{${\red rr}$}
\def\rrb{${\red \overline{rr}}$}
\def\bba{$bb$}
\def\bbb{$\overline{bb}$}
\def\br{$b{\red r}$}
\def\rb{${\red r}b$}
\cps{the_eight_kinks}
But a moment's thought shows that two of these kinks are \EM{impossible} if $f$ is a finite face of
the planar graph
$\overline{\mathcal F}$ corresponding to some folded overlay
${\mathcal F}$ (these are ``crossed out''
in the pictures above), so there are only
\bit
\item two types of kinks consisting of two red edges: denote them by \rra\ and \rrb,
\item two types of kinks consisting of two black edges: denote them by \bba\ and \bbb
\item and one type of kink where a red edge is followed by a black edge (denote this
	by \rb,
\item and one type of kink where a black edge is followed by a red edge (denote this
	by \br),
\eit
which can occur in the boundary $\partial f$ of some finite face $f$ of $\overline{\mathcal F}$. Note that
when traversing such boundary in counterclockwise orientation,
\bit
\item \rra\ may follow \rrb, and vice versa; and \rra\ may also be followed by \rb,
\item \bba\ may follow \bbb, and vice versa; and \bba\ may also be followed by \br,
\item \br\ \EM{can only be followed} by \rra,
\item \rb\ \EM{can only be followed} by \bba,
\eit
and there is no other possible succession of types of kinks.

Clearly, the boundary $\partial f$ must contain at least
\bit
\item one red path segment (consisting only of red edges)
\item and one  black path segment (consisting only of black edges)
\eit
(since otherwise $\partial f$ would be just a \EM{path} in $\Z^2$), and red and black segments must
be ``joined'' by kinks \br\ or \rb. Now consider the maximal kink \br\ (in lexicographic order) in
$\partial f$: the following picture (where the maximal kink \br\ is labeled $P$) makes clear that
there cannot be \EM{another} kink \br\ in $\partial f$.
\cps{kinks_on_face_boundary}
So the boundary $\partial f$ of every finite face of $\overline{\mathcal F}$ consists \EM{precisely} of 
\bit
\item one black segment (we shall call it the \EM{right boundary} of $f$)
\item and one red segment (we shall call it the \EM{left boundary} of $f$).
\eit

\secC{Vertical strips}
\label{sec:vstrips}
For our further considerations we will need a very simple special case:
We call a finite face $f$ of $\overline{\mathcal F}$  a \EM{vertical strip}
if its left and right boundaries both contain \EM{only one} horizontal or
slanted step, and we call the \EM{vertical segment} of the left (or right, respectively)
boundary of $f$ the \EM{left side} (or \EM{right side}, respectively) of $f$.

\def\BB{{\footnotesize ${\blue B} {\blue B}$}}
\def\BG{{\footnotesize ${\blue B} {\green G}$}}
\def\GB{{\footnotesize ${\green G} {\blue B}$}}
\def\GG{{\footnotesize ${\green G} {\green G}$}}

It is easy to see that for any vertical strip $f$ in $\overline{\mathcal F}$, 
\bit
\item the last (vertical) edge $e_v$ of the left (or right) side of $f$
	and the horizontal edge $e_h$ immediately following $e_v$ must have
	\EM{opposite} colours, whence \EM{both} $e_v$ and $e_h$ must be \EM{unicoloured},
\item and all edges of the left
(or right) side of $f$ must have the \EM{same} colour as $e_v$ (but some of them might
have \EM{both} colours);
\eit
see the following picture, which shows the four possible combinations
of the colours of the left and right side of vertical strips,
which we denote by \BB, \BG, \GB\ and \GG:
\cps{vertical_strip_colours}
Note that the horizontal edges in the pictures above could also be \EM{artificial slanted edges},
as explained above.

\secC{Adjacency of vertical strips}
\label{sec:adjacencies}
We say that a pair $\pas{f_l,f_r}$ of vertical strips
is \EM{adjacent} if the intersection
\bit
\item of the right side of $f_l$
\item with the left side of $f_r$
\eit
is \EM{not empty}, but \EM{does not} contain a horizontal edge, and we call $f_l$ the left
strip and $f_r$ the right strip of the adjacent pair.
The following pictures show all the \EM{relative positions} (``modulo
the swap of colours blue and green'') of left
and right strips of an adjacent pair $\pas{f_l,f_r}$, and
make clear that not all of these relative positions are possible
in the planar graph
$\overline{\mathcal F}$ corresponding to some folded overlay
$\mathcal F$: recall that
according to the above considerations, 
\bit
\item \EM{all} horizontal (or artificial slanted) edges are \EM{necessarily unicoloured} edges,
\item and
	\bit
	\item the upper horizontal (or artificial slanted) edge of the left strip $f_l$
		and the uppermost vertical edge of the right side of $f_l$
	\item as well as the the lower horizontal (or artificial slanted) edge of the right strip $f_r$
		and the lowest vertical edge of the left side of $f_r$
	\eit
are both \EM{unicoloured} and of \EM{different} colours.
\eit
In the following pictures, unicoloured horizontal edges whose colour is not already determined
(as just explained) are drawn in gray. However, in some situations, their colour is \EM{implied}
by the situation of the adjacency: this is shown by marking such edges by double exclamation marks
of the respective colour in the following pictures.

\secC{Impossible adjacencies}
A careful inspection reveals that the following situations of adjacencies are impossible
(and therefore are crossed out in the pictures):
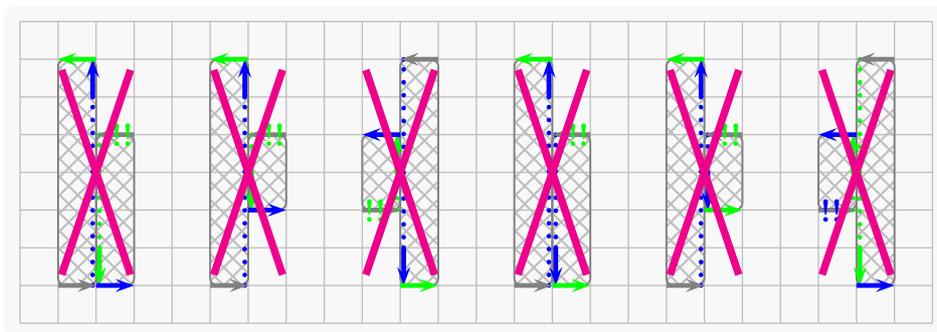
\cps{impossible_adjacencies}
To see this, just note that in a folded overlay
\bit
\item a point is never
incident with two outgoing (or incoming) steps of the \EM{same colour}
as is the case in the first three of the above pictures,
\item and a
path never starts (or ends) in the side of a vertical strip, as is
the case in the last three of the above pictures.
\eit
Also the following situations of adjacencies are impossible (note that the vertical sides
are of \EM{different} colour here):
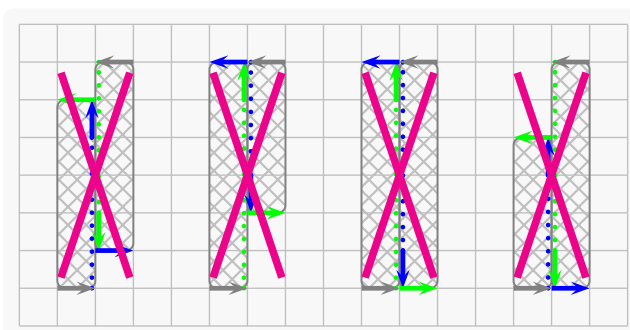
\cps{impossible_ascending_adjacencies}

\secC{Possible adjacencies: Ascending and descending}
However, the similar situations with vertical sides of the \EM{same}
colour \EM{are} possible: we say that an adjacent pair in any of
these relative positions is \EM{ascending}.
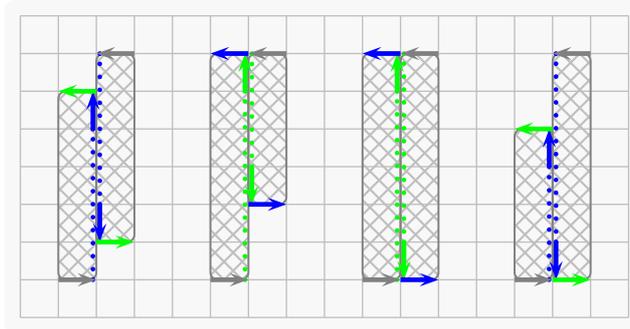
\cps{ascending_adjacencies-0}
Note for later reference that
\bit
\item the horizontal edges involved in an ascending adjacent pair of vertical
strips
\item together with their adjacent vertical edges
\eit
may be viewed as (segments of) \nilp s, see the following picture for an illustration, where
these segments are drawn in black:
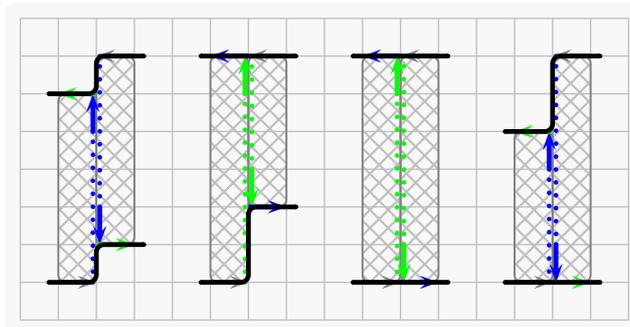
\cps{ascending_adjacencies}
Finally, also the following situations of adjacencies are all possible: we say that an adjacent pair in any of
the relative positions shown below is \EM{descending}.
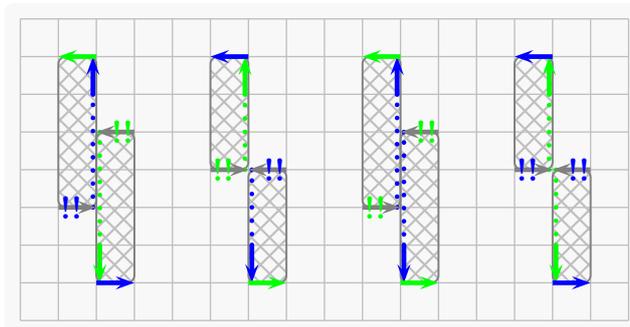
\cps{descending_adjacencies}

The following observation is an immediate consequence of the above considerations:
\bit
\item For any \EM{descending adjacent} pair, the left strip and the right strip must
	be of the \EM{same type} \BB, \BG, \GB\ or \GG.
\item If in an \EM{ascending adjacent} pair the left strip is of type \BB\ or \GB, then
	the right strip must be of the type \BB\ or \BG.
\item If in an \EM{ascending adjacent} pair the left strip is of type \BG\ or \GG, then
	the right strip must be of the type \GB\ or \GG.
\eit
\secC{Rows of columns of vertical strips}
\label{sec:rows-of-columns}
\def\dcovs{descending chains of vertical strips}
\def\asdc{ascending sequence of \dcovs}
\def\dcovs{column}
\def\asdc{row of \dcovs s}
Note that
\EM{if} every finite face of
$\overline{{\mathcal F}}$ is a vertical
strip, \EM{then} we can \EM{partition} the family of these faces in 
\EM{maximal chains of descending adjacent vertical strips}, i.e.,
\bit
\item sequences $\pas{v_1,v_2,\dots,v_l}$ of vertical strips \EM{of the same type}
\item where $\pas{v_{i},v_{i+1}}$ are \EM{descending adjacent} (as explained 
above) 
for $i=1,2,\dots,l-1$, 
\item which cannot be ``extended'' to a longer sequence with the same properties.
\eit
We call such maximal chains \EM{\dcovs s}.

Note that these \dcovs s may be ordered
such that (the vertical strips belonging to) two \EM{consecutive} \dcovs s are
\EM{ascending} adjacent: we call such succession of ascending
adjacent \dcovs s a \EM{\asdc}.

%
%
%

%

\secC{Free kinks and other configurations}
We consider the following four situations that might occur in the simple planar graph
$\overline{\mathcal F}$ corresponding to some folded overlay ${\mathcal F}$:
\bit
\item A \EM{free kink} is a \EM{unicoloured} path segment consisting of a horizontal step
	followed by a vertical step, whose ``middle vertex'' has degree $2$ (i.e., is \EM{not incident}
	to another edge) 
\item a \EM{horizontal bicoloured edge} (as defined above),
\item an \EM{isolated vertex} is a point in the essential region which does not belong to any path,
\item an \EM{uncrossed double step} is a \EM{unicoloured} path segment consisting of two horizontal steps.
\eit
The pictures below illustrate these simple notions:
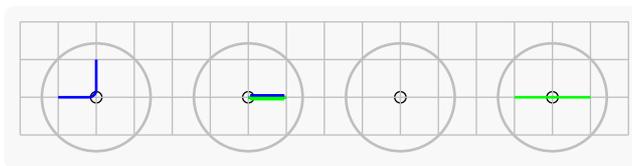
\cps{four_configurations}

\secC{Folded overlays down to some given level $y_0$}
Let  ${\overline{\mathcal F}}$ be the planar graph
corresponding to some folded overlay ${\mathcal F}$. For some fixed $z\in\Z$, we may \EM{remove} all vertices and
edges which lie
\EM{strictly below} the horizontal line $\pas{y=z}$ in both ${\overline{\mathcal F}}$: we denote by $\overline{\mathcal F_{z}}$ the result
of this operation.

This operation might introduce ``new loose ends'' (vertices of degree $1$) in ${\overline{\mathcal F_z}}$:
As before, we ``close'' such ``loose ends'' by adding edges of appropriate colour (again, these
artificial edges are only introduced ``for convenience'', in order to simplify the description).
Note that the boundaries of all finite faces in ${\overline{\mathcal F}_y}$ again
consist of a right segment and a left segment, as described in section~\ref{sec:faces}

\begin{lem}
\label{lem:no-free-kinks-no-configs}
Let  ${\mathcal F}$ be a folded overlay, and let $z\in\Z$. If $\overline{{\mathcal F}}_{z+1}$ does not contain
a \EM{free kink}, then $\overline{{\mathcal F}}_{z}$ does contain
\bit
\item neither an isolated vertex
\item nor a horizontal bicoloured edge
\item nor an uncrossed double step.
\eit
\end{lem}
\bpf
Consider the set 
$$
S=\setof{\zeta\in\Z\fdg {\overline{\mathcal F}}_\zeta\text{ does contain a horizontal step}}\subset\Z.
$$
If $S$ is empty, then $\overline{{\mathcal F}}$ clearly does contain neither a horizontal bicoloured
edge nor an uncrossed double step, and since in this case the whole graph appears as
a ``curtain of vertical paths, pending from the terminal points without any gaps'', there also is no isolated vertex:
So the assertion is true in this case.

If $S\neq\emptyset$, then consider $z_0=\max\of S\in\Z$.

For $z\geq z_0$,
$\overline{\mathcal F}_{z}$ might be
empty: in this case the assertion is trivially true. Otherwise, 
$\overline{\mathcal F}_{z}$ consists of a ``curtain of vertical paths''
(as above),
and some (\EM{unicoloured}!) horizontal steps
on level $z_0$ (only if $z=z_0$): so every vertex $\pas{x,y}$
in $\overline{{\mathcal F}}_{z}$ is \EM{adjacent} to a vertical step, and this implies that
$\overline{{\mathcal F}}_{z}$
does contain
\bit
\item neither an isolated vertex
\item nor a horizontal bicoloured edge $b$ (since the vertical \EM{upwards}
step incident with the left point of this edge $b$ would have the
\EM{same} colour as $b$) 
\item nor an uncrossed double step.
\eit

Now assume that the assertion is \EM{false}, and
let $z_1<z_0$ be the maximal level for which this is the case, i.e.,
$z_1$ is the \EM{maximal} integer such that
$\overline{{\mathcal F}}_{z_1+1}$ does not contain a free kink, but $\overline{{\mathcal F}}_{z_1}$ contains
\bit
\item either an isolated vertex
\item or a horizontal bicoloured edge
\item or an uncrossed double step.
\eit
But a close look at these situations reveals that this impossible:
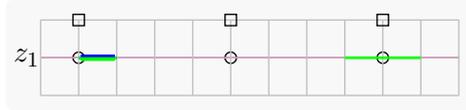
\cps{four_configurations2}
Note that in all three situations at level $z_1$ depicted above, the vertex
directly above (indicated by a small box in the pictures) \EM{cannot} be entered by a vertical
step from below. But it is not an isolated vertex (since there are none on
level $z_1+1$ by assumption), so
\bit
\item it must be entered from the left by a horizontal step \EM{of unique colour}, since
there are no bicoloured edges on level $z_1+1$,
\item and it must be left by a vertical step (of the same colour, of course),
since there are no uncrossed horizontal steps at level $z_1+1$.
\eit
But this would give a free kink on level $z_1+1$, a contradiction.
\epf

\begin{lem}
\label{lem:no-free-kinks-vertical-strips}
Let  ${\mathcal F}$ be a folded overlay, and let $z\in\Z$. If $\overline{{\mathcal F}}_{z+1}$
does not contain a free kink,
then every \EM{finite} face of $\overline{{\mathcal F}_{z}}$ is a \EM{vertical strip}.
\end{lem}
\bpf
Let $f$ be an arbitrary finite face of the planar simple graph $\overline{{\mathcal F}_{z}}$
and consider the left side $l$ of $f$, which we traverse now in the \EM{normal} direction
(i.e., rightwards and upwards).

Recall that $l$ starts with a \EM{vertical step}.
If $l$ does not contain a horizontal step, it must end with a
\EM{single} slanted step.

By Lemma~\ref{lem:no-free-kinks-no-configs}, $l$ does contain
\bit
\item neither a horizontal bicoloured edge
\item nor an uncrossed double step:
\eit
So every horizontal step $h$ belonging to $l$ 
is \EM{unicoloured}: if $h$ is not the \EM{last} step in $l$, then it must be
followed by a \EM{vertical edge} $v$ in $l$, and since $l$ starts with
a \EM{vertical} step, this horizontal step $h$ must appear on some level $>z$.

Since $\overline{{\mathcal F}}_{z+1}$ does not contain a free kink, $v$ must be \EM{bicoloured}. But the segment $\pas{h,v}$
cannot belong to the left side of $f$, as the following picture shows:
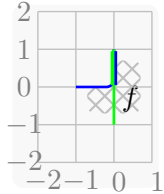
\cps{not_part_of_boundary}
So \EM{only the last step} in $l$ may be horizontal or slanted: in both cases, the right side of
$f$ is a path from the initial point $\pas{x,y}$ of $l$ to its terminal point $\pas{x+1,y+h}$ (for some
integer $h$) which starts with a horizontal or slanted step, and $f$ is a vertical strip.
\epf

\begin{lem}
If $\mathcal F$ is a folded \EM{survivor}, then $\overline{\mathcal F}$ cannot contain free kinks \EM{strictly above} the forbidden
line (i.e., the ``middle vertex'' of any free kink \EM{must lie on the forbidden line}).
\end{lem}

\bpf
Consider again the set 
$$
S=\setof{\zeta\in\Z\fdg \overline{\mathcal F}_\zeta\text{ does contain a horizontal step}}\subset\Z.
$$
If $S$ is empty, then $\overline{\mathcal F}$ clearly does not contain a free kink:
So the assertion is true in this case.

If $S\neq\emptyset$, then consider $z_0=\max\of S\in\Z$: clearly, $\overline{\mathcal F}_{z_0+1}$ does not
contain a free kink. Now assume that the assertion is false, and let $z\leq z_0$ be the maximal number
such that $\overline{{\mathcal F}_{z}}$ does contain a free kink \EM{above} the forbidden line.

By Lemma~\ref{lem:no-free-kinks-vertical-strips}, every finite face of
$\overline{{\mathcal F}_{z}}$ is a vertical
strip, and we may arrange the family of these faces in a \asdc.
Clearly, a free kink can only occur in \dcovs s of type \BG\ or \GB.
The following
picture shows the possibilities for such \dcovs s: 
\cps{bicoloured_chains_of_strips}
Picture $5$ (counted from the left) shows an impossible
chain (since blue paths never intersect the forbidden line) and is thus crossed out.

Picture $4$ shows a \dcovs\ of type \BG\ which ends with a (green) free kink on the forbidden line.

Picture $3$ shows a \dcovs\ of type \GB\ which does not end in a free kink.

Pictures $1$ and $2$ show free kinks strictly above the forbidden line
(marked by small rhombi in the pictures): it is easy to see that these chains would give
\EM{involutive bicoloured connections} of two terminal points of different colour,
which contradicts our assumption that $\mathcal F$ is a folded \EM{survivor}.
\epf

\begin{lem}
If $\mathcal F$ is a folded \EM{survivor} which contains a \dcovs\ $c$ of type \GB,
then there cannot be a \dcovs\  $c^\prime$ of type \BB\ which lies \EM{to the right of $c$}.

In particular, any \dcovs\  which is ascending adjacent ``to the right'' of $c$
\EM{must be of type \BG}.
\end{lem}
\bpf
If $c$ would end above the forbidden line, then it would have a free kink
above the forbidden line and thus yield an involutive bicoloured connection:
So $c$ ends on the forbidden line, and the same must be true for $c^\prime$:
But it is impossible for a \dcovs\  of type \BB\ to run into the forbidden line, see the middle picture in the following
illustration:
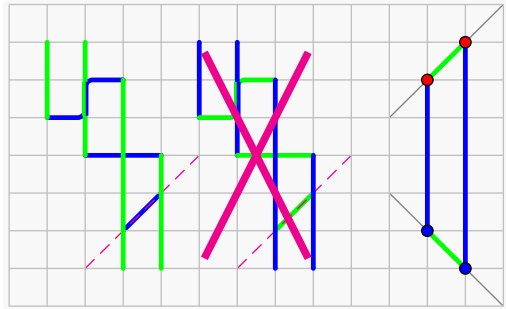
\cps{unicoloured_chains_of_strips}
(It is no problem for a \dcovs\  of type \GG\ to run into the forbidden line,
see the left picture above, and a \dcovs\  of type \BB\ must
end in blue initial points, see the right picture above.)

By the considerations in section~\ref{sec:adjacencies}, a \dcovs\  $c^\prime$ which is
ascending adjacent to the right of $c$ must
be of type \BG\ or \BB, and we just made clear that type \BB\ is impossible.
\epf

Our considerations make clear that certain distributions of colours for
the terminal points (which are not also initial points) are \EM{impossible}
for folded survivors:
\begin{cor}
\label{cor:colour-distribution}
If $\mathcal F$ is a folded \EM{survivor}, then the sequence of terminal points which are
not also initial points
\bit
\item may start with a sequence of consecutive blue points,
\item and every blue point not belonging to this starting sequence must be immediately followed
by a green point.
\eit

(In other words: after the first green terminal point, there can never follow a pair
of consecutive blue terminal points.)
\end{cor}

\secC{The mapping $\xi$ from the \lhs\ to the \rhs}
\label{sec:lhs-rhs}
We showed that the planar graph $\overline{\mathcal F}$
corresponding to a \EM{folded survivor} $\mathcal F$ necessarily appears
as a \asdc\ 
(where consecutive columns are \EM{ascending adjacent}):
It is easy to see that every \asdc\ 
\bit
\item corresponds bijectively to a planar graph $\overline{\mathcal F}$
	of some folded overlay ${\mathcal F}$
\item which does not have an involutive
bicoloured connection (in fact, bicoloured connections correspond
\EM{precisely} to \dcovs s of type \BG\ or \GB\ in a
folded survivor, and these
are \EM{not} involutive since they end on the forbidden line),
\eit
so we have the
following \EM{characterization}: 
a folded overlay $\mathcal F$ is a \EM{folded survivor} if and only if
its corresponding planar graph $\overline{\mathcal F}$ appears
as a \asdc.

The mapping $\xi$ from the \lhs\ to the \rhs\ now is best conceived by looking at
a sequence of pictures: the first picture shows our running example, where the maximal descending 
chains of vertical strips
are coloured:
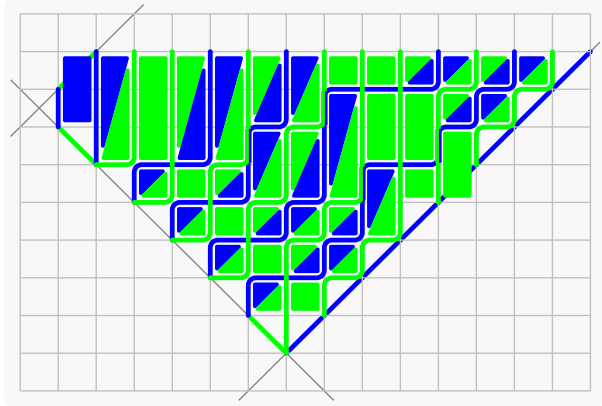
\cps{myex1_13_-9_19}
In the next picture, chains of type \GB\ with their mandatory ``right neighbour''
of type \BG\ are indicated as crosshatched:
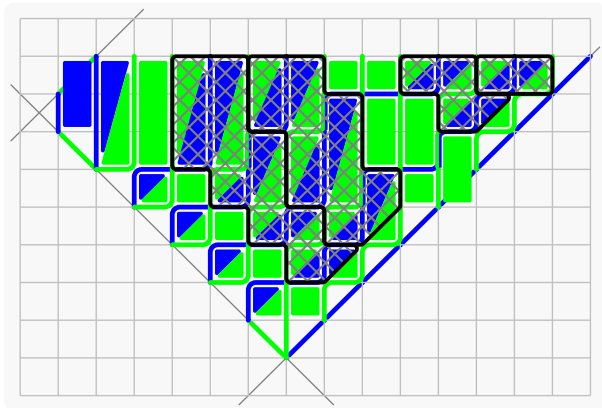
\cps{myex2_13_-9_19}
The next picture visualizes that the \asdc\
may be viewed as
(segments of) lattice paths (drawn in red):
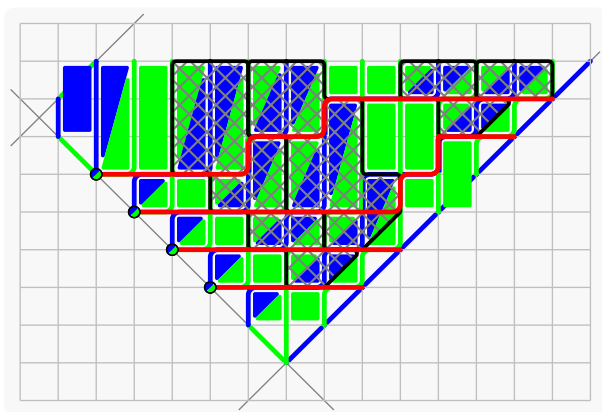
\cps{myex3_13_-9_19}
Note that these red paths
\bit
\item are nonintersecting,
\item and end on the forbidden line,
\item and can never intersect the diagonal $\pas{y=x}$:
\eit
So now
the diagonal plays the role of the forbidden line for these red path
segments, and the forbidden line plays the role of the diagonal (in the
sense that it contains
the terminal points for these paths). 

In the next picture, the artificial slanted edges are removed and
enforced segments of green paths
reaching the  forbidden line \EM{from below} are added:
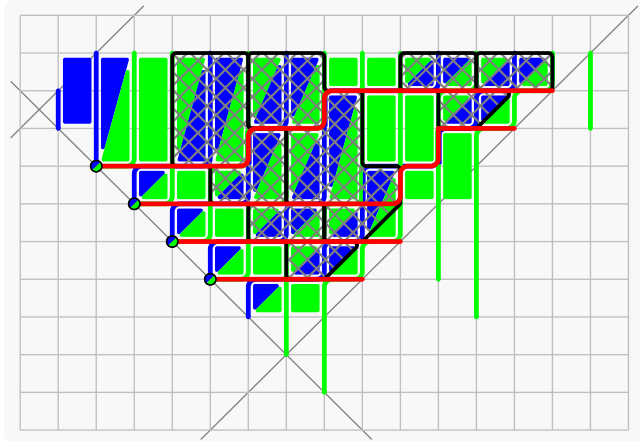
\cps{myex4_13_-9_19}
This  shows how the colour distribution of green or blue terminal points
for a folded survivor
\EM{determines} the green points on the forbidden line which must be reached 
\bit
\item from below (by the green segments of paths just added),
\item or from the left (by the red segments of paths).
\eit

Altogether, this shows that each \asdc\
which corresponds to
some folded survivor \EM{uniquely determines}
segments of \nilp s (the red paths in the above picture) connecting
	\bit
	\item certain 
	blue initial points
	\item with certain points on the forbidden line,
	\eit
(as illustrated in the above pictures).

In the following picture we added
\bit
\item the $01$--code for the \lhs\ (corresponding to the colours of the
blue and green terminal points)
\item and the \EM{enforced terminal segments}
	for the red and green paths ending on the forbidden line
	together with the corresponding \EM{reflected} $01$--code
	for these terminal segments:
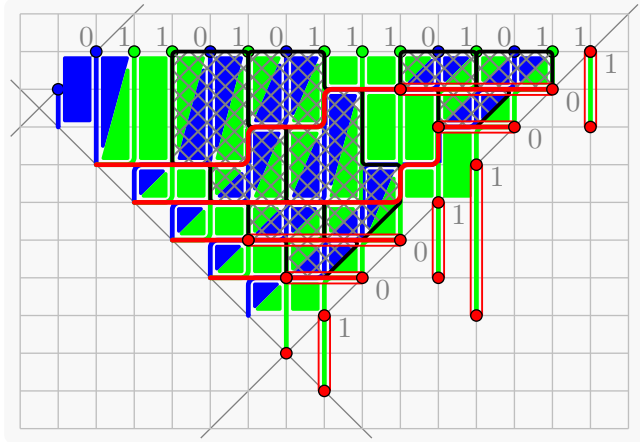
\eit	
\cps{myex5_13_-9_19}
This picture shows that
\bit
\item the (original) diagonal $\pas{y=x}$ plays the role of the forbidden line
\item and the (original) forbidden line plays the role of the diagonal
(in the sense that it contains the terminal points)
\eit
for the collection of red and green segments of paths.

\secC{The effect of $\xi$ on $01$--codes and their signs}
\label{sec:codes-and-signs}
The above picture also shows
\bit
\item the $01$--code $C$ for the \lhs\ (implied by the blue and green terminal
	points, in our example: $c = 01101011101011$)
\item \EM{and} the \EM{reflected} $01$--code $c$ corresponding to the
	black and red or green paths reaching the forbidden line, i.e.
	\bit
	\item $1$ for points reached \EM{from below},
	\item $0$ for points reached \EM{from the left}
	\eit
	(in our example: $c=10011001$).
\eit
(This \EM{swap} of zeros and ones corresponds to the fact that we consider
the \EM{reflected} situation for the \rhs\ here.)

The bijective correspondence between these two $01$--codes is easy to see:
For the $i$--th zero in the $01$--code $C$ corresponding to the green and blue terminal
points, we count the number $z_i$ of ones to
its left (in our example, these numbers are $\pas{0,2,3,6,7}$) and omit
all zeros among these numbers (in our example: $\pas{2,3,6,7}$). Then these
numbers are \EM{precisely} the \EM{indices} (starting with one) of the zeros
in the \EM{reflected} $01$--code $c$. But this
simple observation already proves that the mapping $\xi$ \EM{changes the sign}
of (the permutations corresponding to) these $01$--codes by a constant factor:
Recall (see section~\ref{sec:rhs-sign}) that the number of inversions
of the permutation
\bit
\item corresponding to some $01$--code $c$ for the
	\rhs\ (which is \EM{reflected} in our case, but we already accounted for
	the reflection by swapping ones and zeros) is equal
	to the sum of the indices of the zeros in $c$ (in our example:
	the permutation corresponding to $10011001$ has $2+3+6+7=18$ inversions),
	see \eqref{eq:sign-for-rhs},
\item corresponding to some $01$--code $C$ for the
	\lhs\ is
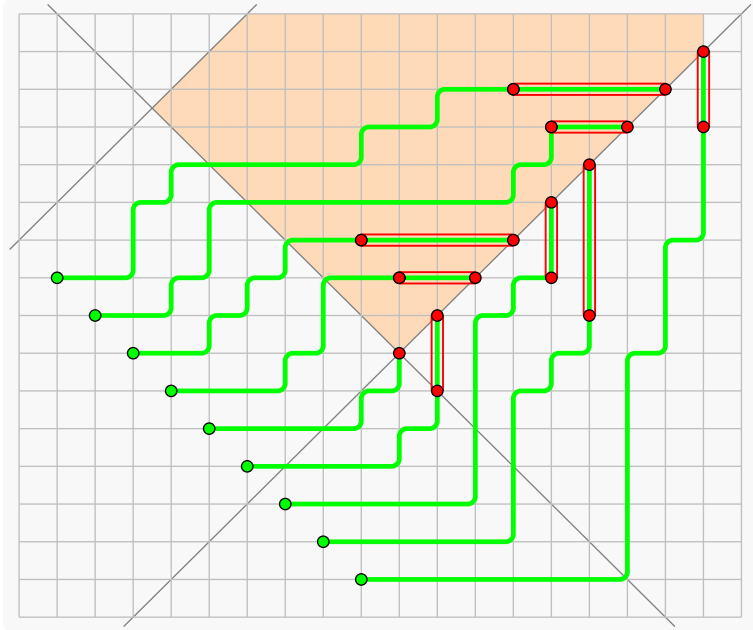
	$$
	\binom{m+k-1+\Iverson{\text{even}}}{2} + \inv\of{C},
	$$
	see \eqref{eq:sign-for-lhs}, and $\inv\of C$ (by definition) equals
	the sum of the numbers
	$\sum_i z_i$ defined above (in our example: $0+2+3+6+7=18$).
\eit

\secC{Completion of the mapping $\xi$}
Now by simply
\bit
\item omitting the folded survivor's overlay of green and blue paths
\item and changing the colours of the red paths to green
\eit
in the above picture,
we obtain the restriction to the essential region of some survivor from the
\EM{reflected} \rhs\ (in the picture, also segments of green paths \EM{outside} the essential
regions are shown, in order to indicate that we view the essential
region now as part of a bigger picture):
\cps{myex-rhs-0_13_-9_19}
If we \EM{restore} the segments of green paths \EM{outside} the essential regions (see the first picture of our running example), 
then we obtain the following ``normal'' survivor:
\cps{myex-rhs_13_-9_19}
It is easy to see that the above  picture corresponds to the situation of the
reflected \rhs\ with
$01$--code $01100110$, and
it is obvious that the mapping $\xi$ thus constructed is \EM{injective}.

\secC{Surjectivity of the mapping $\xi$}
To see
that $\xi$ is also \EM{surjective}, we shall exhibit its inverse
mapping: consider a ``normal'' survivor
with $01$--code $c$ from the \EM{reflected} \rhs, and 
\EM{interpret} the segments of its green paths inside the essential region
as \EM{red paths}
of a \EM{folded survivor} with code $C$ corresponding to $c$
(as described above). Then a \EM{unique} \asdc\ (and
thus a unique \EM{folded survivor})
is obtained ``step by step'', 
starting at the diagonal $\pas{y=x}$,
as illustrated in the following example:
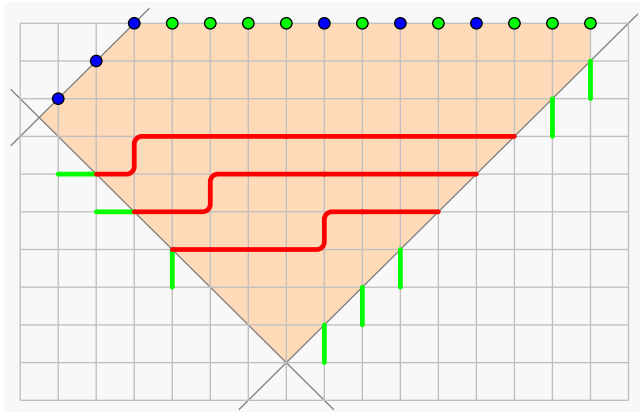
\cps{myex_reverse_0_13_-9_19}
The picture shows a section of
a (``normal'') survivor for the \EM{reflected} \rhs\ with (reflected)
$01$--code
$00011100$, where green path segments in the essential region are coloured red. According to the considerations in section~\ref{sec:codes-and-signs},
the colouring of terminal points for the corresponding folded survivor is
uniquely determined by this $01$--code 
(since the graphics is cut at
level $6$, the terminal blue and green points \EM{above} level $6$ are ``shifted down'').

If we draw the \asdc\ 
``pending from the blue and green
terminal points'' down to the the uppermost red path, then we obtain
the following picture:
\cps{myex_reverse_1_13_-9_19}

Now recall that in a folded survivor, every descending chain must contain
vertical strips \EM{of the same type}: so the second step of completing the
descending chains ``rightwards--downwards'' yields
the following situation:
\cps{myex_reverse_2_13_-9_19}
The third step yields:
\cps{myex_reverse_3_13_-9_19}
And the fourth step finishes the construction, and yields
a \asdc\ which we may view as the restriction to the essential region
of a \EM{folded survivor} :
\cps{myex_reverse_4_13_-9_19}

So the mapping $\xi$ is indeed a bijection which changes the sign of
survivors by the constant factor
$$
\pas{-1}^{\binom{m+k-1+\Iverson{\text{even}}}{2}}
$$
This finishes our proof of \eqref{eq:cigler-even} and \eqref{eq:cigler-odd}.


\end{document}